\documentstyle[11pt]{article}
\begin{document}

\newcommand{\qed}{\hphantom{.}\hfill $\Box$\medbreak}
\newcommand{\proof}{\noindent{\bf Proof \ }}
\newtheorem{Theorem}{Theorem}[section]
\newtheorem{Lemma}[Theorem]{Lemma}
\newtheorem{Corollary}[Theorem]{Corollary}
\newtheorem{Remark}[Theorem]{Remark}
\newtheorem{Example}[Theorem]{Example}
\newtheorem{Definition}[Theorem]{Definition}
\newtheorem{Construction}[Theorem]{Construction}

\thispagestyle{empty}

\begin{center}
{\Large\bf The fine triangle intersections for
maximum kite packings \footnote{This work was supported by the Fundamental Research Funds for the Central Universities under Grant $2011$JBZ$012$ (Y. Chang), $2011$JBM$298$(T. Feng), and the NSFC under Grant $61071221$ (Y. Chang), the NSFC under Grant $10901016$ (T. Feng).}}

\vskip24pt

Guizhi Zhang, Yanxun Chang, Tao Feng\\ Institute of Mathematics\\ Beijing Jiaotong
University\\ Beijing 100044, P. R. China\\ {\tt
yxchang@bjtu.edu.cn}\\ {\tt zgz\_hlbr@163.com}
 \\{\tt tfeng@bjtu.edu.cn}
\vskip12pt

\end{center}

\vskip12pt

\noindent {\bf Abstract:}  In this paper the fine triangle intersection problem for a pair of maximum kite packings is investigated. Let $Fin(v)=\{(s,t):$ $\exists$ a pair of
maximum kite packings of order $v$ intersecting in $s$ blocks and $s+t$
triangles$\}$. Let $Adm(v)=\{(s,t):\ s+t\leq b_v, s,t$ are non-negative integers$\}$,
where $b_v=\lfloor v(v-1)/8\rfloor$. It is established that $Fin(v)= Adm(v)\setminus \{(b_v-1,0),(b_v-1,1)\}$ for any integer $v\equiv 0,1\ ({\rm mod }\ 8)$ and $v\geq 8$; $Fin(v)=Adm(v)$ for any integer $v\equiv 2,3,4,5,6,7\ ({\rm mod }\ 8)$ and $v\geq 4$.



\section{Introduction}

Let $K_{v}$ be the complete graph with $v$ vertices and $\lambda K_{v}$
denote the graph $K_{v}$ with each of its edges replicated $\lambda$ times.
Given a family ${\cal G}$ of graphs each of which is simple and connected,
a \emph{$\lambda$-fold ${\cal G}$-design of order $v$}, denoted by
$(\lambda K_{v},{\cal G})$-design, is a pair $(X, {\cal B})$
where $X$ is the vertex set of $K_{v}$ and ${\cal B}$ is a collection
of subgraphs (called \emph{blocks}) of $\lambda K_{v}$,
such that each block is isomorphic to a graph in ${\cal G}$, and each edge
of $\lambda K_{v}$ belongs to exactly $\lambda$ blocks of ${\cal B}$.

If in the definition of $\cal G$-designs we replace the term ``exactly" with ``at most" (or ``at least"), we have a {\em $(\lambda K_{v},{\cal G})$-packing} (or {\em covering}).
When $\lambda =1$, a $(K_{v},{\cal G})$-packing (or covering) is called a ${\cal G}$-packing (or covering) of order $v$. When ${\cal G}$ contains a single graph $G$, i.e., ${\cal G}=\{G\}$, a $(\lambda K_{v},\{G\})$-design (packing or covering) is simply written as a $(\lambda K_{v},G)$-design (packing or covering). If $G$ is the
complete graph $K_k$, a $(K_{v},K_k)$-design is called a {\em
Steiner system} $S(2,k,v)$.

A $(K_{v},G)$-packing (or covering) $(X, {\cal B})$
is called \emph{maximum $($or minimum$)$} if there does not exist any
$(K_{v},G)$-packing (or covering) $(X, {\cal B'})$
with $|\cal B|<|\cal B'|$ (or $|\cal B|>|\cal B'|$).

Two $(K_{v},G)$-packings (or coverings) $(X, {\cal B}_{1})$ and $(X, {\cal B}_{2})$
are said to \emph{intersect} in $s$ blocks provided $|{\cal B}_{1}\cap {\cal B}_{2}|=s$.
If $s=0$, $(X,{\cal B}_{1} )$ and $(X, {\cal B}_{2})$ are said to be \emph{disjoint}.
The\emph{ intersection problem }for $(K_{v},G)$-packings (or coverings) is the determination of all integral
pairs $(v,s)$ such that there exists a pair of $(K_{v},G)$-packings (or coverings)
intersecting in $s$ blocks.

The intersection problem for $S(2,k,v)$s was first introduced by Kramer
and Mesner in \cite{km}. A complete solution to the intersection
problem for $S(2,3,v)$s was made by Lindner and Rosa \cite{lr}. The
intersection problem for $S(2,4,v)$s was dealt with by Colbourn et
al. \cite{chl}, apart from three undecided values for $v=25$, $28$
and $37$. Billington and Kreher \cite{bk} completed the intersection
problem for all connected simple graphs $G$ where the minimum of the
number of vertices and the number of edges of $G$ is not bigger than
$4$. The intersection problem is also considered for many other
types of combinatorial structures. The interested reader may refer
to \cite{b,bgl,bh,cl,cl1,fu,gl,hl}.

Let $B$ be a simple graph. Denote by $T(B)$ the set of all triangles
of the graph $B$. For example, if $B$ is the graph with vertices
$a$, $b$, $c$, $d$ and edges $ab$, $ac$, $bc$, $cd$, (such a
graph is called a {\em kite}), then
$T(B)=\{\{a,b,c\}\}$. Two $(K_{v},G)$-packings (or coverings)
$(X,{\cal B}_1)$ and $(X,{\cal B}_2)$ are said to {\em intersect} in $t$
triangles provided $|T({\cal B}_1)\cap T({\cal B}_2)|=t$, where
$T({\cal B}_i)=\bigcup_{B\in {\cal B}_i} T(B)$, $i=1$, $2$. The {\em
triangle intersection problem} for $(K_{v},G)$-packings (or coverings) is the determination
of all integer pairs $(v,t)$ such that there exists a pair of
$(K_{v},G)$-packings (or coverings) intersecting in $t$ triangles.

The triangle intersection problem was first considered by Lindner and Yazici in
\cite{ly}, who made a complete solution to the triangle intersection
problem for kite systems. Billington et al. \cite{bly} solved the
triangle intersection problem for $(K_4-e)$-designs. Chang et al.
\cite{cfl} investigated the triangle intersection problem for
$S(2,4,v)$s.

Every block $B$ in a $(K_{v},G)$-packing (or covering) contributes $|T(B)|=|T(G)|$ triangles.
If two $(K_{v},G)$-packings (or coverings) $(X, {\cal B}_{1})$ and $(X, {\cal B}_{2})$ intersect in $s$ blocks, then they intersect in
at least $s|T(G)|$ triangles. It is natural to ask how about the triangle intersection problem for a pair
of $(K_{v},G)$-packings (or coverings) intersecting in $s$ blocks. Thus the fine triangle intersection problem was introduced in \cite{cflt}.
Define $Fin_{G}(v)=\{(s,t):\exists$ a pair of $(K_{v},G)$-packings (or coverings) intersecting in
$s$ blocks and $t+s|T(G)|$ triangles\}. The {\em fine triangle intersection problem} for $(K_{v},G)$-packings (or coverings) is
to determine $Fin_{G}(v)$.

Chang et al. has completely solved the fine triangle intersection problems for kite systems \cite{cflt} and $(K_4-e)$-designs \cite{cflt1,cflt2}. The purpose of this paper is to example the fine triangle intersection problem for maximum kite packings. In what
follows we always write $Fin_{G}(v)$ simply as $Fin(v)$ when $G$ is a kite, that is,
$Fin(v)=\{(s,t):\exists$ a pair of maximum kite packings of order $v$ intersecting in
$s$ blocks and $t+s$ triangles\}. It is known that for any positive integer $v$ there is a maximum kite packings of order $v$ with $\lfloor v(v-1)/8\rfloor$ blocks \cite{cl2}.

Let $Adm(v)=\{(s,t): s+t\leq b_{v}, s,t$ are non-negative integers$\}$, where $b_{v}=\lfloor v(v-1)/8\rfloor$. Since every block in a maximum kite packing contributes only one triangle, we have that $Fin(v)\subseteq Adm(v)$. In the following we always denote the copy of the kite with vertices $a,b,c,d$ and edges $ab,ac,bc,cd$ by $[a,b,c-d]$.

\begin{Example}\label{Fin-4}
$Fin(4)= Adm(4)$.
\end{Example}

\proof Take the vertex set $X=\{0,1,2,3\}$. Let ${\cal B}=\{[0,1,3-2]\}$.
 Then $(X,{\cal B})$ is a maximum kite packing of order $4$. Consider the following
permutations on $X$.
\begin{center}\small
\begin{tabular}{lll}
$\pi_{0,0}=(2\ 3)$,&
$\pi_{0,1}=(0\ 1\ 3)$,&
$\pi_{1,0}=(1)$.
\end{tabular}
\end{center}
\noindent We have that for each $(s,t)\in Adm(4)$,  $|\pi_{s,t}{\cal B}\cap {\cal B}|=s$ and $|T(\pi_{s,t}{\cal B}\setminus{\cal B})\cap T({\cal
B}\setminus\pi_{s,t}{\cal B})|=t$. \qed

As the main result of the present paper, we are to prove the following theorem.

\begin{Theorem}\label{main theorem}
$Fin(v)= Adm(v)\setminus \{(b_v-1,0),(b_v-1,1)\}$ for any integer $v\equiv 0,1\ ({\rm mod }\ 8)$ and $v\geq 8$.
$Fin(v)=Adm(v)$ for any integer $v\equiv 2,3,4,5,6,7\ ({\rm mod }\ 8)$ and $v\geq 4$.
\end{Theorem}

\section{Basic design constructions}

Let $K$ be a set of positive integers. A {\em group divisible
design} (GDD) $K$-GDD is a triple ($X, {\cal G},{\cal A}$)
satisfying the following properties: ($1$) $\cal G$ is a partition
of a finite set $X$ into subsets (called {\em groups}); ($2$) $\cal
A$ is a set of subsets of $X$ (called {\em blocks}), each of
cardinality from $K$, such that every $2$-subset of $X$ is either
contained in exactly one block or in exactly one group, but not in
both. If $\cal G$ contains $u_i$ groups of size $g_i$ for $1\leq
i\leq r$, then we call $g_1^{u_1}g_2^{u_2}\cdots g_r^{u_r}$ the {\em
group type} (or {\em type}) of the GDD. If $K=\{k\}$, we write
a $\{k\}$-GDD as a $k$-GDD.

  Let ${\cal H}=\{H_1, H_2, \ldots, H_m\}$ be a partition of a finite
set $X$ into subsets (called {\em holes}), where $|H_i|=n_i$ for
$1\leq i\leq m$. Let $K_{n_1, n_2, \ldots, n_m}$ be the complete
multipartite graph on $X$ with the $i$-th part on $H_i$, and $G$ be
a subgraph of $K_{n_1, n_2, \ldots, n_m}$. A {\em holey $G$-design} is a triple $(X,{\cal
H},{\cal B})$ such that $(X,{\cal B})$ is a
$(K_{n_1,n_2,\ldots,n_m},G)$-design. The {\em hole type} (or {\em
type}) of the holey $G$-design is $\{n_1, n_2, \ldots, n_m\}$. We use an ``exponential"
notation to describe hole types: the hole
type $g_1^{u_1}g_2^{u_2}\cdots g_r^{u_r}$ denotes $u_i$ occurrences
of $g_i$ for $1\leq i\leq r$. Obviously if $G$ is the complete graph
$K_k$, a holey $K_k$-design is just a $k$-GDD. When $G$ is kite,
a holey $G$-design is said to be a \emph{kite-GDD}.

A pair of holey $G$-designs $(X,{\cal H},{\cal B}_1)$ and $(X,{\cal
H},{\cal B}_2)$ of the same type is said to {\em intersect in $s$
blocks} if $|{\cal B}_1\cap {\cal B}_2|=s$. A pair of holey
$G$-designs $(X,{\cal H},{\cal B}_1)$ and $(X,{\cal H},{\cal B}_2)$
of the same type is said to {\em intersect in $t$ triangles} if
$|T({\cal B}_1)\cap T({\cal B}_2)|=t$, where $T({\cal
B}_i)=\bigcup_{B\in {\cal B}_i} T(B)$, $i=1,2$.

Wilson's fundamental construction on GDDs \cite{wilson} has been adapted to obtain the following construction on the intersection numbers of holey $G$-designs.

\begin{Construction}
\label{WeightingConstruction}  {\rm (\cite{cflt})}{\rm (Weighting Construction)}  Suppose
that $(X, {\cal G}, {\cal A})$ is a $K$-GDD, and let $\omega:\,
X\longmapsto Z^{+}\cup \{ 0 \}$ be a weight function. For every
block $A\in {\cal A}$, suppose that there is a pair of holey
$G$-designs of type $\{\omega(x) : \ x\in A\}$, which intersect in
$b_A$ blocks and $t_A$ triangles. Then there exists a pair of holey
$G$-designs of type $\{\sum_{x\in H}\omega(x): \ H\in {\cal G} \}$,
which intersect in $\sum_{A\in {\cal A}}b_A$ blocks and $\sum_{A\in
{\cal A}}t_A$ triangles.
\end{Construction}

The following construction is simple but very useful, which is a variation of Construction $2.2$ in \cite{cflt}.

\begin{Construction}
\label{FillingSubdesigns} {\rm (Filling Construction)} Let $a$ be a
nonnegative integer. Suppose that there exists a pair of holey
$G$-designs of type $\{g_1,g_2,\ldots,g_s\}$, which intersect in $b$
blocks and $t$ triangles. If there is a pair of
$((K_{g_i+a}\setminus K_a),G)$-designs with the same subgraph $K_a$
removed for each $1\leq i\leq s-1$, which intersect in $b_i$ blocks
and $t_i$ triangles, and there is a pair of $(K_{g_s+a},G)$-packings $($or coverings$)$,
which intersect in $b_s$ blocks and $t_s$ triangles, then there
exists a pair of $(K_{v+a},G)$-packings $($or coverings$)$ intersecting in
$b+\sum_{i=1}^s b_i$ blocks and $t+\sum_{i=1}^s t_i$ triangles.
\end{Construction}

We quote the following result for later use.

\begin{Lemma}
\label{3-GDD}{\rm (\cite{chr})} Let $g$, $t$ and $u$ be nonnegative
integers. There exists a $3$-$GDD$ of type $g^{t}u^{1}$ if and only
if the following conditions are all satisfied:

{\em (1)}\ if\ $g>0$, then $t\geq 3$, or  $t=2$ and $u=g$, or $t=1$
and $u=0$, or $t=0$;

{\em (2)}\ $u\leq g(t-1)$ or $gt=0$;

{\em (3)}\ $g(t-1)+u\equiv 0\ (\rm{mod}\ 2)$ or $gt=0$;

{\em (4)}\ $gt\equiv 0\ (\rm{mod}\ 2)$ or $u=0$;

{\em (5)}\ $g^{2}t(t-1)/2+gtu\equiv 0\ (\rm{mod}\ 3)$.
\end{Lemma}

\section{Fin$(v)$ for $5\leq v\leq 23$ and $v\neq 8,9,16,17$}

When $v\in\{8,9,16,17\}$, it is shown that $Fin(v)= Adm(v)\setminus\{(b_v-1,0),(b_v-1,1)\}$ in \cite{cflt}, where $b_{v}=\lfloor v(v-1)/8\rfloor$. In this section, we shall examine $Fin(v)$ for $5\leq v\leq 23$ and $v\neq 8,9,16,17$, by {\em ad hoc} methods.

\begin{Lemma}\label{Fin-5}
$Fin(5)= Adm(5)$.
\end{Lemma}

\proof Take the vertex set $X=\{0,1,2,3,4\}$. Let ${\cal B}=\{[0,4,3-1],[0,1,2-3]\}$.
 Then $(X,{\cal B})$ is a maximum kite packing of order $5$. Consider the following
permutations on $X$.
\begin{center}\small
\begin{tabular}{lll}
$\pi_{0,0}=(0\ 4)(2\ 3)$,&
$\pi_{0,1}=(0\ 1)(3\ 4)$,&
$\pi_{0,2}=(1\ 4)(2\ 3)$,\\
$\pi_{1,0}=(0\ 1\ 3\ 2\ 4)$,&
$\pi_{1,1}=(1\ 3\ 2\ 4)$,&
$\pi_{2,0}=(1)$.
\end{tabular}
\end{center}
\noindent We have that for each $(s,t)\in Adm(5)$,  $|\pi_{s,t}{\cal B}\cap {\cal B}|=s$ and $|T(\pi_{s,t}{\cal B}\setminus{\cal B})\cap T({\cal
B}\setminus\pi_{s,t}{\cal B})|=t$. \qed

\begin{Lemma}\label{Fin-6}
$Fin(6)= Adm(6)$.
\end{Lemma}

\proof Take the vertex set $X=\{0,1,2,3,4,5\}$. Let ${\cal B}_1=\{[0,1,2-3],[0,3,4-5],[1,3,5-0]\}$ and
${\cal B}_2=\{[0,1,2-3],[0,3,5-1],[2,5,4-3]\}$. Then $(X,{\cal B}_i)$ is a maximum kite packing of order $6$ for each $i=1,2$. Consider the following permutations on $X$.
\begin{center}\small
\begin{tabular}{lllll}
$\pi_{0,0}=(0\ 4)(2\ 3)$,&
$\pi_{0,1}=(0\ 4)(1\ 2)$,&
$\pi_{0,2}=(0\ 4)(1\ 5)$,&
$\pi_{0,3}=(0\ 2)(3\ 4)$,&
$\pi_{1,0}=(1\ 3)$,\\
$\pi_{1,1}=(0\ 5\ 4\ 2\ 1\ 3)$,&
$\pi_{1,2}=(0\ 1)(4\ 5)$,&
$\pi_{2,0}=(0\ 1)$,&
$\pi_{2,1}=(0\ 3)(2\ 5)$,&
$\pi_{3,0}=(1)$.
\end{tabular}
\end{center}
\noindent Let $M_1=Adm(6)\setminus\{(0,3),(2,0)\}$ and $M_2=\{(0,3),(2,0)\}$. We have that for each $(s,t)\in M_i$, $i=1,2$, $|\pi_{s,t}{\cal B}_i\cap {\cal B}_i|=s$ and
 $|T(\pi_{s,t}{\cal B}_i\setminus{\cal B}_i)\cap T({\cal B}_i\setminus\pi_{s,t}{\cal B}_i)|=t$. \qed

\begin{Lemma}\label{Fin-7}
$Fin(7)= Adm(7)$.
\end{Lemma}
\proof Take the vertex set $X=\{0,1,2,3,4,5,6\}$. Let ${\cal B}_1=\{[0,2,1-6],[2,4,3-1],[1,4,5-0],[2,6,5-3],[6,4,0-3]\}$, ${\cal B}_2=({\cal B}_1\setminus\{[6,4,0-3]\})\cup\{[6,3,0-4]\}$ and ${\cal B}_3=({\cal B}_1\setminus\{[2,6,5-3]\}\cup\{[2,5,6-3]\}$. Then $(X,{\cal B}_i)$ is a maximum kite packing of order $7$ for each $i=1,2,3$. Consider the following
permutations on $X$.
\begin{center}\small
\begin{tabular}{llll}
$\pi_{0,0}=(0\ 4)(2\ 6\ 3)$,&
$\pi_{0,1}=(0\ 6\ 5\ 3\ 1\ 4)$,&
$\pi_{0,2}=(0\ 6\ 5\ 2)(1\ 4\ 3)$,&
$\pi_{0,3}=(0\ 1\ 2)(3\ 4\ 6)$,\\
$\pi_{0,4}=(3\ 5)(4\ 6)$,&
$\pi_{0,5}=(1\ 6)(2\ 4)$,&
$\pi_{1,0}=(0\ 2\ 4\ 1)(3\ 5\ 6)$,&
$\pi_{1,1}=(0\ 1\ 3\ 5)(2\ 6\ 4)$,\\
$\pi_{1,2}=(0\ 6\ 5\ 1)(2\ 3\ 4)$,&
$\pi_{1,3}=(0\ 2\ 4)(1\ 3\ 6)$,&
$\pi_{1,4}=(0\ 6\ 5\ 1)(2\ 4)$,&
$\pi_{2,0}=(0\ 6\ 3\ 5\ 1)(2\ 4)$,\\
$\pi_{2,1}=(0\ 5\ 3\ 1\ 6\ 4)$,&
$\pi_{2,2}=(0\ 1\ 4\ 2\ 5\ 3\ 6)$,&
$\pi_{2,3}=(0\ 5)(2\ 4)$,&
$\pi_{3,0}=(0\ 3\ 1\ 6\ 4\ 2)$,\\
$\pi_{3,1}=(0\ 5\ 3\ 1\ 6\ 4)$,&
$\pi_{3,2}=(0\ 5)(2\ 4)$,&
$\pi_{4,0}=\pi_{4,1}=\pi_{5,0}=(1)$.
\end{tabular}
\end{center}
\noindent We have that for each $(s,t)\in Adm(7)\setminus\{(4,0),(3,1),(2,3),(4,1)\}$,
 $|\pi_{s,t}{\cal B}_1\cap {\cal B}_1|=s$ and
 $|T(\pi_{s,t}{\cal B}_1\setminus{\cal B}_1)\cap T({\cal B}_1\setminus\pi_{s,t}{\cal B}_1)|=t$.
For each $(s,t)\in\{(4,0),(3,1)\}$,
 $|\pi_{s,t}{\cal B}_2\cap {\cal B}_1|=s$ and
 $|T(\pi_{s,t}{\cal B}_2\setminus{\cal B}_1)\cap T({\cal B}_1\setminus\pi_{s,t}{\cal B}_2)|=t$.
For each $(s,t)\in\{(2,3),(4,1)\}$,
 $|\pi_{s,t}{\cal B}_3\cap {\cal B}_1|=s$ and
 $|T(\pi_{s,t}{\cal B}_3\setminus{\cal B}_1)\cap T({\cal B}_1\setminus\pi_{s,t}{\cal B}_3)|=t$.\qed

\begin{Lemma}\label{10}
$Fin(10)= Adm(10)$.
\end{Lemma}

\proof Take the vertex set $X=\{0,1,\ldots,9\}$. Let
\begin{center}
\small
 \begin{tabular}{llllllll}
${\cal B}_1:$ &$[2,0,1-9]$,&$[4,7,0-5]$,&$[1,4,8-2]$,&$[1,7,5-8]$,&$[2,4,3-1]$,&$[7,8,6-5]$,&
$[2,9,5-3]$,\\\vspace{0.1cm}&$[9,0,8-3]$,&$[0,3,6-1]$,&$[6,9,4-5]$,&$[9,3,7-2]$;\\

${\cal B}_2:$ &$[2,0,1-9]$,&$[4,7,0-5]$,&$[1,4,8-2]$,&$[1,7,5-8]$,&$[2,4,3-1]$,&$[7,8,3-9]$,&
$[2,9,5-3]$,\\\vspace{0.1cm}&$[8,6,9-0]$,&$[0,3,6-1]$,&$[5,6,4-9]$,&$[6,2,7-9]$;\\

${\cal B}_3:$ &$[7,9,5-8]$,&$[1,4,8-2]$,&$[2,5,0-9]$,&$[1,5,3-0]$,&$[2,9,4-7]$,&$[1,6,9-3]$,&
$[8,6,0-4]$,\\\vspace{0.1cm}&$[7,3,8-9]$,&$[0,7,1-2]$,&$[2,7,6-5]$,&$[6,4,3-2]$;\\

${\cal B}_4:$ &$[1,2,0-3]$,&$[3,2,4-5]$,&$[4,0,7-8]$,&$[1,4,8-2]$,&$[1,7,5-8]$,&$[9,0,8-3]$,&
$[2,5,9-1]$,\\\vspace{0.1cm}&$[4,9,6-3]$,&$[7,2,6-1]$,&$[9,7,3-1]$,&$[0,5,6-8]$;\\

${\cal B}_{5}:$ &$[0,1,2-7]$,&$[4,7,0-5]$,&$[1,4,8-3]$,&$[1,7,5-3]$,&$[4,3,2-6]$,&$[7,8,6-1]$,&
$[2,9,5-8]$,\\\vspace{0.1cm}&$[9,0,8-2]$,&$[0,3,6-5]$,&$[6,9,4-5]$,&$[3,7,9-1]$;\\

${\cal B}_{6}:$ &$[2,0,1-9]$,&$[4,7,0-5]$,&$[1,4,8-3]$,&$[1,7,5-3]$,&$[2,4,3-1]$,&$[7,8,6-1]$,&
$[2,9,5-8]$,\\\vspace{0.1cm}&$[9,0,8-2]$,&$[0,3,6-5]$,&$[6,9,4-5]$,&$[9,3,7-2]$;\\

${\cal B}_{7}:$ &$[2,0,1-3]$,&$[4,0,7-2]$,&$[1,4,8-3]$,&$[1,7,5-3]$,&$[3,4,2-6]$,&$[7,8,6-1]$,&
$[2,9,5-8]$,\\\vspace{0.1cm}&$[9,0,8-2]$,&$[0,3,6-5]$,&$[6,9,4-5]$,&$[7,3,9-1]$;\\

${\cal B}_{8}:$ &$[7,9,5-8]$,&$[1,4,8-2]$,&$[2,5,0-7]$,&$[1,5,3-0]$,&$[2,9,4-7]$,&$[1,6,9-3]$,&
$[8,6,0-4]$,\\\vspace{0.1cm}&$[7,3,8-9]$,&$[7,2,1-0]$,&$[2,3,6-7]$,&$[6,5,4-3]$.\\
\end{tabular}
\end{center}

\noindent Let ${\cal B}_{9}=({\cal B}_1\setminus\{[1,4,8-2],[9,0,8-3]\})\cup\{[1,4,8-3],[9,0,8-2]\}$, ${\cal B}_{10}=({\cal B}_{9}\setminus\{[0,3,6-1]\})\cup\{[0,3,6-2]\}$, ${\cal B}_{11}=({\cal B}_{9}\setminus\{[1,7,5-8],[2,9,5-3]\})\cup\{[1,7,5-3],[2,9,5-8]\}$; ${\cal B}_{12}=({\cal B}_{1}\setminus\{[1,4,8-2],[7,8,6-5],[6,9,4-5],[9,3,7-2]\})\cup\{[1,4,8-6],[7,8,2-6],[4,5,6-9],
[3,7,9-4]\}$, ${\cal B}_{13}=({\cal B}_{12}\setminus\{[7,8,2-6]\})\cup\{[7,6,2-8]\}$, ${\cal B}_{14}=({\cal B}_{13}\setminus\{[1,4,8-6]\})\cup\{[1,4,8-7]\}$, ${\cal B}_{15}=({\cal B}_{14}\setminus\{[9,0,8-3]\})\cup\{[9,0,8-6]\}$, ${\cal B}_{16}=({\cal B}_{15}\setminus\{[1,4,8-7],[3,7,9-4]\})\cup\{[1,4,8-3],[3,9,7-8]\}$, ${\cal B}_{17}=({\cal B}_{16}\setminus\{[1,7,5-8],[2,9,5-3],[3,9,7-8]\})\cup\{[1,7,5-3],[2,9,5-8],[3,7,9-4]\}$, ${\cal B}_{18}=({\cal B}_{17}\setminus\{[3,7,9-4]\})\cup\{[3,9,7-8]\}$; ${\cal B}_{19}=({\cal B}_3\setminus\{[0,7,1-2],[2,7,6-5],[6,4,3-2]\})\cup\{[7,2,1-0],[2,3,6-7],[6,5,4-3]\}$; ${\cal B}_{20}=({\cal B}_4\setminus\{[1,2,0-3],[3,2,4-5]\})\cup\{[3,2,0-1],[3,5,4-2]\}$, ${\cal B}_{21}=({\cal B}_4\setminus\{[1,2,0-3],[3,2,4-5],[1,4,8-2]\})\cup\{[3,2,0-1],[3,5,4-2],[8,4,1-2]\}$.
Then $(X,{\cal B}_i)$ is a maximum kite packing of order $10$ for each $1\leq i\leq 21$. Consider the following permutations on $X$.
\begin{center}\small
\begin{tabular}{lll}
$\pi_{0,0}=(0\ 8\ 4)(2\ 9\ 6\ 7\ 3)$,&
$\pi_{0,1}=(0\ 4\ 1\ 3\ 2)(6\ 9\ 8\ 7)$,&
$\pi_{0,2}=(0\ 7\ 1)(2\ 4)(3\ 5\ 9\ 6\ 8)$,\\
$\pi_{0,3}=(0\ 5\ 4\ 1)(2\ 6\ 9\ 3)$,&
$\pi_{0,4}=(0\ 7\ 1\ 6\ 4\ 5\ 9\ 2\ 3)$,&
$\pi_{0,5}=(0\ 4\ 3\ 9)(2\ 8\ 7)$,\\
$\pi_{0,6}=(0\ 8\ 7)(1\ 3)(2\ 5)(4\ 9\ 6)$,&
$\pi_{0,7}=(1\ 4\ 9\ 3\ 5)(6\ 7\ 8)$,&
$\pi_{0,8}=(0\ 9\ 7)(1\ 2)(3\ 4\ 8)$,\\
$\pi_{0,9}=(0\ 7)(2\ 5)(3\ 8)$,&
$\pi_{0,10}=(0\ 9)(1\ 5)(3\ 4)$,&
$\pi_{0,11}=(0\ 6)(1\ 5)(2\ 4)(8\ 9)$,\\
$\pi_{1,0}=(0\ 3\ 2\ 4\ 6\ 7\ 9\ 1\ 5\ 8)$,&
$\pi_{1,1}=(0\ 9\ 7\ 5\ 4)(1\ 3)(2\ 8\ 6)$,&
$\pi_{1,2}=(1\ 3)(2\ 6\ 8)$,\\
$\pi_{1,3}=(0\ 6\ 1\ 8\ 4)(2\ 7)(3\ 5)$,&
$\pi_{1,4}=(0\ 7\ 3\ 8\ 2)(1\ 5)(4\ 9)$,&
$\pi_{1,5}=(0\ 9\ 2\ 3\ 1\ 7)(4\ 8\ 5)$,\\
$\pi_{1,6}=(0\ 8\ 9)(1\ 4\ 5\ 2)(3\ 6\ 7)$,&
$\pi_{1,7}=(0\ 4\ 9\ 7\ 6\ 3\ 2\ 1\ 8)$,&
$\pi_{1,8}=(0\ 6\ 7\ 9)(1\ 5\ 2)(3\ 8)$,\\
$\pi_{1,9}=(0\ 1)(3\ 5)(4\ 9)(6\ 7)$,&
$\pi_{1,10}=(0\ 4)(1\ 3)(5\ 9)(6\ 8)$,&
$\pi_{2,0}=(0\ 5\ 6\ 7)$,\\
$\pi_{2,1}=(0\ 7\ 9\ 4\ 3\ 6\ 5\ 2)$,&
$\pi_{2,2}=(0\ 4\ 1)(6\ 8)$,&
$\pi_{2,3}=(0\ 6)(2\ 3)(5\ 7)$,\\
$\pi_{2,4}=(1\ 3\ 5)(4\ 9)(7\ 8)$,&
$\pi_{2,5}=(0\ 9\ 3\ 2\ 5\ 1\ 8\ 7\ 4\ 6)$,&
$\pi_{2,6}=(1\ 5)$,\\
$\pi_{2,7}=(0\ 7\ 9)(1\ 5\ 2)(3\ 8\ 4)$,&
$\pi_{2,8}=(0\ 3\ 9)(1\ 4\ 5)(6\ 7\ 8)$,&
$\pi_{3,0}=(4\ 6\ 7)$,\\
$\pi_{3,1}=(5\ 8\ 7)$,&
$\pi_{3,2}=(0\ 6\ 8\ 9\ 7)(1\ 2\ 5)(3\ 4)$,&
$\pi_{3,3}=(0\ 3\ 8\ 9\ 7\ 6\ 4)(1\ 2)$,\\
$\pi_{3,4}=(5\ 9)$,&
$\pi_{3,5}=(0\ 9)(6\ 7)$,&
$\pi_{3,6}=(0\ 6\ 7\ 9)(1\ 2)(3\ 8)$,\\
$\pi_{3,7}=(0\ 5)(1\ 9)(3\ 4)(7\ 8)$,&
$\pi_{4,0}=(0\ 2\ 1\ 3\ 8\ 6\ 4)(7\ 9)$,&
$\pi_{4,1}=(1\ 6)$,\\
$\pi_{4,2}=(0\ 1)$,&
$\pi_{4,3}=(1\ 2)$,&
$\pi_{4,4}=(0\ 2)(3\ 7)$,\\
$\pi_{4,5}=(1\ 2)(3\ 8)(6\ 9)$,&
$\pi_{4,6}=(0\ 9)(1\ 5)(3\ 4)$,&
$\pi_{5,0}=(0\ 9)$,\\
$\pi_{5,1}=(2\ 4)$,&
$\pi_{5,2}=(6\ 8)$,&
$\pi_{5,3}=(0\ 2)(4\ 6)(5\ 8)$,\\
$\pi_{6,0}=(7\ 8)$,&
$\pi_{6,1}=(6\ 8)$.\\
\end{tabular}
\end{center}

\noindent Let $S_1=\{(x,y): x+y\leq 11, x\in \{6,7,8,9,10,11\}, y$ is a non-negative integer$\}\setminus \{(6,0),(6,1)\}$, and $S_2=\{(2,9),(3,8),(4,7),(5,4),(5,5),(5,6)\}$. For each $(x',y')\in S_1\cup S_2$, take the identity permutation $\pi_{x'y'}=(1)$. It is readily checked that for each $i,j,s,t$ in Table $1$,
 $|{\cal B}_i\cap \pi_{s,t}{\cal B}_j|=s$ and
 $|T({\cal B}_i\setminus\pi_{s,t}{\cal B}_j)\cap T(\pi_{s,t}{\cal B}_j\setminus{\cal B}_i)|=t$. \qed

\begin{center}
Table $1$.  \ Fine triangle intersections for maximum kite packings of order $10$
  \begin{tabular}{ccc|ccc}\hline
$i$ &  $j$ &  $(s,t)$ & $i$ &  $j$ &  $(s,t)$ \\\hline
1  &   1   &  (0,0),(0,1),$\ldots$,(0,10), & 1  &   2   &  (7,0)\\
&& (1,0), (1,1),$\ldots$,(1,8),&  1  &   5    &   (2,9)   \\
&& (2,0), (2,1),$\ldots$,(2,7),&  1  &   6    &   (5,6)   \\
&& (3,0), (3,1),$\ldots$,(3,6),&  1  &   9    &   (9,2)   \\
&& (4,0), (4,1),(4,2), & 1  &   10    &   (8,3)   \\
&& (4,3), (4,5),(5,0), & 1  &   11    &   (7,4)   \\
&& (5,1), (5,2),(6,0), & 3  &   8    &   (7,1)   \\
&& (11,0). & 3  &  19    &  (8,0)    \\
4  &   20   &   (9,0)  & 4  &   21    &   (8,1)   \\
9  &   7    &   (3,8),(4,6) & 9  &   15    &   (7,2)\\
10  &   5    &   (4,7)& 10  &   9   &   (10,1)  \\
10  &   13    &   (5,4) & 11  &   5    &   (6,5) \\
12  &   2    &   (6,2)& 12  &   12   &   (3,7) \\
12  &   13    &   (1,10),(2,8),(6,1),(10,0) & 12  &   14    &   (9,1)\\
12  &   15    &   (8,2) & 12  &   16    &   (7,3) \\
12  &   17    &   (6,4) & 12  &   18    &   (5,5) \\
13  &   2    &   (0,11),(6,3) & 13  &   13    &   (1,9),(4,4),(5,3)   \\
\hline
\end{tabular}
\end{center}

For counting $Fin(v)$ for $11\leq v\leq 23$ and $v\neq 16,17$, we need to search for a large number of instances of maximum kite packings of order $v$ as we have done in Lemma \ref{10}. To reduce the computation, when $v\neq 22$, we shall first try to determine the fine triangle intersection numbers of a pair of maximum $(K_{v}\setminus K_{h_{v}},G)$-packings with the same vertex set and the same subgraph $K_{h_{v}}$ removed, where $G$ is a kite and
$$h_{v}=\left\{
    \begin{array}{ll}
        6,   &\mbox{if $v=11,14$},\\
        5,   &\mbox{if $v=12,13$},\\
        7,   &\mbox{if $v=15$},\\
        10,  &\mbox{if $v=18,19,20,23$},\\
        12,  &\mbox{if $v=21$}.
   \end{array}
    \right.$$
When $v=22$, we shall try to determine the fine triangle intersection numbers of a pair of kite-GDDs of type $8^{2}6^{1}$ with the same group set. These results will be listed in Lemmas \ref{11-6}-\ref{23-10}. Here in order to save space, we give only the detail of the proof of Lemma \ref{11-6}. For other details, the interested reader may find them in Appendix of this paper.

\begin{Lemma}\label{11-6}
Let $M_{11}=\{(s,t): s+t\leq 10, s\in \{0,1,10\}, t$ is a non-negative integer$\}$
and $N_{11}=\{(4,6),(5,0),(5,2),(6,4),(7,0),(8,1)\}$. Let $G$ be a kite and $(s,t)\in M_{11}\cup N_{11}$.
Then there is a pair of $(K_{11}\setminus K_{6},G)$-designs with the same vertex set
and the same subgraph $K_{6}$ removed, which intersect in $s$ blocks and
$t+s$ triangles.
\end{Lemma}

\proof Take the vertex set $X=\{0,1,\ldots,10\}$. Let
\begin{center}\small
 \begin{tabular}{lllllll}
${\cal B}_1:$ &$[0,10,9-3]$,&$[5,9,7-1]$,&$[1,9,8-5]$,&$[6,1,10-2]$,&$[2,8,7-0]$,&$[3,10,8-4]$,\\\vspace{0.1cm}&
$[3,7,6-4]$,&$[6,2,9-4]$,&$[4,7,10-5]$,&$[8,0,6-5]$;\\
${\cal B}_2:$ &$[0,10,9-4]$,&$[5,9,7-0]$,&$[1,9,8-4]$,&$[6,1,10-2]$,&$[2,8,7-1]$,&$[10,7,3-8]$,\\&
$[6,2,9-3]$,&$[8,0,6-3]$,&$[5,8,10-4]$,&$[7,4,6-5]$.
\end{tabular}
\end{center}

\noindent Let ${\cal B}_3=({\cal B}_1\setminus\{[3,10,8-4],[3,7,6-4],[4,7,10-5]\}\cup\{[10,7,3-8],[4,8,10-5],[7,4,6-3]\}$, and ${\cal B}_4=({\cal B}_3\setminus\{[1,9,8-5],[4,8,10-5]\}\cup\{[1,9,8-4],[5,8,10-4]\}$. Then $(X,{\cal B}_{i})$ is a $(K_{11}\setminus K_{6},G)$-design for each $1\leq i\leq 4$,
where the removed subgraph $K_{6}$ is constructed on $Y=\{0,1,2,3,4,5\}$. Consider the following permutations on $X$.

\begin{center}\small
\begin{tabular}{lll}
$\pi_{0,0}=(7\ 8\ 9\ 10)$,&
$\pi_{0,1}=(6\ 7\ 8\ 10\ 9)$,&
$\pi_{0,2}=(7\ 8)(9\ 10)$,\\
$\pi_{0,3}=(6\ 7\ 9)(8\ 10)$,&
$\pi_{0,4}=(6\ 7)(8\ 9)$,&
$\pi_{0,5}=(4\ 5)(6\ 8)(7\ 9\ 10)$,\\
$\pi_{0,6}=(6\ 9)(8\ 10)$,&
$\pi_{0,7}=(0\ 2)(1\ 5\ 3)(7\ 10)(8\ 9)$,&
$\pi_{0,8}=(0\ 2\ 4\ 1\ 5)(6\ 8\ 7\ 10)$,\\
$\pi_{0,9}=(0\ 1)(3\ 5)(6\ 9)$,&
$\pi_{0,10}=(2\ 3)(4\ 5)(6\ 8)(9\ 10)$,&
$\pi_{1,0}=(8\ 9\ 10)$,\\
$\pi_{1,1}=(0\ 1)(2\ 3)$,&
$\pi_{1,2}=(1\ 2)(3\ 4)$,&
$\pi_{1,3}=(0\ 3\ 2\ 1\ 5\ 4)(7\ 9\ 8)$,\\
$\pi_{1,4}=(0\ 2\ 1)(3\ 5\ 4)(6\ 9\ 8)$,&
$\pi_{1,5}=(2\ 4\ 5)$,&
$\pi_{1,6}=(2\ 3\ 4\ 5)(6\ 10)(8\ 9)$,\\
$\pi_{1,7}=(0\ 3)(1\ 4\ 5\ 2)(6\ 10\ 7\ 9)$,&
$\pi_{1,8}=(0\ 5\ 4)(1\ 2\ 3)(7\ 10\ 9)$,&
$\pi_{1,9}=(2\ 3\ 5\ 4)(6\ 8)(9\ 10)$,\\
$\pi_{4,6}=(1)$,&
$\pi_{5,0}=(0\ 4)(1\ 5\ 2)(6\ 9)(7\ 10)$,&
$\pi_{5,2}=(3\ 4)$,\\
$\pi_{6,4}=(3\ 4)$,&
$\pi_{7,0}=\pi_{8,1}=\pi_{10,0}=(1)$.
\end{tabular}
\end{center}
\noindent We have that for each row in Table $2$,
 $\pi_{s,t}Y=Y$, $|{\cal B}_i\cap \pi_{s,t}{\cal B}_j|=s$ and
 $|T({\cal B}_i\setminus\pi_{s,t}{\cal B}_j)\cap T(\pi_{s,t}{\cal B}_j\setminus{\cal B}_i)|=t$.  \qed

\begin{center}\small
Table $2$.  \ Fine triangle intersection numbers for $(K_{11}\setminus K_{6},G)$-designs
  \begin{tabular}{llc}\hline
$i$ &  $j$ &  $(s,t)$ \\
\hline     1  &   1   &  $(M_{11}\cup N_{11})\setminus \{(1,9),(4,6),(5,0),(6,4),(7,0),(8,1)\}$     \\
     1  &   3   &   (1,9),(6,4),(7,0)   \\
     2  &   4   &   (4,6)    \\
     3  &   4    &  (5,0),(8,1)    \\
\hline
\end{tabular}
\end{center}

\begin{Lemma}\label{12-5}
Let $M_{12}=\{(s,t): s+t\leq 14, s\in \{0,3,8,12,14\}, t$ is a non-negative integer$\}\setminus \{(3,9),(3,10),(8,2),(8,4),(8,5),(12,1)\}$
and $N_{12}=\{(4,10),(5,9),(6,0),(6,2),(6,5),(6,8),(7,4),\\
(7,7),(9,5),(10,0),(10,2),(10,4)\}$. Let $G$ be a kite and $(s,t)\in M_{12}\cup N_{12}$.
Then there is a pair of $(K_{12}\setminus K_{5},G)$-designs with the same vertex set
and the same subgraph $K_{5}$ removed, which intersect in $s$ blocks and
$t+s$ triangles.
\end{Lemma}

\begin{Lemma}\label{13-5}
Let $M_{13}=\{(s,t): s+t\leq 17, s\in \{0,3,6,9,15,17\} , t$ is a non-negative integer$\}$ $\setminus \{(0,15),(3,12),(6,8),(6,9),(9,4),(9,5),(9,7),(15,1)\}$, and
$N_{13}=\{(2,15),(5,11),(7,7),(8,9),\\(10,5),(11,4),(11,6),(12,0),(12,3),(13,0),(13,2),(14,3)\}$.
Let $G$ be a kite and $(s,t)\in M_{13}\cup N_{13}$.
Then there is a pair of $(K_{13}\setminus K_{5},G)$-designs with the same vertex set
and the same subgraph $K_{5}$ removed, which intersect in $s$ blocks and
$t+s$ triangles.
\end{Lemma}

\begin{Lemma}\label{14-6}
Let $M_{14}=\{(s,t): s+t\leq 19, s\in \{0,4,8,17,19\}, t$ is a non-negative integer$\}$ $\setminus \{(4,13),(8,8),(8,9),(8,10),(17,1)\}$,
and $N_{14}=\{(7,9),(7,11),(10,5),(10,9),(11,6),(12,0),(12,\\1),(12,2),(13,6),(14,0),(14,3),(15,2),(16,3)\}$.
Let $G$ be a kite and $(s,t)\in M_{14}\cup N_{14}$.
Then there is a pair of  $(K_{14}\setminus K_{6},G)$-designs with the same vertex set
and the same subgraph $K_{6}$ removed, which intersect in $s$ blocks and
$t+s$ triangles.
\end{Lemma}

\begin{Lemma}\label{15-7}
Let $M_{15}=\{(0,0),(0,1),\ldots,(0,16),(0,18),(0,21),(21,0)\}$
and $N_{15}=\{(4,17),(6,0),\\(6,1),\ldots,(6,8),(6,12),(8,13),(9,8),(12,0),(12,4),(12,6),(12,9),
(13,2),(13,8),(17,1),(17,4),\\(18,0)\}$.
Let $G$ be a kite and $(s,t)\in M_{15}\cup N_{15}$.
Then there is a pair of $(K_{15}\setminus K_{7},G)$-designs with the same vertex set
and the same subgraph $K_{7}$ removed, which intersect in $s$ blocks and
$t+s$ triangles.
\end{Lemma}

\begin{Lemma}\label{18-10}
Let $M_{18}=\{(0,0),(0,1),\ldots,(0,17),(0,27),(27,0)\}$ and $N_{18}=\{(3,24),(6,21),(8,\\0),(9,18),(12,5),(12,15),(15,0),(15,12),(18,9),(21,0),(24,3)\}$.
Let $G$ be a kite and $(s,t)\in M_{18}\cup N_{18}$.
Then there is a pair of $(K_{18}\setminus K_{10},G)$-designs with the same vertex set
and the same subgraph $K_{10}$ removed, which intersect in $s$ blocks and
$t+s$ triangles.
\end{Lemma}

\begin{Lemma}\label{19-10}
Let $M_{19}=\{(0,0),(0,1),\ldots,(0,18),(0,21),(0,31),(31,0)\}$
and $N_{19}=\{(2,17),(3,\\28),(6,25),(9,0),(9,22),(11,20),(12,4),(12,12),(15,3),(17,14),
(18,0),(18,3),(20,11),(22,9),\\(25,0),(25,6)\}$.
Let $G$ be a kite and $(s,t)\in M_{19}\cup N_{19}$.
Then there is a pair of  maximum $(K_{19}\setminus K_{10},G)$-packings with the same vertex set
and the same subgraph $K_{10}$ removed, which intersect in $s$ blocks and
$t+s$ triangles.
\end{Lemma}

\begin{Lemma}\label{20-10}
Let $M_{20}=\{(0,0),(0,1),\ldots,(0,17),(0,22),(0,24),(0,28),(0,36),(36,0)\}$
and $N_{20}$ $=\{(2,16),(5,31),(7,21),(10,18),(10,26),(12,0),(12,3),(12,6),(13,15),(13,23),
(15,13),(18,10),\\(18,18),(21,1),(21,7),(23,13),(24,0),(26,10),(28,0),(31,5)\}$.
Let $G$ be a kite and $(s,t)\in M_{20}\cup N_{20}$.
Then there is a pair of  maximum $(K_{20}\setminus K_{10},G)$-packings with the same vertex set
and the same subgraph $K_{10}$ removed, which intersect in $s$ blocks and
$t+s$ triangles.
\end{Lemma}

\begin{Lemma}\label{21-12}
Let $M_{21}=\{(0,0),(0,1),\ldots,(0,19),(0,24),(0,36),(36,0)\}$ and $N_{21}=\{(6,30),(7,\\0),\ldots,(7,7),(12,5),(12,24),(13,11),(18,18),(24,0),(24,12),(30,6)\}$.
Let $G$ be a kite and $(s,t)\in M_{21}\cup N_{21}$.
Then there is a pair of  $(K_{21}\setminus K_{12},G)$-designs with the same vertex set
and the same subgraph $K_{12}$ removed, which intersect in $s$ blocks and
$t+s$ triangles.
\end{Lemma}

\begin{Lemma}\label{22-8}
Let $M_{22}=\{(0,0),(0,1),\ldots,(0,20),(10,0),(10,1)\ldots,(10,6),(13,7),(40,0)\}$ and
$N_{22}=\{(0,28),(0,40),(7,33),(11,17),(14,26),(17,11),(21,19),(26,14),(28,0),(33,7)\}$.
Let $(s,t)$ $\in M_{22}\cup N_{22}$.
Then there is a pair of  kite-GDDs of type $8^{2}6^{1}$ with the same group set,
which intersect in $s$ blocks and $t+s$ triangles.
\end{Lemma}

\begin{Lemma}\label{23-10}
Let $M_{23}=\{(0,0),(0,1),\ldots,(0,18),(7,0),(7,6),(7,7),(7,10),(10,6)\}$ and
$N_{23}=\{(0,22),(0,26),(0,34),(0,42),(0,52),(3,49),(5,21),(5,29),(5,37),(6,46),
(7,19),(9,43),(10,24),\\(10,32),(12,14),(12,40),(14,0),(14,12),
(15,19),(15,27),(15,37),(18,8),(18,34),(19,15),(20,\\22),(21,5),(21,31),(22,0),
(22,20),(24,10),(24,28),(26,0),(27,15),(27,25),(28,24),(29,5),(31,\\21),(32,10),(34,0),(34,18),
(37,5),(37,15),(40,12),(42,0),(43,9),(46,6),(49,3),(52,0)\}$. Let $G$ be a kite and $(s,t)\in M_{23}\cup N_{23}$. Then there is a pair of  $(K_{23}\setminus K_{10},G)$-designs with the same vertex set and the same subgraph $K_{10}$ removed, which intersect in $s$ blocks and $t+s$ triangles.
\end{Lemma}

\begin{Lemma}\label{other small values}
$Fin(v)=Adm(v)$ for $11\leq v\leq 23$ and $v\neq 16,17,22$.
\end{Lemma}

\proof Obviously $Fin(v)\subseteq Adm(v)$. We need to show that $Adm(v)\subseteq Fin(v)$.
When $v=21$, our proof will rely on the fact that $Fin(12)=Adm(12)$. Thus one can fist make use of the following procedure to obtain $Fin(v)=Adm(v)$ for $11\leq v\leq 23$ and $v\neq 16,17,21,22$. Then the same procedure guarantees $Fin(21)=Adm(21)$.

For any $11\leq v\leq 23$ and $v\neq 16,17,22$, take the corresponding $M_v$ and $N_v$ from Lemmas \ref{11-6}-\ref{23-10}. Let $G$ be a kite and $(\alpha_{v},\beta_{v})\in M_{v}\cup N_{v}$. Let
$$h_{v}=\left\{
    \begin{array}{ll}
        6,   &\mbox{if $v=11,14$},\\
        5,   &\mbox{if $v=12,13$},\\
        7,   &\mbox{if $v=15$},\\
        10,  &\mbox{if $v=18,19,20,23$},\\
        12,  &\mbox{if $v=21$}.
   \end{array}
    \right.$$
By Lemmas \ref{11-6}-\ref{23-10}, there is a pair of $(K_{v} \setminus K_{h_{v}},G)$-designs (or maximum $(K_{v} \setminus K_{h_{v}},G)$-packings)
$(X,{\cal B}^{(v)}_{1})$ and $(X,{\cal B}^{(v)}_{2})$
with the same subgraph $K_{h_{v}}$ removed, which intersect in $\alpha_{v}$ blocks and
$\beta_{v}+\alpha_{v}$ triangles. Here the subgraph $K_{h_{v}}$ is constructed on $Y\subset X$.
Let $(\gamma_{v},\omega_{v})\in Adm(h_{v})$. By Lemmas \ref{Fin-5}-\ref{10}, there is a pair of maximum kite packings of order $h_{v}$,
$(Y,{\cal B}'^{(v)}_{1})$ and $(Y,{\cal B}'^{(v)}_{2})$, with $\gamma_{v}$ common blocks and
$\gamma_{v}+\omega_{v}$ common triangles. Then $(X,{\cal B}^{(v)}_{1}\cup {\cal B}'^{(v)}_{1})$ and
$(X,{\cal B}^{(v)}_{2}\cup {\cal B}'^{(v)}_{2})$ are both maximum kite packings of order $v$ with
$\alpha_{v}+\gamma_{v}$ common blocks and $\alpha_{v}+\beta_{v}+\gamma_{v}+\omega_{v}$ common triangles.
Thus we have
$$ Fin(v)\supseteq \{(\alpha_{v}+\gamma_{v},\beta_{v}+\omega_{v}): (\alpha_{v},\beta_{v})\in M_{v}\cup N_{v},(\gamma_{v},\omega_{v})\in Adm(h_{v})\}.$$
It is readily checked that for any pair of integers $(s,t)\in Adm(v)$,
we have $(s,t)\in Fin(v)$.  \qed

\begin{Lemma}\label{Fin(22)}
$Fin(22)= Adm(22)$.
\end{Lemma}

\proof Take the same sets $M_{22}$ and $N_{22}$ as in Lemma \ref{22-8}. Let $(\alpha,\beta)\in M_{22}\cup N_{22}$. By Lemma \ref{22-8}, there is
a pair of kite-GDDs of type $8^{2}6^{1}$ with the same group set,
which intersect in $\alpha$ blocks and $\beta+\alpha$ triangles. Let $(\gamma_1,\omega_1)$ and $(\gamma_2,\omega_2)\in$ $Adm(8)\setminus\{(6,0),(6,1)\}$. \cite{cflt} shows that there is a pair of kite systems of order $8$ intersecting in $\gamma_i$ common blocks and
$\omega_i+\gamma_i$ common triangles for each $i=1,2$. Let $(\gamma_3,\omega_3)\in Adm(6)$.  By Lemma \ref{Fin-6}, there is a pair of kite systems of order $6$ with $\gamma_3$ common blocks and $\omega_3+\gamma_3$ common triangles. Now applying Construction
\ref{FillingSubdesigns}, we obtain a pair of maximum kite packings of order $22$
with $\alpha+\sum_{i=1}^3\gamma_i$ common
blocks and $\beta+\sum_{i=1}^3\omega_i+(\alpha+\sum_{i=1}^3\gamma_i)$ common
triangles. Thus we have
\begin{eqnarray*}
Fin(22)\supseteq \hspace{-3mm}&\hspace{-3mm}& \{(\alpha+\sum_{i=1}^3\gamma_i,\beta+\sum_{i=1}^3\omega_i):(\alpha,\beta)\in M_{22}\cup N_{22},(\gamma_3,\omega_3)\in Adm(6),
\\\hspace{-3mm}&\hspace{-3mm}&\hspace{4.2cm}(\gamma_1,\omega_1),(\gamma_2,\omega_2)\in Adm(8)\setminus\{(6,0),(6,1)\}\}.
\end{eqnarray*}
It is readily checked that for any pair of integers $(s,t)\in Adm(22)$,
we have $(s,t)\in Fin(22)$.  \qed

\section{Input designs}

\begin{Lemma}\label{10-2}
Let $A_{10}=\{(0,0),(0,1),\ldots,(0,8),(0,10),(0,11),(11,0)\}$. Let $G$ be a kite and $(s,t)\in A_{10}$.
Then there is a pair of $(K_{10}\setminus K_{2},G)$-designs with the same vertex set
and the same subgraph $K_{2}$ removed, which intersect in $s$ blocks and
$t+s$ triangles.
\end{Lemma}

\proof Take the vertex set $X=\{0,1,\ldots,9\}$. Let
\begin{center}\small
 \begin{tabular}{lllllll}
${\cal B}_1:$ &$[2,0,1-9]$,&$[4,7,0-5]$,&$[1,4,8-2]$,&$[1,7,5-8]$,&$[2,4,3-1]$,&$[7,8,6-5]$,\\\vspace{0.1cm}&$[2,9,5-3]$,
&$[9,0,8-3]$,&$[0,3,6-1]$,&$[6,9,4-5]$,&$[9,3,7-2]$;\\

${\cal B}_2:$ &$[1 ,2 ,0- 5]$,&$[4 ,0 ,7- 2]$,&$[1 ,4 ,8- 3]$,&$[1 ,7 ,5 -4]$,&$[2 ,4, 3 -5]$,&$[7 ,8 ,6 -1]$,\\&$[2 ,9, 5- 8]$,
&$[9 ,0, 8- 2]$,&$[0 ,3, 6- 5]$,&$[4, 6 ,9- 1]$,&$[7 ,9, 3- 1]$.
\end{tabular}
\end{center}
\noindent Then $(X,{\cal B}_{i})$ is a $(K_{10}\setminus K_{2},G)$-design for $i=1,2$,
where the removed subgraph $K_{2}$ is constructed on $Y=\{2,6\}$. Consider the following permutations on $X$.

\begin{center}\small
\begin{tabular}{lll}
$\pi_{0,0}=(0\ 4\ 3\ 7\ 5\ 9)(1\ 8)(2\ 6)$,&
$\pi_{0,1}=(1\ 4)(3\ 5)(7\ 8)$,&
$\pi_{0,2}=(0\ 4)(1\ 3)(5\ 7\ 8\ 9)$,\\
$\pi_{0,3}=(0\ 9\ 8\ 3\ 1\ 7\ 4\ 5)(2\ 6)$,&
$\pi_{0,4}=(0\ 5\ 8\ 1\ 9\ 7)$,&
$\pi_{0,5}=(0\ 7\ 9\ 1\ 8\ 5)(3\ 4)$,\\
$\pi_{0,6}=(0\ 4\ 5\ 3\ 9)(7\ 8)$,&
$\pi_{0,7}=(3\ 5)(4\ 9)(7\ 8)$,&
$\pi_{0,8}=(0\ 1)(5\ 9)(7\ 8)$,\\
$\pi_{0,10}=(0\ 9)(1\ 5)(3\ 4)$,&
$\pi_{0,11}=\pi_{11,0}=(1)$.\\
\end{tabular}
\end{center}

\noindent We have that $\pi_{s,t}Y=Y$. For each $(s,t)\in A_{10}\setminus (0,11)$,
$|{\cal B}_1\cap \pi_{s,t}{\cal B}_1|=s$ and
 $|T({\cal B}_1\setminus\pi_{s,t}{\cal B}_1)\cap T(\pi_{s,t}{\cal B}_1\setminus{\cal B}_1)|=t$.
$|{\cal B}_1\cap \pi_{s,t}{\cal B}_2|=0$ and
 $|T({\cal B}_1\setminus\pi_{s,t}{\cal B}_2)\cap T(\pi_{s,t}{\cal B}_2\setminus{\cal B}_1)|=11$.\qed

\begin{Lemma}\label{11-3}
Let $A_{11}=\{(0,0),(0,1),\ldots,(0,9),(0,11),(0,13),(13,0)\}$. Let $G$ be a kite and $(s,t)\in A_{11}$.
Then there is a pair of $(K_{11}\setminus K_{3},G)$-designs with the same vertex set
and the same subgraph $K_{3}$ removed, which intersect in $s$ blocks and
$t+s$ triangles.
\end{Lemma}

\proof Take the vertex set $X=\{0,1,\ldots,10\}$. Let
\begin{center}
\small
 \begin{tabular}{lllllll}

${\cal B}_1:$ &$[0 ,10, 9- 3]$,&$[5 ,9, 7- 1]$,&$[1 ,9, 8 -5]$,&$[6 ,1 ,10- 2]$,&$[2 ,8, 7- 0]$,&$[3 ,10, 8 -4]$,\\&$[3 ,7, 6 -4]$,
&$[6, 2 ,9- 4]$,&$[4 ,7 ,10- 5]$,&$[8, 0 ,6- 5]$,&$[0 ,1, 2 -5]$,&$[2 ,3 ,4 -1]$,\\\vspace{0.1cm}&$[4, 5, 0 -3]$;\\

${\cal B}_2:$ &$[0 ,10, 9 -4]$,&$[5 ,9, 7- 0]$,&$[1 ,9 ,8- 4]$,&$[6 ,1, 10 -5]$,&$[2, 8, 7- 1]$,&$[3 ,10, 8 -5]$,\\&$[3, 7 ,6- 5]$,
&$[6 ,2 ,9 -3]$,&$[4 ,7, 10- 2]$,&$[8, 0, 6 -4]$,&$[2, 1, 0- 3]$,&$[4 ,3, 2- 5]$,\\&$[0 ,5 ,4- 1]$.
\end{tabular}
\end{center}
\noindent Then $(X,{\cal B}_{i})$ is a $(K_{11}\setminus K_{3},G)$-design for $i=1,2$, where the removed subgraph $K_{3}$ is constructed on $Y=\{1,3,5\}$. Consider the following permutations on $X$.

\begin{center}
\small
 \begin{tabular}{lll}
$\pi_{0,0}=(0\ 9\ 8\ 2\ 6\ 4\ 10)$,&
$\pi_{0,1}=(0\ 6\ 7\ 10\ 8\ 4\ 9\ 2)(1\ 3)$,&
$\pi_{0,2}=(0\ 6\ 10\ 9)(1\ 5\ 3)(2\ 7)$,\\
$\pi_{0,3}=(0\ 10\ 4\ 8\ 7)(2\ 6\ 9)(3\ 5)$,&
$\pi_{0,4}=(0\ 9\ 4\ 7\ 10\ 8\ 2)(1\ 3)$,&
$\pi_{0,5}=(0\ 9\ 6)(2\ 8)(3\ 5)$,\\
$\pi_{0,6}=(0\ 2)(1\ 3\ 5)(4\ 10\ 7)(8\ 9)$,&
$\pi_{0,7}=(0\ 4)(1\ 5\ 3)(6\ 10\ 8\ 9\ 7)$,&
$\pi_{0,8}=(0\ 6)(2\ 10)(3\ 5)(4\ 7)$,\\
$\pi_{0,9}=(0\ 8)(1\ 3)(2\ 10)(4\ 6)(7\ 9)$,&
$\pi_{0,11}=(0\ 2)(3\ 5)(7\ 10)(8\ 9)$,&
$\pi_{0,13}=\pi_{13,0}=(1)$.
\end{tabular}
\end{center}

\noindent We have that $\pi_{s,t}Y=Y$. For each $(s,t)\in A_{11}\setminus (0,13)$,
$|{\cal B}_1\cap \pi_{s,t}{\cal B}_1|=s$ and
 $|T({\cal B}_1\setminus\pi_{s,t}{\cal B}_1)\cap T(\pi_{s,t}{\cal B}_1\setminus{\cal B}_1)|=t$.
$|{\cal B}_1\cap \pi_{s,t}{\cal B}_2|=0$ and
 $|T({\cal B}_1\setminus\pi_{s,t}{\cal B}_2)\cap T(\pi_{s,t}{\cal B}_2\setminus{\cal B}_1)|=13$.\qed

\begin{Lemma}\label{12-4}
Let $A_{12}=\{(0,0),(0,1),\ldots,(0,11),(0,14),(0,15),(15,0)\}$. Let $G$ be a kite and $(s,t)\in A_{12}$.
Then there is a pair of $(K_{12}\setminus K_{4},G)$-designs with the same vertex set
and the same subgraph $K_{4}$ removed, which intersect in $s$ blocks and
$t+s$ triangles.
\end{Lemma}

\proof Take the vertex set $X=\{0,1,\ldots,11\}$. Let
\begin{center}
\small
 \begin{tabular}{lllllll}
${\cal B}_1:$ &$[8 ,5, 6- 7]$,&$[6, 11, 0- 8]$,&$[11 ,10, 5 -2]$,&$[4 ,9 ,10 -3]$,&$[9, 8, 3- 5]$,&$[8 ,7, 2- 4]$,\\&$[1, 11, 7- 4]$,
&$[7 ,10, 0 -9]$,&$[2 ,9 ,11- 3]$,&$[10 ,8 ,1- 5]$,&$[9 ,5 ,7- 3]$,&$[1 ,9, 6 -3]$,\\\vspace{0.1cm}&$[2 ,10 ,6 -4]$,&$[11, 8, 4- 3]$,
&$[0 ,5, 4 -1]$;\\

${\cal B}_2:$ &$[8 ,5, 6- 4]$,&$[6 ,11, 0- 9]$,&$[11 ,10 ,5 -1]$,&$[9 ,10 ,4 -3]$,&$[9 ,8 ,3 -11]$,&$[8 ,7, 2- 5]$,\\&$[1 ,11, 7- 3]$,
&$[7 ,10, 0- 8]$,&$[9, 11 ,2- 4]$,&$[8 ,1, 10 -3]$,&$[9 ,5 ,7 -4]$,&$[1, 9, 6 -7]$,\\&$[2 ,10 ,6- 3]$,&$[11, 8 ,4- 1]$,
&$[0 ,4 ,5 -3]$.
\end{tabular}
\end{center}

\noindent Let ${\cal B}_3=({\cal B}_2\setminus\{[9 ,10 ,4 -3],[0 ,4 ,5 -3]\}\cup\{[9 ,10, 4- 0],[4, 3 ,5- 0]\}$. Then $(X,{\cal B}_{i})$ is a $(K_{12}\setminus K_{4},G)$-design for each $1\leq i\leq 3$, where the removed subgraph $K_{4}$ is constructed on $Y=\{0,1,2,3\}$. Consider the following permutations on $X$.

\begin{center}\small\tabcolsep 0.05in
 \begin{tabular}{lll}
$\pi_{0,0}=(2\ 3)(4\ 11\ 9\ 10\ 7\ 8\ 5)$,&
$\pi_{0,1}=(0\ 3)(1\ 2)(4\ 5\ 9\ 8\ 10\ 11\ 7\ 6)$,&
$\pi_{0,2}=(1\ 2)(5\ 8\ 9\ 10\ 6\ 7\ 11)$,\\
$\pi_{0,3}=(0\ 3\ 1)(4\ 5)(7\ 9\ 11\ 10\ 8)$,&
$\pi_{0,4}=(0\ 3)(4\ 5\ 8\ 7)(6\ 9)(10\ 11)$,&
$\pi_{0,5}=(0\ 2\ 3)(6\ 11\ 9\ 7\ 8\ 10)$,\\
$\pi_{0,6}=(0\ 1\ 3\ 2)(4\ 8\ 6\ 7\ 9)(5\ 10\ 11)$,&
$\pi_{0,7}=(1\ 2)(7\ 11\ 9\ 10)$,&
$\pi_{0,8}=(0\ 2\ 1)(7\ 10\ 9\ 11)$,\\
$\pi_{0,9}=(0\ 3)(1\ 2)(6\ 8)(7\ 9)$,&
$\pi_{0,10}=(1\ 2)(6\ 8)(7\ 9)$,&
$\pi_{0,11}=(0\ 2)(4\ 5)(6\ 9)(8\ 10)$,\\
$\pi_{0,14}=\pi_{0,15}=\pi_{15,0}=(1)$.
\end{tabular}
\end{center}

\noindent We have that $\pi_{s,t}Y=Y$. For each $(s,t)\in A_{12}\setminus \{(0,14),(0,15)\}$,
$|{\cal B}_1\cap \pi_{s,t}{\cal B}_1|=s$ and
 $|T({\cal B}_1\setminus\pi_{s,t}{\cal B}_1)\cap T(\pi_{s,t}{\cal B}_1\setminus{\cal B}_1)|=t$.
$|{\cal B}_1\cap \pi_{s,t}{\cal B}_2|=0$ and
 $|T({\cal B}_1\setminus\pi_{s,t}{\cal B}_2)\cap T(\pi_{s,t}{\cal B}_2\setminus{\cal B}_1)|=15$.
$|{\cal B}_1\cap \pi_{s,t}{\cal B}_3|=0$ and
 $|T({\cal B}_1\setminus\pi_{s,t}{\cal B}_3)\cap T(\pi_{s,t}{\cal B}_3\setminus{\cal B}_1)|=14$.\qed

\begin{Lemma}{\rm (\cite{cflt})}\label{kite-GDD-4^3}
Let $(s,t)\in \{(0,0),(0,1),(0,12),(12,0)\}$. There exists a pair of kite-GDDs of type
$4^{3}$ with the same group set intersecting in $s$ blocks and $t+s$ triangles.
\end{Lemma}

\section{Putting the pieces together}

First we need the following definition. Let $(s_{1},t_{1})$ and  $(s_{2},t_{2})$ be two pairs of non-negative
integers. Define $(s_{1},t_{1})+(s_{2},t_{2})=(s_{1}+s_{2},t_{1}+t_{2})$. If $X$ and $Y$ are two sets of pairs of
non-negative integers, then $X+Y$ denotes the set $\{(s_{1},t_{1})+(s_{2},t_{2}):(s_{1},t_{1})\in X,(s_{2},t_{2})\in Y\}$.
If $X$ is a set of pairs of non-negative integers and $h$ is some positive integer, then $h\ast X$ denotes the
set of all pairs of non-negative integers which can be obtained by adding any $h$ elements of $X$ together
(repetitions of elements of $X$ allowed).

\begin{Lemma}\label{L_1}
For any integer $v\equiv 2,3,4,5,6,7,10,11,12,13,14,15\ ({\rm mod }\ 24)$ and $v\geq 4$, $Fin(v)=Adm(v)$.
\end{Lemma}

\proof  The cases of $v=4,5,6,7,10,11,12,13,14,15$ follows from Example \ref{Fin-4}, Lemmas \ref{Fin-5}-\ref{10} and Lemma \ref{other small values}.
Assume that $v\geq 24$. Obviously $Fin(v)\subseteq Adm(v)$. We need to show that $Adm(v)\subseteq Fin(v)$.
Let $v=8u+a$ with $u\equiv 0,1$ (mod 3), $u\geq 3$ and $a\in \{2,3,4,5,6,7\}$.

Start from a $3$-GDD of type $2^{u}$
from Lemma \ref{3-GDD}. Give each point of the GDD weight $4$. By Lemma \ref{kite-GDD-4^3}, there is a pair of kite-GDDs of type $4^{3}$
with $\alpha$ common blocks and $\alpha+\beta$ common triangles, $(\alpha,\beta)\in \{(0,0),(0,12),(12,0)\}$.
Then apply Construction \ref{WeightingConstruction} to obtain a pair of kite-GDDs of type $8^{u}$ with $\sum^{x}_{i=1}\alpha_{i}$ common blocks
and $\sum^{x}_{i=1}(\alpha_{i}+\beta_{i})$ common triangles, where $x=2u(u-1)/3$ is the number of blocks of the $3$-GDD of
type $2^{u}$ and $(\alpha_{i},\beta_{i})\in \{(0,0),(0,12),(12,0)\}$ for $1\leq i\leq x$.

Define $A(8+c)=A_{8+c}$, where $c\in \{2,3,4\}$ and $A_{8+c}$ is taken from Lemmas \ref{10-2}-\ref{12-4}. When $c\in \{5,6,7\}$, define $A(8+c)=M_{8+c}\cup N_{8+c}$, where $M_{8+c}$ and $N_{8+c}$ are taken from Lemmas \ref{13-5}-\ref{15-7}. By Lemmas \ref{13-5}-\ref{15-7} and Lemmas \ref{10-2}-\ref{12-4}, there is a pair of $(K_{8+a}\setminus K_{a},G)$-designs with the same vertex set and the same subgraph $K_{a}$ removed with $\gamma_{j}$ common blocks and $\gamma_{j}+\omega_{j}$ common triangles, where $(\gamma_{j},\omega_{j})\in A(8+a)$, $1\leq j\leq u-1$, and $a\in \{2,3,4,5,6,7\}$. By Lemmas \ref{10} and \ref{other small values}, there is a pair of maximum kite packings of order $8+a$
with $\gamma_{u}$ common blocks and $\gamma_{u}+\omega_{u}$ common triangles, $(\gamma_{u},\omega_{u})\in Adm(8+a)$.
Apply Construction \ref{FillingSubdesigns} to obtain a pair of maximum kite packings of order $v=8u+a$ with
$\sum^{x}_{i=1}\alpha_{i}+\sum^{u}_{j=1}\gamma_{j}$ common blocks and
$\sum^{x}_{i=1}(\alpha_{i}+\beta_{i})+\sum^{u}_{j=1}(\gamma_{j}+\omega_{j})$ common triangles. Thus we have
\begin{eqnarray*}
Fin(v)\supseteq \hspace{-3mm}&\hspace{-3mm}& \{(\sum^{x}_{i=1}\alpha_{i}+\sum^{u}_{j=1}\gamma_{j},\sum^{x}_{i=1}\beta_{i}+\sum^{u}_{j=1}\omega_{j}):
(\alpha_{i},\beta_{i})\in \{(0,0),(0,12),(12,0)\}, 1\leq i\leq x, \\\hspace{-3mm}&\hspace{-3mm}&\hspace{3cm}
(\gamma_{j},\omega_{j})\in A(8+a), 1\leq j\leq u-1, (\gamma_{u},\omega_{u})\in Adm(8+a)\}\\
=\hspace{-3mm}&\hspace{-3mm}&
\{\sum^{x}_{i=1}(\alpha_{i},\beta_{i})+\sum^{u}_{j=1}(\gamma_{j},\omega_{j}):
(\alpha_{i},\beta_{i})\in \{(0,0),(0,12),(12,0)\}, 1\leq i\leq x, \\\hspace{-3mm}&\hspace{-3mm}&\hspace{3cm}
(\gamma_{j},\omega_{j})\in A(8+a), 1\leq j\leq u-1, (\gamma_{u},\omega_{u})\in Adm(8+a)\}\\
=\hspace{-3mm}&\hspace{-3mm}&
x\ast \{(0,0),(0,12),(12,0)\}+(u-1)\ast A(8+a)+Adm(8+a).
\end{eqnarray*}
Let $n$ be the number of blocks of a $(K_{8+a}\setminus K_{a},G)$-design. Then $n=7+2a$. Let
$m$ be the number of blocks of a maximum kite packing of order $8+a$.
Define $S(z)=\{(s,t): s+t\leq z,s,t\geq0\}$. Then it is readily checked that
\begin{eqnarray*}
& &(u-1)\ast A(8+a)+Adm(8+a)\\
& &=(u-1)\ast A(8+a)+S(m)\\
& &=(u-2)\ast A(8+a)+[A(8+a)+S(m)]\\
& &=(u-2)\ast A(8+a)+S(m+n)\\
& &=S((u-1)n+m).
\end{eqnarray*}
Note that
\begin{equation}
\label{formu}
\{(0,0),(0,12),(12,0)\}+S(z)=S(z+12)
\end{equation}
for any integer $z\geq 22$. Since $(u-1)n+m\geq 22$, using Formula $(1)$ inductively $x$
times, we have
\begin{eqnarray*}
Fin(v)\supseteq \hspace{-3mm}&\hspace{-3mm}&
(x-1)*\{(0,0),(0,12),(12,0)\}+(\{(0,0),(0,12),(12,0)\}+
S((u-1)n+m))\\
=\hspace{-3mm}&\hspace{-3mm}&
(x-1)*\{(0,0),(0,12),(12,0)\}+S((u-1)n+m+12)\\
=\hspace{-3mm}&\hspace{-3mm}&
S((u-1)n+m+12x)=S(b_v)=Adm(v).
\end{eqnarray*}
This completes the proof. \qed

\begin{Lemma}\label{L_2}
For any integer $v\equiv 18,19,20,21,22,23\ ({\rm mod }\ 24)$ and $v\geq 18$, $Fin(v)=Adm(v)$.
\end{Lemma}

\proof  The cases of $v=18,19,20,21,22,23$ follows from Lemmas \ref{other small values} and \ref{Fin(22)}. Assume that $v\geq 42$. Obviously $Fin(v)\subseteq Adm(v)$. We need to show that $Adm(v)\subseteq Fin(v)$.
Let $v=8u+a$ with $u\equiv 2$ (mod $3$), $u\geq 5$ and $a\in \{2,3,4,5,6,7\}$.

Start from a $3$-GDD of type $2^{u-2}4^{1}$
from Lemma \ref{3-GDD}. Give each point of the GDD weight $4$. By Lemma \ref{kite-GDD-4^3}, there is a pair of kite-GDDs of type $4^{3}$
with $\alpha$ common blocks and $\alpha+\beta$ common triangles, $(\alpha,\beta)\in \{(0,0),(0,12),(12,0)\}$.
Then apply Construction \ref{WeightingConstruction} to obtain a pair of kite-GDDs of type $8^{u-2}16^{1}$ with $\sum^{x}_{i=1}\alpha_{i}$ common blocks
and $\sum^{x}_{i=1}(\alpha_{i}+\beta_{i})$ common triangles, where $x=2(u+1)(u-2)/3$ is the number of blocks of the $3$-GDD of
type $2^{u-2}4^{1}$ and $(\alpha_{i},\beta_{i})\in \{(0,0),(0,12),(12,0)\}$ for $1\leq i\leq x$.

Define $A(8+c)=A_{8+c}$, where $c\in \{2,3,4\}$ and $A_{8+c}$ is taken from Lemmas \ref{10-2}-\ref{12-4}. When $c\in \{5,6,7\}$, define $A(8+c)=M_{8+c}\cup N_{8+c}$, where $M_{8+c}$ and $N_{8+c}$ are taken from Lemmas \ref{13-5}-\ref{15-7}. By Lemmas \ref{13-5}-\ref{15-7} and Lemmas \ref{10-2}-\ref{12-4}, there is a pair of $(K_{8+a}\setminus K_{a},G)$-designs with the same vertex set and the same subgraph $K_{a}$ removed with $\gamma_{j}$ common blocks and $\gamma_{j}+\omega_{j}$ common triangles, where $(\gamma_{j},\omega_{j})\in A(8+a)$, $1\leq j\leq u-2$, and $a\in \{2,3,4,5,6,7\}$. By Lemmas \ref{other small values} and \ref{Fin(22)}, there is a pair of maximum kite packings of order $16+a$
with $\gamma_{u-1}$ common blocks and $\gamma_{u-1}+\omega_{u-1}$ common triangles, $(\gamma_{u-1},\omega_{u-1})\in Adm(16+a)$.
Apply Construction \ref{FillingSubdesigns} to obtain a pair of maximum kite packings of order $v=8u+a$ with
$\sum^{x}_{i=1}\alpha_{i}+\sum^{u-1}_{j=1}\gamma_{j}$ common blocks and
$\sum^{x}_{i=1}(\alpha_{i}+\beta_{i})+\sum^{u-1}_{j=1}(\gamma_{j}+\omega_{j})$ common triangles. Thus we have
\begin{eqnarray*}
Fin(v)\supseteq \hspace{-3mm}&\hspace{-3mm}& \{(\sum^{x}_{i=1}\alpha_{i}+\sum^{u-1}_{j=1}\gamma_{j},\sum^{x}_{i=1}\beta_{i}+\sum^{u-1}_{j=1}\omega_{j}):
(\alpha_{i},\beta_{i})\in \{(0,0),(0,12),(12,0)\}, 1\leq i\leq x, \\\hspace{-3mm}&\hspace{-3mm}&\hspace{3cm}
(\gamma_{j},\omega_{j})\in A(8+a), 1\leq j\leq u-2, (\gamma_{u-1},\omega_{u-1})\in Adm(16+a)\}\\
=\hspace{-3mm}&\hspace{-3mm}&
\{\sum^{x}_{i=1}(\alpha_{i},\beta_{i})+\sum^{u-1}_{j=1}(\gamma_{j},\omega_{j}):
(\alpha_{i},\beta_{i})\in \{(0,0),(0,12),(12,0)\}, 1\leq i\leq x, \\\hspace{-3mm}&\hspace{-3mm}&\hspace{3cm}
(\gamma_{j},\omega_{j})\in A(8+a), 1\leq j\leq u-2, (\gamma_{u-1},\omega_{u-1})\in Adm(16+a)\}\\
=\hspace{-3mm}&\hspace{-3mm}&
x\ast \{(0,0),(0,12),(12,0)\}+(u-2)\ast A(8+a)+Adm(16+a).
\end{eqnarray*}
Let $n$ be the number of blocks of a $(K_{8+a}\setminus K_{a},G)$-design. Then $n=7+2a$. Let
$m$ be the number of blocks of a maximum kite packing of order $16+a$.
Define $S(z)=\{(s,t): s+t\leq z,s,t\geq0\}$. Then it is readily checked that
\begin{eqnarray*}
& &(u-2)\ast A(8+a)+Adm(16+a)\\
& &=(u-2)\ast A(8+a)+S(m)\\
& &=(u-3)\ast A(8+a)+[A(8+a)+S(m)]\\
& &=(u-3)\ast A(8+a)+S(m+n)\\
& &=S((u-2)n+m).
\end{eqnarray*}
Note that
\begin{equation}
\label{formu}
\{(0,0),(0,12),(12,0)\}+S(z)=S(z+12)
\end{equation}
for any integer $z\geq 22$. Since $(u-1)n+m\geq 22$, using Formula $(2)$ inductively $x$
times, we have
\begin{eqnarray*}
Fin(v)\supseteq \hspace{-3mm}&\hspace{-3mm}&
(x-1)*\{(0,0),(0,12),(12,0)\}+(\{(0,0),(0,12),(12,0)\}+
S((u-2)n+m))\\
=\hspace{-3mm}&\hspace{-3mm}&
(x-1)*\{(0,0),(0,12),(12,0)\}+S((u-2)n+m+12)\\
=\hspace{-3mm}&\hspace{-3mm}&
S((u-2)n+m+12x)=S(b_v)=Adm(v).
\end{eqnarray*}
This completes the proof. \qed

\noindent \textbf{Proof of Theorem \ref{main theorem}:} When $v\equiv 0,1\ ({\rm mod }\ 8)$ and $v\geq 8$, \cite{cflt} shows that $Fin(v)= Adm(v)\setminus \{(b_v-1,0),(b_v-1,1)\}$. When $v\equiv 2,3,4,5,6,7\ ({\rm mod }\ 8)$ and $v\geq 4$, combining the results of Lemmas \ref{L_1} and \ref{L_2}, we have $Fin(v)=Adm(v)$.  \qed

\newpage

\appendix
\section{Appendix}


\noindent {\bf Lemma $3.6$} {\em
Let $M_{12}=\{(s,t): s+t\leq 14, s\in \{0,3,8,12,14\}, t$ is a non-negative integer$\}\setminus \{(3,9),(3,10),(8,2),(8,4),(8,5),(12,1)\}$
and $N_{12}=\{(4,10),(5,9),(6,0),(6,2),(6,5),(6,8),(7,4),\\
(7,7),(9,5),(10,0),(10,2),(10,4)\}$. Let $G$ be a kite and $(s,t)\in M_{12}\cup N_{12}$.
Then there is a pair of $(K_{12}\setminus K_{5},G)$-designs with the same vertex set
and the same subgraph $K_{5}$ removed, which intersect in $s$ blocks and
$t+s$ triangles.}

\proof Take the vertex set $X=\{0,1,\ldots,11\}$. Let
\begin{center}\small
\begin{tabular}{lllllll}
${\cal B}_1:$&$[8,5,6-7]$,&$[6,11,0-8]$,&$[11,10,5-0]$,&$[4,9,10-3]$,&$[9,8,3-5]$, &$[8,7,2-5]$,\\
&$[1,11,7-4]$,&$[7,10,0-9]$,&$[2,9,11-3]$,&$[10,8,1-5]$,&$[9,5,7-3]$,&$[1,9,6-3]$,\\\vspace{0.1cm}
&$[2,10,6-4]$,&$[11,8,4-5]$;\\

${\cal B}_2:$&$[8,5,6-4]$,&$[6,11,0-9]$,&$[11,10,5-3]$,&$[9,10,4-5]$,&$[9,8,3-11]$, &$[8,2,7-3]$,\\
&$[11,7,1-5]$,&$[7,10,0-8]$,&$[9,11,2-5]$,&$[8,1,10-3]$,&$[9,7,5-0]$,&$[1,9,6-7]$,\\\vspace{0.1cm}
&$[2,10,6-3]$,&$[8,11,4-7]$;\\

${\cal B}_3:$&$[0,11,10-1]$,&$[10,9,8-2]$,&$[6,3,11-9]$,&$[8,7,6-0]$,&$[9,5,2-7]$, &$[10,6,2-11]$,\\&
$[4,8,11-7]$,&$[1,6,9-0]$,&$[0,5,8-1]$,&$[11,5,1-7]$,&$[5,10,3-8]$,&$[5,6,4-9]$,\\
&$[10,4,7-0]$,&$[9,3,7-5]$.\\
\end{tabular}
\end{center}

\noindent Let ${\cal B}_4=({\cal B}_1\setminus\{[8,5,6-7],[8,7,2-5]\})\cup\{[8,7,6-5],[8,5,2-7]\}$, ${\cal B}_5=({\cal B}_2\setminus\{[8,5,6-4],[1,9,6-7],[2,10,6-3]\})\cup\{[8,5,6-7],[1,9,6-3],[2,10,6-4]\}$,
${\cal B}_6=({\cal B}_5\setminus\{[6,11,0-9],[7,10,0-8]\})\cup\{[6,11,0-8],[7,10,0-9]\}$,
${\cal B}_7=({\cal B}_2\setminus\{[9,10,4-5],[11,7,1-5],[8,1,10-3],[8,11,4-7]\})\cup\{[4,9,10-3],[1,11,7-4],
[10,8,1-5],[8,11,4-5]\}$, ${\cal B}_8=({\cal B}_3\setminus\{[6,3,11-9],[10,6,2-11],[5,10,3-8],[9,3,7-5]\})\cup
\{[3,5,7-9],[11,6,2-10],[11,9,3-8],[3,6,10-5]\}$, ${\cal B}_9=({\cal B}_4\setminus\{[1,9,6-3],[2,10,6-4]\})\cup\{[1,9,6-4],[2,10,6-3]\}$, and ${\cal B}_{10}=({\cal B}_4\setminus\{[6,11,0-8],[7,10,0-9],[1,11,7-4],[9,5,7-3]\})\cup\{[6,11,0-9],[7,10,0-8],
[1,11,7-3],[9,5,7-4]\}$.
Then $(X,{\cal B}_{i})$ is a $(K_{12}\setminus K_{5},G)$-design for each $1\leq i\leq 10$,
where the removed subgraph $K_{5}$ is constructed on $Y=\{0,1,2,3,4\}$. Consider the following permutations on $X$.

\begin{center}\small
\begin{tabular}{lll}
$\pi_{0,0}=(8\ 9\ 10\ 11)$,&
$\pi_{0,1}=(7\ 8\ 9\ 10)$,&
$\pi_{0,2}=(8\ 9\ 11\ 10)$,\\
$\pi_{0,3}=(7\ 8\ 10)(9\ 11)$,&
$\pi_{0,4}=(6\ 7\ 8)(9\ 11\ 10)$,&
$\pi_{0,5}=(7\ 10)(9\ 11)$,\\
$\pi_{0,6}=(6\ 7\ 10\ 8)(9\ 11)$,&
$\pi_{0,7}=(6\ 8)(7\ 10)(9\ 11)$,&
$\pi_{0,8}=(0\ 4)(6\ 11\ 8)(7\ 10\ 9)$,\\
$\pi_{0,9}=(2\ 3)(6\ 8)(7\ 11)(9\ 10)$,&
$\pi_{0,10}=(0\ 4)(1\ 3\ 2)(6\ 8)(7\ 9)$,&
$\pi_{0,11}=(0\ 2)(6\ 8)(7\ 11)(9\ 10)$,\\
$\pi_{0,12}=(0\ 3)(1\ 4)(6\ 8)(7\ 10)(9\ 11)$,&
$\pi_{0,13}=(0\ 4)(1\ 2)(6\ 8)(7\ 9)$,&
$\pi_{0,14}=(1)$,\\
$\pi_{3,0}=(6\ 7\ 8)$,&
$\pi_{3,1}=(0\ 1)(2\ 4)$,&
$\pi_{3,2}=(3\ 4)(10\ 11)$,\\
$\pi_{3,3}=(0\ 4)(7\ 11)$,&
$\pi_{3,4}=(6\ 7\ 11\ 10)(8\ 9)$,&
$\pi_{3,5}=(0\ 4)(2\ 3)(7\ 9)$,\\
$\pi_{3,6}=(0\ 3\ 4)$,&
$\pi_{3,7}=(0\ 4\ 2\ 1\ 3)(7\ 10\ 9\ 11)$,&
$\pi_{3,8}=(0\ 1)(2\ 4)(6\ 8)(10\ 11)$,\\
$\pi_{3,11}=(1)$,&
$\pi_{4,10}=(1\ 2)(6\ 7)(8\ 9)(10\ 11)$,&
$\pi_{5,9}=(1)$,\\
$\pi_{6,0}=(1\ 4\ 2)(10\ 11)$,&
$\pi_{6,2}=(2\ 3)(7\ 9)$,&
$\pi_{6,5}=(3\ 4)$,\\
$\pi_{6,8}=(1\ 2)(6\ 7)(8\ 9)(10\ 11)$,&
$\pi_{7,4}=(5\ 7)$,&
$\pi_{7,7}=(1)$,\\
$\pi_{8,0}=(1\ 2)$,&
$\pi_{8,1}=(0\ 1)$,&
$\pi_{8,3}=(2\ 3)$,\\
$\pi_{8,6}=\pi_{9,5}=\pi_{10,0}=(1)$,&
$\pi_{10,2}=\pi_{10,4}=\pi_{12,0}=(1)$,&
$\pi_{12,2}=\pi_{14,0}=(1)$.
\end{tabular}
\end{center}
\noindent We have that for each $i,j,s,t$ in Table I,
 $\pi_{s,t}Y=Y$, $|{\cal B}_i\cap \pi_{s,t}{\cal B}_j|=s$ and
 $|T({\cal B}_i\setminus\pi_{s,t}{\cal B}_j)\cap T(\pi_{s,t}{\cal B}_j\setminus{\cal B}_i)|=t$.\qed

\newpage
\begin{center}
\small  Table I.  \ Fine triangle intersections for $(K_{12}\setminus K_{5},G)$-designs
  \begin{tabular}{ccc|ccc|ccc}\hline
$i$ &  $j$ &  $(s,t)$ & $i$ &  $j$ &  $(s,t)$& $i$ &  $j$ &  $(s,t)$ \\\hline
1  &   1   &  $(M_{12}\cup N_{12})\setminus \{$(0,14),&1 &   2   &  (0,14),(4,10) & 1 &   4   &  (7,4),(12,0)  \\
&&(3,11),(4,10),(5,9),&1 &   5   &  (3,11)&1 &   6  &   (5,9) \\
&&(7,4),(7,7),(8,3),& 1 &   9  &   (10,2)&2 &   6  &   (9,5)\\
&&(8,6),(9,5),(10,0),& 2 &   7  &   (10,4) & 3 &   8   &  (8,3),(10,0)\\
&&(10,2),(10,4),(12,0), & 4 &   9  &   (12,2) &5 &   7  &   (7,7) \\
&&(12,2)$\}$& 9 &   10  &   (8,6) \\\hline
\end{tabular}
\end{center}

\noindent {\bf Lemma $3.7$} {\em
Let $M_{13}=\{(s,t): s+t\leq 17, s\in \{0,3,6,9,15,17\} , t$ is a non-negative integer$\}$ $\setminus \{(0,15),(3,12),(6,8),(6,9),(9,4),(9,5),(9,7),(15,1)\}$, and
$N_{13}=\{(2,15),(5,11),(7,7),(8,9),\\(10,5),(11,4),(11,6),(12,0),(12,3),(13,0),(13,2),(14,3)\}$.
Let $G$ be a kite and $(s,t)\in M_{13}\cup N_{13}$.
Then there is a pair of $(K_{13}\setminus K_{5},G)$-designs with the same vertex set
and the same subgraph $K_{5}$ removed, which intersect in $s$ blocks and
$t+s$ triangles.}

\proof Take the vertex set $X=\{0,1,\ldots,12\}$. Let
\begin{center}\small\tabcolsep 0.07in
 \begin{tabular}{llllllll}
${\cal A}_1:$ &$[0,5,6-12]$,&$[0,7,8-11]$,&$[0,9,10-8]$,&$[0,11,12-4]$,&$[1,5,7-4]$,&$[1,6,8-4]$,\\&
$[1,9,11-4]$,&$[1,10,12-5]$,&$[2,5,8-9]$,&$[2,6,7-9]$,&$[2,12,9-6]$,&$[2,10,11-6]$,\\\vspace{0.1cm}&
$[3,5,9-4]$,&$[3,10,6-4]$,&$[3,7,11-5]$,&$[3,8,12-7]$,&$[4,5,10-7]$;\\

${\cal A}_2:$ &$[0,5,6-4]$,&$[0,7,8-9]$,&$[0,9,10-7]$,&$[0,11,12-7]$,&$[1,5,7-9]$,&$[1,6,8-11]$,\\&
$[1,9,11-5]$,&$[1,10,12-4]$,&$[2,5,8-4]$,&$[2,6,7-4]$,&$[2,12,9-4]$,&$[2,10,11-4]$,\\\vspace{0.1cm}&
$[3,5,9-6]$,&$[3,10,6-12]$,&$[3,7,11-6]$,&$[3,8,12-5]$,&$[4,5,10-8]$;\\

${\cal A}_3:$ &$[0,5,6-4]$,&$[0,7,8-9]$,&$[0,9,10-8]$,&$[0,11,12-7]$,&$[1,5,7-4]$,&$[1,6,8-11]$,\\&
$[1,9,11-4]$,&$[1,10,12-4]$,&$[2,5,8-4]$,&$[2,6,7-9]$,&$[2,12,9-6]$,&$[2,10,11-6]$,\\\vspace{0.1cm}&
$[3,5,9-4]$,&$[3,10,6-12]$,&$[3,7,11-5]$,&$[3,8,12-5]$,&$[4,5,10-7]$;\\

${\cal A}_{4}:$ &$[7,5,1-10]$,&$[9,3,5-12]$,&$[10,6,3-12]$,&$[10,5,4-12]$,&$[4,7,9-1]$,&$[6,7,12-1]$,\\\vspace{0.1cm}&
$[12,10,11-8]$,&$[8,12,9-10]$;\\

${\cal A}_5:$ &$[6,5,0-10]$,&$[0,7,8-3]$,&$[11,9,0-12]$,&$[1,6,8-4]$,&$[2,5,8-10]$,&$[9,6,2-11]$,\\\vspace{0.1cm}&
$[10,7,2-12]$,&$[3,7,11-1]$,&$[4,6,11-5]$;\\

${\cal A}_6:$ &$[6,5,0-12]$,&$[0,7,8-3]$,&$[11,9,0-10]$,&$[1,6,8-4]$,&$[2,5,8-10]$,&$[9,6,2-12]$,\\\vspace{0.1cm}&
$[10,7,2-11]$,&$[3,7,11-5]$,&$[4,6,11-1]$;\\

${\cal A}_7:$ &$[6,5,0-12]$,&$[0,7,8-10]$,&$[11,9,0-10]$,&$[1,6,8-3]$,&$[2,5,8-4]$,&$[9,6,2-11]$,\\\vspace{0.1cm}&
$[10,7,2-12]$,&$[3,7,11-1]$,&$[4,6,11-5]$;\\

${\cal A}_8:$ &$[6,5,0-10]$,&$[0,7,8-10]$,&$[11,9,0-12]$,&$[1,6,8-3]$,&$[2,5,8-4]$,&$[9,6,2-11]$,\\\vspace{0.1cm}&
$[10,7,2-12]$,&$[3,7,11-1]$,&$[4,6,11-5]$;\\

${\cal A}_9:$ &$[6,5,0-12]$,&$[0,7,8-3]$,&$[11,9,0-10]$,&$[1,6,8-4]$,&$[2,5,8-10]$,&$[9,6,2-12]$,\\\vspace{0.1cm}&
$[10,7,2-11]$,&$[3,7,11-1]$,&$[4,6,11-5]$;\\

${\cal A}_{10}:$ &$[6,5,0-12]$,&$[0,7,8-3]$,&$[11,9,0-10]$,&$[1,6,8-4]$,&$[2,5,8-10]$,&$[9,6,2-11]$,\\\vspace{0.1cm}&
$[10,7,2-12]$,&$[3,7,11-1]$,&$[4,6,11-5]$;\\\vspace{0.1cm}

${\cal A}_{11}:$ &$[4,12,5-6]$,&$[6,12,7-8]$,&$[8,12,9-10]$,&$[10,12,11-4]$;\\\vspace{0.1cm}

${\cal A}_{12}:$ &$[6,12,5-4]$,&$[8,12,7-6]$,&$[10,12,9-8]$,&$[4,12,11-10]$;\\

${\cal A}_{13}:$ &$[7,5,1-10]$,&$[9,3,5-12]$,&$[10,6,3-12]$,&$[10,5,4-12]$,&$[4,7,9-1]$,&$[6,7,12-1]$,\\\vspace{0.1cm}&
$[12,8,11-10]$,&$[10,12,9-8]$;\\

${\cal A}_{14}:$ &$[0,5,7-4]$,&$[8,0,6-4]$,&$[11,9,0-12]$,&$[8,5,1-11]$,&$[9,6,1-12]$,&$[1,7,10-0]$,\\&
$[9,5,2-8]$,&$[10,6,2-12]$,&$[2,7,11-5]$,&$[10,5,3-8]$,&$[11,6,3-12]$,&$[3,7,9-4]$,\\\vspace{0.1cm}&
$[4,10,8-11]$.
\end{tabular}
\end{center}
\noindent Let ${\cal B}_{i}={\cal A}_{i}$ when $i=1,2,3$. Let ${\cal B}_{4}={\cal A}_{4}\cup {\cal A}_{5}$.
Let ${\cal B}_{i}={\cal A}_{i}\cup {\cal A}_{13}$ when $5\leq i\leq 10$.
Let ${\cal B}_{i}={\cal A}_{i}\cup {\cal A}_{14}$ when $i=11,12$.
Then $(X,{\cal B}_{i})$ is a $(K_{13}\setminus K_{5},G)$-design for each $1\leq i\leq 12$,
where the removed subgraph $K_{5}$ is constructed on $Y=\{0,1,2,3,4\}$. Consider the following permutations on $X$.

\begin{center}\small\tabcolsep 0.02in
\begin{tabular}{lll}
$\pi_{0,0}=(0\ 1)(2\ 3\ 4)$,&
$\pi_{0,1}=(7\ 8\ 9)(10\ 12\ 11)$,&
$\pi_{0,2}=(7\ 8)(9\ 10\ 12)$,\\
$\pi_{0,3}=(7\ 8\ 9\ 10\ 11\ 12)$,&
$\pi_{0,4}=(7\ 8)(9\ 11\ 10\ 12)$,&
$\pi_{0,5}=(7\ 8)(9\ 10\ 11\ 12)$,\\
$\pi_{0,6}=(7\ 8\ 11\ 12)(9\ 10)$,&
$\pi_{0,7}=(6\ 7\ 8)(9\ 10)(11\ 12)$,&
$\pi_{0,8}=(7\ 8)(9\ 12)(10\ 11)$,\\
$\pi_{0,9}=(5\ 10)(7\ 12)(8\ 11)$,&
$\pi_{0,10}=(7\ 8)(9\ 10)(11\ 12)$,&
$\pi_{0,11}=(0\ 1)(2\ 3\ 4)(5\ 8)(9\ 10\ 11)$,\\
$\pi_{0,12}=(3\ 4)(9\ 12)(10\ 11)$,&
$\pi_{0,13}=(0\ 1\ 2)(3\ 4)(5\ 7\ 6)(10\ 11\ 12)$,&
$\pi_{0,14}=(3\ 4)(5\ 8)(6\ 7)(9\ 11)(10\ 12)$,\\
$\pi_{0,16}=(5\ 7)(6\ 8)(9\ 11)(10\ 12)$,&
$\pi_{0,17}=(1)$,&
$\pi_{2,15}=(1\ 2)(5\ 10)(6\ 9)(7\ 11)(8\ 12)$,\\
$\pi_{3,0}=(2\ 4)(7\ 9)(8\ 10\ 11)$,&
$\pi_{3,1}=(1\ 2)(3\ 4)$,&
$\pi_{3,2}=(0\ 3)(8\ 12\ 11)$,\\
$\pi_{3,3}=(0\ 3)(6\ 8\ 11)$,&
$\pi_{3,4}=(0\ 4\ 3\ 1)(5\ 11\ 7)(8\ 10)$,&
$\pi_{3,5}=(0\ 2\ 4)$,\\
$\pi_{3,6}=(5\ 10)(6\ 7\ 12\ 9)(8\ 11)$,&
$\pi_{3,7}=(5\ 6)(7\ 8)(9\ 10)$,&
$\pi_{3,8}=(0\ 1)(5\ 10)(6\ 12\ 7\ 9)(8\ 11)$,\\
$\pi_{3,9}=(5\ 9\ 6\ 10)(7\ 11\ 8\ 12)$,&
$\pi_{3,10}=(0\ 2)(6\ 8)(9\ 11)$,&
$\pi_{3,11}=(1\ 2)(3\ 4)(5\ 6)(11\ 12)$,\\
$\pi_{3,13}=(1\ 2)(5\ 6)(9\ 10)$,&
$\pi_{3,14}=(1\ 2)(7\ 8)(11\ 12)$,&
$\pi_{5,11}=(0\ 2\ 1)(6\ 8\ 7)(10\ 12\ 11)$,\\
$\pi_{6,0}=(6\ 10\ 12)$,&
$\pi_{6,1}=(0\ 2)(8\ 12)$,&
$\pi_{6,2}=(5\ 11)$,\\
$\pi_{6,3}=(5\ 10)(6\ 9)(8\ 11)$,&
$\pi_{6,4}=(5\ 9)$,&
$\pi_{6,5}=(6\ 7)$,\\
$\pi_{6,6}=(2\ 4)$,&
$\pi_{6,7}=(1\ 2)(3\ 4)(5\ 6)(9\ 11)$,&
$\pi_{6,10}=(0\ 2)(5\ 10\ 7\ 12)(6\ 11\ 8\ 9)$,\\
$\pi_{6,11}=(5\ 10)(6\ 9)(7\ 12)(8\ 11)$,&
$\pi_{7,7}=(0\ 1)(2\ 3)(6\ 7)(9\ 11)$,&
$\pi_{8,9}=(1)$,\\
$\pi_{9,0}=(0\ 2)$,&
$\pi_{9,1}=(9\ 11)$,&
$\pi_{9,2}=(3\ 4)$,\\
$\pi_{9,3}=(1\ 4)$,&
$\pi_{9,6}=\pi_{9,8}=\pi_{10,5}=(1)$,&
$\pi_{11,4}=\pi_{11,6}=(1)$\\
$\pi_{12,0}=(9\ 11)$,&
$\pi_{12,3}=\pi_{13,0}=(1)$,&
$\pi_{13,2}=\pi_{14,3}=\pi_{15,0}=(1)$,\\
$\pi_{15,2}=\pi_{17,0}=(1)$.
\end{tabular}
\end{center}
\noindent We have that for each row in Table II,
 $\pi_{s,t}Y=Y$, $|{\cal B}_i\cap \pi_{s,t}{\cal B}_j|=s$ and
 $|T({\cal B}_i\setminus\pi_{s,t}{\cal B}_j)\cap T(\pi_{s,t}{\cal B}_j\setminus{\cal B}_i)|=t$. \qed

\begin{center}
\small  Table II.  \ Fine triangle intersections for $(K_{13}\setminus K_{5},G)$-designs
  \begin{tabular}{ccc|ccc|ccc}\hline
$i$ &  $j$ &  $(s,t)$ & $i$ &  $j$ &  $(s,t)$ & $i$ &  $j$ &  $(s,t)$\\\hline
1  &   1   &  $(M_{13}\cup N_{13})\setminus \{$(0,11),&  1  &   2   &   (0,17) & 1&      3&  (3,14),(5,11),(6,10),(9,8)      \\
&& (0,17),(3,14),(5,11), &  2&3&(8,9) & 4& 4&(0,11)  \\
&& (6,7),(6,10),(7,7),  &4&5& (6,7),(7,7),(12,0),(15,0)&  4&6&   (9,6) \\
&& (8,9),(9,1),(9,2),  & 4&7&  (10,5)&4&  8&  (12,3)   \\
&& (9,6),(9,8),(10,5),  & 4&  9&  (11,4)&4& 10&  (13,2)   \\
&& (11,4),(11,6),(12,0),& 5&      6&    (11,6) &  5&   8&   (14,3)  \\
&&  (12,3),(13,0),(13,2),  & 6&      6&   (9,1),(9,2)&7&   8&   (15,2)     \\
&&  (14,3),(15,0),$ (15,2)\}$ &  11&  12&  (13,0)    \\
\hline
\end{tabular}
\end{center}

\noindent {\bf Lemma $3.8$} {\em
Let $M_{14}=\{(s,t): s+t\leq 19, s\in \{0,4,8,17,19\}, t$ is a non-negative integer$\}$ $\setminus \{(4,13),(8,8),(8,9),(8,10),(17,1)\}$,
and $N_{14}=\{(7,9),(7,11),(10,5),(10,9),(11,6),(12,0),(12,\\1),(12,2),(13,6),(14,0),(14,3),(15,2),(16,3)\}$.
Let $G$ be a kite and $(s,t)\in M_{14}\cup N_{14}$.
Then there is a pair of  $(K_{14}\setminus K_{6},G)$-designs with the same vertex set
and the same subgraph $K_{6}$ removed, which intersect in $s$ blocks and
$t+s$ triangles.}

\proof Take the vertex set $X=\{0,1,\ldots,13\}$. Let
\begin{center}
\small
 \begin{tabular}{lllllll}
${\cal A}_1:$ &$[0,10,11-5]$,&$[0,12,13-5]$,&$[1,8,6-5]$,&$[1,10,12-4]$,&$[1,11,13-6]$,&$[2,9,6-12]$,\\\vspace{0.1cm}&
$[2,10,13-7]$,&$[2,11,12-5]$,&$[3,7,11-8]$,&$[3,8,12-9]$,&$[4,6,11-9]$;\\
\end{tabular}
\end{center}

\begin{center}
\small
 \begin{tabular}{lllllll}
${\cal A}_{2}:$ &$[0,6,7-12]$,&$[0,8,9-5]$,&$[1,9,7-5]$,&$[2,7,8-5]$,&$[3,6,10-5]$,&$[3,13,9-4]$,\\\vspace{0.1cm}&
$[4,7,10-9]$,&$[4,13,8-10]$;\\

${\cal A}_{3}:$ &$[0,10,11-9]$,&$[0,12,13-7]$,&$[1,8,6-12]$,&$[1,10,12-9]$,&$[1,11,13-5]$,&$[2,9,6-5]$,\\\vspace{0.1cm}&
$[2,10,13-6]$,&$[2,11,12-4]$,&$[3,7,11-5]$,&$[3,8,12-5]$,&$[4,6,11-8]$;\\

${\cal A}_{4}:$ &$[0,10,12-4]$,&$[3,7,11-0]$,&$[3,8,12-1]$,&$[4,6,11-1]$,&$[5,6,12-2]$,&$[5,9,11-2]$,\\\vspace{0.1cm}&
$[11,13,12-9]$,&$[9,13,10-11]$;\\

${\cal A}_{5}:$ &$[0,10,12-4]$,&$[3,7,11-0]$,&$[3,8,12-1]$,&$[4,6,11-1]$,&$[5,6,12-2]$,&$[5,9,11-2]$,\\\vspace{0.1cm}&
$[9,13,12-11]$,&$[11,13,10-9]$;\\

${\cal A}_{6}:$ &$[0,10,12-1]$,&$[3,7,11-0]$,&$[3,8,12-4]$,&$[4,6,11-1]$,&$[5,6,12-2]$,&$[5,9,11-2]$,\\\vspace{0.1cm}&
$[9,13,12-11]$,&$[11,13,10-9]$;\\

${\cal A}_{7}:$ &$[0,10,12-4]$,&$[3,7,11-2]$,&$[3,8,12-1]$,&$[4,6,11-0]$,&$[5,6,12-2]$,&$[5,9,11-1]$,\\\vspace{0.1cm}&
$[9,13,12-11]$,&$[11,13,10-9]$;\\

${\cal A}_{8}:$ &$[0,10,12-2]$,&$[3,7,11-2]$,&$[3,8,12-4]$,&$[4,6,11-0]$,&$[5,6,12-1]$,&$[5,9,11-1]$,\\\vspace{0.1cm}&
$[9,13,12-11]$,&$[11,13,10-9]$;\\\vspace{0.1cm}

${\cal A}_{9}:$ &$[12,13,0-11]$,&$[11,13,1-10]$,&$[10,13,2-9]$,&$[9,13,3-8]$,&$[8,13,4-12]$;\\\vspace{0.1cm}

${\cal A}_{10}:$ &$[11,13,0-12]$,&$[10,13,1-11]$,&$[9,13,2-10]$,&$[8,13,3-9]$,&$[12,13,4-8]$;\\

${\cal A}_{11}:$ &$[0,6,7-5]$,&$[0,8,9-4]$,&$[1,9,7-12]$,&$[2,7,8-5]$,&$[3,6,10-5]$,&$[3,13,9-5]$,\\\vspace{0.1cm}&
$[4,7,10-9]$,&$[4,13,8-10]$;\\

${\cal A}_{12}:$ &$[0,6,7-12]$,&$[0,8,9-4]$,&$[1,8,6-13]$,&$[1,7,9-3]$,&$[9,6,2-13]$,&$[2,7,8-11]$,\\\vspace{0.1cm}&
$[6,10,3-13]$,&$[4,7,10-2]$,&$[4,8,13-1]$,&$[5,7,13-0]$,&$[5,8,10-1]$;\\

${\cal A}_{13}:$ &$[0,6,7-12]$,&$[0,8,9-12]$,&$[1,6,8-12]$,&$[1,7,9-10]$,&$[2,6,11-8]$,&$[8,7,2-12]$,\\&
$[3,6,10-0]$,&$[11,7,3-12]$,&$[6,9,4-11]$,&$[4,7,10-11]$,&$[5,6,12-1]$,&$[5,7,13-6]$,\\\vspace{0.1cm}&
$[5,8,10-12]$,&$[5,9,11-12]$.\\

\end{tabular}
\end{center}
\noindent Let ${\cal B}_{i}={\cal A}_{i}\cup {\cal A}_{11}$ when $i=1,3$. Let ${\cal B}_{2}={\cal A}_{2}\cup {\cal A}_{3}$.
Let ${\cal B}_{i}={\cal A}_{i}\cup {\cal A}_{12}$ when $4\leq i\leq 8$.
Let ${\cal B}_{i}={\cal A}_{i}\cup {\cal A}_{13}$ when $i=9,10$.
Then $(X,{\cal B}_{i})$ is a $(K_{14}\setminus K_{6},G)$-design for each $1\leq i\leq 10$,
where the removed subgraph $K_{6}$ is constructed on $Y=\{0,1,2,3,4,5\}$. Consider the following permutations
on $X$.

\begin{center}
\small
\begin{tabular}{lll}\tabcolsep 0.02in
$\pi_{0,0}=(0\ 1)(2\ 3)(4\ 5)$,&
$\pi_{0,1}=(7\ 8)(9\ 10\ 11\ 13)$,&
$\pi_{0,2}=(7\ 8)(9\ 10)(11\ 12\ 13)$,\\
$\pi_{0,3}=(7\ 8)(9\ 10)(12\ 13)$,&
$\pi_{0,4}=(1\ 2)(3\ 4\ 5)$,&
$\pi_{0,5}=(7\ 8)(9\ 12)(10\ 13)$,\\
$\pi_{0,6}=(7\ 8)(10\ 11\ 12\ 13)$,&
$\pi_{0,7}=(7\ 8)(9\ 10\ 13)(11\ 12)$,&
$\pi_{0,8}=(7\ 8)(10\ 11)(12\ 13)$,\\
$\pi_{0,9}=(6\ 7)(8\ 9\ 10\ 11)(12\ 13)$,&
$\pi_{0,10}=(1\ 2\ 5)(8\ 9)$&
$\pi_{0,11}=(7\ 10)(8\ 13)(9\ 12)$,\\
$\pi_{0,12}=(7\ 8)(10\ 13)(11\ 12)$,&
$\pi_{0,13}=(1\ 2)(3\ 4\ 5)(6\ 7)(10\ 12\ 13)$,&
$\pi_{0,14}=(1\ 5\ 2)(8\ 9)(12\ 13)$,\\
$\pi_{0,15}=(3\ 4\ 5)(10\ 11)(12\ 13)$,&
$\pi_{0,16}=(1\ 2)(4\ 5)(8\ 9)(12\ 13)$,&
$\pi_{0,17}=(1\ 2)(3\ 4)(6\ 7)(12\ 13)$,\\
$\pi_{0,18}=(6\ 9)(7\ 8)(10\ 13)(11\ 12)$,&
$\pi_{0,19}=(6\ 11)(7\ 10)(8\ 13)(9\ 12)$,&
$\pi_{4,0}=(0\ 1)(3\ 4)$,\\
$\pi_{4,1}=(9\ 10\ 12\ 13)$,&
$\pi_{4,2}=(0\ 3)(12\ 13)$,&
$\pi_{4,3}=(1\ 2)(3\ 5)$,\\
$\pi_{4,4}=(0\ 1)(4\ 5)$,&
$\pi_{4,5}=(1\ 2)(6\ 13\ 7)(10\ 11)$,&
$\pi_{4,6}=(6\ 13)(8\ 11)(9\ 10)$,\\
$\pi_{4,7}=(6\ 11)(7\ 10)$,&
$\pi_{4,8}=(3\ 4\ 5)$,&
$\pi_{4,9}=(0\ 2)(7\ 9)(11\ 13)$,\\
$\pi_{4,10}=(0\ 4)(7\ 11)(8\ 12)$,&
$\pi_{4,11}=(0\ 3)(7\ 10)(9\ 12)$,&
$\pi_{4,12}=(0\ 2\ 1)(6\ 11\ 8\ 10\ 7\ 12)$\\
$\pi_{4,14}=(6\ 10)(7\ 11)(8\ 12)(9\ 13)$,&
$\pi_{4,15}=(1)$,&
\hspace{1.2cm}$(9\ 13)$,\\
$\pi_{7,9}=(1\ 2)(4\ 5)(6\ 7)(10\ 11)$,&
$\pi_{7,11}=(1\ 3)(2\ 4)(6\ 8\ 12\ 10)$&
$\pi_{8,0}=(0\ 1\ 4)$,\\
$\pi_{8,1}=(8\ 10)$,&
\hspace{1.2cm}$(7\ 9\ 13\ 11)$,&
$\pi_{8,2}=(7\ 12)$,\\
$\pi_{8,3}=(6\ 8)$,&
$\pi_{8,4}=(7\ 9)$,&
$\pi_{8,5}=(12\ 13)$,\\
$\pi_{8,6}=(0\ 1)$,&
$\pi_{8,7}=(1\ 5)$,&
$\pi_{8,11}=(1)$,\\
$\pi_{10,5}=(0\ 5)$,&
$\pi_{10,9}=(1\ 2)(6\ 7)(10\ 11)$,&
$\pi_{11,6}=(1)$,\\
$\pi_{12,0}=(2\ 5)$,&
$\pi_{12,1}=(2\ 5)$,&
$\pi_{12,2}=(9\ 11)$,\\
$\pi_{13,6}=\pi_{14,0}=\pi_{14,3}=(1)$£¬&
$\pi_{15,2}=\pi_{16,3}=\pi_{17,0}=(1)$,&
$\pi_{17,2}=\pi_{19,0}=(1)$.
\end{tabular}
\end{center}

\noindent We have that for each row in Table III,
 $\pi_{s,t}Y=Y$, $|{\cal B}_i\cap \pi_{s,t}{\cal B}_j|=s$ and
 $|T({\cal B}_i\setminus\pi_{s,t}{\cal B}_j)\cap T(\pi_{s,t}{\cal B}_j\setminus{\cal B}_i)|=t$. \qed

\begin{center}
\small  Table III.  \ Fine triangle intersections for $(K_{14}\setminus K_{6},G)$-designs
  \begin{tabular}{ccc|ccc|ccc}\hline
$i$ &  $j$ &  $(s,t)$ & $i$ &  $j$ &  $(s,t)$&$i$ &  $j$ &  $(s,t)$ \\\hline
1  &   1   &  $(M_{14}\cup N_{14})\setminus \{$(0,13), & 1& 2&  (4,15) & 1&  3&  (8,11)     \\
&&(0,17),(4,9),(4,15),& 4&  4&  (0,17),(4,9),(8,5),(8,6),(12,1) & 4&  5&  (0,13),(12,2),(17,0)   \\
&&(8,5),(8,6),(8,11),& 4& 6&   (15,2)& 4&  7&   (14,3)    \\
&&(11,6),(12,0),(12,1),& 4& 8&   (11,6)   &  5&      6&   (17,2)  \\
&&(12,2),(13,6),(14,0),& 5&  7&   (16,3)    &  5&      8&    (13,6)  \\
&&(14,3),(15,2),(16,3), &  9&      9&    (12,0)   & 9&      10&   (14,0)    \\
&&(17,0),(17,2)$\}$ &&&&&&\\
\hline
\end{tabular}
\end{center}

\noindent {\bf Lemma $3.9$} {\em
Let $M_{15}=\{(0,0),(0,1),\ldots,(0,16),(0,18),(0,21),(21,0)\}$
and $N_{15}=\{(4,17),(6,0),\\(6,1),\ldots,(6,8),(6,12),(8,13),(9,8),(12,0),(12,4),(12,6),(12,9),
(13,2),(13,8),(17,1),(17,4),\\(18,0)\}$.
Let $G$ be a kite and $(s,t)\in M_{15}\cup N_{15}$.
Then there is a pair of $(K_{15}\setminus K_{7},G)$-designs with the same vertex set
and the same subgraph $K_{7}$ removed, which intersect in $s$ blocks and
$t+s$ triangles.}

\proof Take the vertex set $X=\{0,1,\cdots,14\}$. Let
\begin{center}\small
 \begin{tabular}{lllllll}
${\cal A}_1:$ &$[6 ,14 ,7- 5]$,&$[8, 1, 9 -4]$,&$[10 ,2 ,9 -6]$,&$[10 ,3, 11- 8]$,&$[13 ,6, 12- 5]$,&$[13, 5 ,14 -1]$,\\\vspace{0.1cm}&
$[9 ,5 ,11- 6]$,&$[10 ,12, 7 -4]$,&$[8 ,3, 7- 0]$,&$[11 ,1 ,7- 2]$;\\

${\cal A}_2:$ &$[6 ,14 ,7 -4]$,&$[8 ,1, 9- 6]$,&$[10, 2 ,9 -4]$,&$[11 ,4, 12 -2]$,&$[13 ,6 ,12- 1]$,&$[13, 5 ,14- 3]$,\\\vspace{0.1cm}&
$[9 ,3 ,12- 5]$,&$[10 ,12, 7- 5]$,&$[8 ,3, 7- 2]$,&$[11 ,1, 7 -0]$;\\

${\cal A}_3:$ &$[5 ,14, 7 -6]$,&$[8 ,1, 9- 4]$,&$[10, 2 ,9 -6]$,&$[10 ,3 ,11- 8]$,&$[13 ,5, 12- 6]$,&$[13 ,6 ,14- 1]$,\\\vspace{0.1cm}&
$[9 ,5 ,11- 6]$,&$[10 ,12 ,7 -4]$,&$[8 ,3 ,7- 0]$,&$[11 ,1, 7- 2 ]$;\\

${\cal A}_4:$ &$[6, 14, 7 -4]$,&$[8 ,1 ,9 -4]$,&$[10 ,2 ,9- 6]$,&$[10 ,3 ,11- 8]$,&$[13 ,6 ,12- 5]$,&$[13 ,5 ,14- 1]$,\\\vspace{0.1cm}&
$[9 ,5 ,11 -6]$,&$[10 ,12, 7- 5]$,&$[8 ,3 ,7 -2]$,&$[11 ,1, 7- 0]$;\\

${\cal A}_5:$ &$[5 ,14, 7 -6]$,&$[8 ,1, 9- 6]$,&$[10, 2 ,9 -4]$,&$[10 ,3 ,11- 6]$,&$[13 ,5, 12- 6]$,&$[13 ,6 ,14- 1]$,\\\vspace{0.1cm}&
$[9 ,5 ,11- 8]$,&$[10 ,12 ,7 -4]$,&$[8 ,3 ,7- 2]$,&$[11 ,1, 7- 0 ]$;\\

${\cal A}_6:$ &$[6, 14, 7 -4]$,&$[8 ,1 ,9 -6]$,&$[10 ,2 ,9- 4]$,&$[10 ,3 ,11- 6]$,&$[13 ,6 ,12- 5]$,&$[13 ,5 ,14- 1]$,\\\vspace{0.1cm}&
$[9 ,5 ,11 -8]$,&$[10 ,12, 7- 5]$,&$[8 ,3 ,7 -2]$,&$[11 ,1, 7- 0]$;\\

${\cal A}_7:$ &$[6 ,14 ,7 -4]$,&$[8 ,1, 9- 6]$,&$[10, 2 ,9 -4]$,&$[11 ,4, 12 -2]$,&$[13 ,6 ,12- 1]$,&$[13, 5 ,14- 1]$,\\\vspace{0.1cm}&
$[9 ,3 ,12- 5]$,&$[10 ,12, 7- 5]$,&$[8 ,3, 7- 2]$,&$[11 ,1, 7 -0]$;\\

${\cal A}_8:$ &$[5 ,14, 7- 4]$,&$[8 ,1 ,9- 4]$,&$[10, 2, 9 -6]$,&$[11 ,4 ,12- 1]$,&$[13 ,5 ,12- 6]$,&$[13 ,6, 14- 1]$,\\\vspace{0.1cm}&
$[9 ,3, 12- 2]$,&$[10 ,12, 7 -6]$,&$[8 ,3, 7- 0]$,&$[11 ,1 ,7- 2]$;\\

${\cal A}_9:$ &$[11 ,4 ,12- 1]$,&$[ 0 ,9 ,14- 12]$,&$[2 ,11 ,14 -3]$,&$[13 ,1 ,10 -5]$,&$[12, 0, 8- 14]$,&$[13, 2, 8- 4]$,\\\vspace{0.1cm}&
$[14 ,4, 10- 0]$,&$[9 ,7 ,13- 3]$,&$[10 ,6, 8- 5]$,&$[0 ,11, 13- 4]$,&$[9, 3, 12- 2]$;\\

${\cal A}_{10}:$ &$[10, 3 ,11- 8]$,&$[0 ,9 ,14- 12]$,&$[2 ,11, 14 -3]$,&$[13 ,1 ,10 -5]$,&$[12, 0 ,8 -14]$,&$[13 ,2, 8- 4]$,\\\vspace{0.1cm}&
$[14, 4 ,10- 0]$,&$[9 ,7 ,13- 3]$,&$[9 ,5, 11- 6]$,&$[10 ,6 ,8- 5]$,&$[0 ,11, 13 -4]$;\\

${\cal A}_{11}:$ &$[10 ,3 ,11 -6]$,&$[0 ,9 ,14- 1]$,&$[2 ,11, 14- 12]$,&$[13, 1, 10- 0]$,&$[12 ,0, 8 -5]$,&$[13 ,2, 8- 14]$,\\\vspace{0.1cm}&
$[14 ,4 ,10- 5]$,&$[9 ,7, 13 -4]$,&$[9 ,5, 11 -8]$,&$[10 ,6 ,8 -4]$,&$[0 ,11, 13- 3]$.
\end{tabular}
\end{center}

\noindent Let ${\cal B}_{i}={\cal A}_{i}\cup {\cal A}_{9}$ when $i=1,3,4,5,6$. Let ${\cal B}_{2}={\cal A}_{2}\cup {\cal A}_{11}$. Let
${\cal B}_{i}={\cal A}_{i}\cup {\cal A}_{10}$ when $i=7,8$. Then $(X,{\cal B}_{i})$ is a $(K_{15}\setminus K_{7},G)$-design for each $1\leq i\leq 8$,
where the removed subgraph $K_{7}$ is constructed on $Y=\{0,1,2,3,4,5,6\}$. Consider the following permutations
on $X$.

\begin{center}\small\tabcolsep 0.15in
\begin{tabular}{lll}
$\pi_{0,0}=(8\ 9\ 10)(11\ 12\ 13\ 14)$,&
$\pi_{0,1}=(8\ 9)(11\ 12\ 13\ 14)$,&
$\pi_{0,2}=(0\ 1\ 2\ 3)(4\ 5\ 6)$,\\
$\pi_{0,3}=(8\ 9)(11\ 12)(13\ 14)$,&
$\pi_{0,4}=(0\ 1\ 2\ 3)(4\ 6)$,&
$\pi_{0,5}=(0\ 1\ 2\ 3)(4\ 5)$,\\
$\pi_{0,6}=(8\ 9)(10\ 11)(12\ 14)$,&
$\pi_{0,7}=(0\ 1\ 3)(4\ 6)$,&
$\pi_{0,8}=(0\ 2\ 3)(4\ 5)$,\\
$\pi_{0,9}=(8\ 10\ 12)(9\ 13)(11\ 14)$,&
$\pi_{0,10}=(1\ 4)(3\ 6)(7\ 12)$,&
$\pi_{0,11}=(8\ 10)(9\ 13)(11\ 14)$,\\
$\pi_{0,12}=(0\ 4)(1\ 6)(3\ 5)(9\ 10)$&
$\pi_{0,13}=(1\ 3)(2\ 5\ 4)(8\ 11)$&
$\pi_{0,14}=(1\ 2\ 4)(3\ 6)(7\ 10)$\\
\hspace{1.2cm}$(12\ 13)$,&\hspace{1.2cm}$(9\ 10)(12\ 13)$,&\hspace{1.2cm}$(9\ 13)(11\ 14)$,\\
$\pi_{0,15}=(1\ 4)(3\ 6)(7\ 10)(9\ 13)(11\ 14)$,&
$\pi_{0,16}=(0\ 4)$,&
$\pi_{0,18}=\pi_{0,21}=\pi_{4,17}=(1)$,\\
$\pi_{6,0}=(9\ 10\ 13)$,&
$\pi_{6,1}=(7\ 9)(8\ 10)$,&
$\pi_{6,2}=(1\ 4)(9\ 13)$,\\
$\pi_{6,3}=(2\ 4)(7\ 11)$,&
$\pi_{6,4}=(0\ 4)(13\ 14)$,&
$\pi_{6,5}=(3\ 4)(8\ 12)$,\\
$\pi_{6,6}=(1\ 3\ 6)$,&
$\pi_{6,7}=(0\ 3\ 5)$,&
$\pi_{6,8}=(0\ 5\ 4)$,\\
$\pi_{6,12}=\pi_{8,13}=(1)$,&
$\pi_{9,8}=(4\ 5)$,&
$\pi_{12,0}=(7\ 9)$,\\
$\pi_{12,4}=(0\ 4)$,&
$\pi_{12,6}=\pi_{12,9}=(1)$,&
$\pi_{13,2}=(2\ 3)$,\\
$\pi_{13,8}=\pi_{17,1}=\pi_{17,4}=(1)$,&
$=\pi_{18,0}=\pi_{21,0}=(1)$.
\end{tabular}
\end{center}
\noindent We have that for each $i,j,s,t$ in Table IV,
 $\pi_{s,t}Y=Y$, $|{\cal B}_i\cap \pi_{s,t}{\cal B}_j|=s$ and
 $|T({\cal B}_i\setminus\pi_{s,t}{\cal B}_j)\cap T(\pi_{s,t}{\cal B}_j\setminus{\cal B}_i)|=t$. \qed

\begin{center}
\small  Table IV.  \ Fine triangle intersections for $(K_{15}\setminus K_{7},G)$-designs
  \begin{tabular}{ccc|ccc|ccc}\hline
$i$ &  $j$ &  $(s,t)$ & $i$ &  $j$ &  $(s,t)$&$i$ &  $j$ &  $(s,t)$  \\
\hline     1  &   1   &  $(M_{15}\cup N_{15})\setminus \{(0,16),$  & 1& 2 & (0,16),(0,21) &1 &    3   &   (18,0) \\
&&(0,18),(0,21),(4,17),&1  &   4  &   (17,4) & 1  &   5  &   (12,6)  \\
&&(6,12),(8,13),(12,6),&1  &   6  &   (13,8) & 1  &   7  &   (12,9)   \\
&&(12,9),(13,8),(17,1),&1  &   8  &   (17,1) & 2  &   3  &   (0,18)    \\
&&(17,4)(18,0)\}& 2  &   4  &   (4,17) &2  &   5  &   (6,12)    \\
&&& 2  &   6  &   (8,13)    \\
\hline
\end{tabular}
\end{center}

\noindent {\bf Lemma $3.10$} {\em
Let $M_{18}=\{(0,0),(0,1),\ldots,(0,17),(0,27),(27,0)\}$ and $N_{18}=\{(3,24),(6,21),(8,0),\\(9,18),(12,5),(12,15),(15,0),(15,12),(18,9),(21,0),(24,3)\}$.
Let $G$ be a kite and $(s,t)\in M_{18}\cup N_{18}$.
Then there is a pair of $(K_{18}\setminus K_{10},G)$-designs with the same vertex set
and the same subgraph $K_{10}$ removed, which intersect in $s$ blocks and
$t+s$ triangles.}

\proof Take the vertex set $X=\{0,1,\ldots,17\}$. Let
\begin{center}
\small
 \begin{tabular}{llllll}
${\cal A}_1:$ &$[16, 7, 15 -1]$,&$[2 ,14, 12- 7]$,&$[1, 10, 12 -0]$,&$[17 ,2 ,13- 0]$,&$[16 ,1, 13- 9]$,\\
&$[17 ,7, 14- 0]$,&$[8 ,12, 15- 0]$,&$[8 ,11, 14 -1]$, &$[7 ,10, 13- 8]$,&$[17 ,9 ,12- 6]$,\\\vspace{0.1cm}&$[3 ,10, 14- 9]$,&$[4, 10 ,15 -2]$;\\

${\cal A}_2:$ &$[16 ,7 ,15 -0]$, &$[13 ,12 ,4 -11]$,&$[12 ,11, 3- 17]$,&$[1 ,17 ,11- 7]$,&$[16 ,14, 4- 17]$,\\
&$[15, 13, 3 -16]$,&$[1, 10, 12- 7]$, &$[16, 1, 13- 0]$,&$[17 ,7, 14- 9]$,&$[7 ,10 ,13- 9]$,\\\vspace{0.1cm}&$[3 ,10 ,14- 1]$,&$[4 ,10 ,15- 1]$;\\

${\cal A}_3:$ &$[16, 7, 15 -0]$,&$[2 ,14, 12- 7]$,&$[1, 10, 12 -0]$,&$[17 ,2 ,13- 0]$,&$[16 ,1, 13- 9]$,\\
&$[17 ,7, 14- 0]$,&$[8 ,12, 15- 2]$,&$[8 ,11, 14 -1]$, &$[7 ,10, 13- 8]$,&$[17 ,9 ,12- 6]$,\\\vspace{0.1cm}&$[3 ,10, 14- 9]$,&$[4, 10 ,15 -1]$;\\

${\cal A}_4:$ &$[16, 7, 15 -0]$,&$[2 ,14, 12- 7]$,&$[1, 10, 12 -0]$,&$[17 ,2 ,13- 0]$,&$[16 ,1, 13- 9]$,\\
&$[17 ,7, 14- 9]$,&$[8 ,12, 15- 2]$,&$[8 ,11, 14 -0]$, &$[7 ,10, 13- 8]$,&$[17 ,9 ,12- 6]$,\\\vspace{0.1cm}&$[3 ,10, 14- 1]$,&$[4, 10 ,15 -1]$;\\

${\cal A}_5:$ &$[16, 7, 15 -0]$,&$[2 ,14, 12- 7]$,&$[1, 10, 12 -0]$,&$[17 ,2 ,13- 8]$,&$[16 ,1, 13- 0]$,\\
&$[17 ,7, 14- 9]$,&$[8 ,12, 15- 2]$,&$[8 ,11, 14 -0]$, &$[7 ,10, 13- 9]$,&$[17 ,9 ,12- 6]$,\\\vspace{0.1cm}&$[3 ,10, 14- 1]$,&$[4, 10 ,15 -1]$;\\

${\cal A}_6:$ &$[16, 7, 15 -0]$,&$[2 ,14, 12- 6]$,&$[1, 10, 12 -7]$,&$[17 ,2 ,13- 8]$,&$[16 ,1, 13- 0]$,\\
&$[17 ,7, 14- 9]$,&$[8 ,12, 15- 2]$,&$[8 ,11, 14 -0]$, &$[7 ,10, 13- 9]$,&$[17 ,9 ,12- 0]$,\\\vspace{0.1cm}&$[3 ,10, 14- 1]$,&$[4, 10 ,15 -1]$;\\
\end{tabular}
\end{center}

\begin{center}
\small
 \begin{tabular}{llllll}
${\cal A}_7:$ &$[16 ,1 ,15 -0]$, &$[13 ,12 ,7 -11]$,&$[12 ,11, 4- 17]$,&$[3 ,17 ,11- 1]$,&$[16 ,14, 7- 17]$,\\
&$[15, 13, 4 -16]$,&$[3, 10, 12- 1]$, &$[16, 3, 13- 0]$,&$[17 ,1, 14- 9]$,&$[1 ,10 ,13- 9]$,\\\vspace{0.1cm}&$[4 ,10 ,14- 3]$,&$[7 ,10 ,15- 3]$;\\

${\cal A}_8:$ &$[16 ,7 ,15 -0]$, &$[13 ,12 ,3 -11]$,&$[12 ,11, 4- 17]$,&$[1 ,17 ,11- 7]$,&$[16 ,14, 3- 17]$,\\
&$[15, 13, 4 -16]$,&$[1, 10, 12- 7]$, &$[16, 1, 13- 0]$,&$[17 ,7, 14- 9]$,&$[7 ,10 ,13- 9]$,\\\vspace{0.1cm}&$[4 ,10 ,14- 1]$,&$[3 ,10 ,15- 1]$;\\

${\cal A}_9:$ &$[0 ,17, 10 -5]$,&$[8 ,17, 16- 2]$, &$[15 ,14 ,6- 16]$,&$[14 ,13, 5- 16]$,&$[13, 12, 4 -17]$,\\
&$[12 ,11, 3- 16]$,&$[11, 2, 10- 6]$,&$[ 16 ,9 ,10 -8]$, &$[ 1, 17 ,11 -5]$,&$[0 ,11 ,16 -12]$,\\\vspace{0.1cm}
&$[17 ,15, 5 -12]$,&$[ 16 ,14 ,4 -11]$,
&$[ 15 ,13, 3- 17]$,&$[13 ,11, 6 -17]$, &$[15, 9, 11 -7]$;\\

${\cal A}_{10}:$ &$[0 ,17 ,10- 8]$, &$[8, 17 ,16 -12]$,&$[15 ,14 ,6- 17]$,&$[14, 13 ,5 -12]$,&$[11, 2, 10- 5]$,\\&$[16 ,9, 10- 6]$,
&$[ 0 ,11 ,16- 2]$, &$[17 ,15, 5 -16]$,&$[2 ,14 ,12- 6]$,&$[13 ,11, 6 -16]$,\\&$[15, 9 ,11- 5]$,&$[17, 2, 13- 8]$,&$[8 ,12, 15- 2]$, &$[8 ,11 ,14 -0]$,&$[17, 9 ,12- 0]$.\\
\end{tabular}
\end{center}
\noindent Let ${\cal B}_{i}={\cal A}_{i}\cup {\cal A}_{9}$ when $i=1,3,4,5,6$. Let ${\cal B}_{i}={\cal A}_{i}\cup {\cal A}_{10}$ when $i=2,7,8$.
Then $(X,{\cal B}_{i})$ is a $(K_{18}\setminus K_{9},G)$-design for each $1\leq i\leq 8$,
where the removed subgraph $K_{10}$ is constructed on $Y=\{0,1,\ldots,9\}$. Consider the following permutations on $X$.

\begin{center}
\small
\begin{tabular}{ll}
$\pi_{0,0}=(11\ 12\ 13\ 14)(15\ 16\ 17)$,&
$\pi_{0,1}=(11\ 12)(13\ 14\ 17\ 16\ 15)$,\\
$\pi_{0,2}=(11\ 12)(13\ 14\ 15\ 16\ 17)$,,&
$\pi_{0,3}=(11\ 12)(13\ 14)(15\ 16\ 17)$,\\
$\pi_{0,4}=(11\ 12)(13\ 14\ 15)(16\ 17)$,&
$\pi_{0,5}=(11\ 12)(13\ 14\ 16)(15\ 17)$,\\
$\pi_{0,6}=(11\ 12)(13\ 16\ 17)(14\ 15)$,&
$\pi_{0,7}=(10\ 11)(12\ 14\ 17)(15\ 16)$,\\
$\pi_{0,8}=(10\ 11)(14\ 15)(16\ 17)$,&
$\pi_{0,9}=(10\ 11\ 16\ 15)(12\ 17)(13\ 14)$,\\
$\pi_{0,10}=(10\ 11)(12\ 17)(15\ 16)$,&
$\pi_{0,11}=(10\ 17)(11\ 15)(12\ 16)(13\ 14)$,\\
$\pi_{0,12}=(10\ 11)(12\ 17)(13\ 14)(15\ 16)$,\\
$\pi_{0,13}=(0\ 6\ 4\ 3\ 7\ 1\ 2\ 5)(10\ 12\ 14\ 13)(11\ 15)(16\ 17)$,\\
$\pi_{0,14}=(0\ 5\ 9\ 8\ 6\ 7\ 3\ 1\ 2)(10\ 12\ 11\ 17\ 14)(13\ 16\ 15)$,\\
$\pi_{0,15}=(0\ 6\ 3\ 9)(1\ 4\ 8)(2\ 7\ 5)(10\ 15\ 12)(11\ 16)(13\ 17\ 14)$,\\
$\pi_{0,16}=(0\ 9\ 1\ 2)(3\ 4\ 8\ 5\ 6\ 7)(11\ 17\ 13\ 14\ 15\ 16\ 12)$,\\
$\pi_{0,17}=(0\ 7\ 6\ 9\ 5\ 3\ 8\ 1\ 2)(11\ 17\ 13\ 12\ 16\ 14)$,&
$\pi_{0,27}=\pi_{3,24}=\pi_{6,21}=(1)$,\\
$\pi_{8,0}=(3\ 7\ 6\ 4)(11\ 12\ 15)$,&
$\pi_{9,18}=(1)$,\\
$\pi_{12,5}=(2\ 6\ 5\ 3)$,&
$\pi_{12,15}=\pi_{15,0}=\pi_{15,12}=\pi_{18,9}=(1)$,\\
$\pi_{21,0}=\pi_{24,3}=\pi_{27,0}=(1)$.
\end{tabular}
\end{center}

\noindent We have that for each $i,j,s,t$ in Table V,
 $\pi_{s,t}Y=Y$, $|{\cal B}_i\cap \pi_{s,t}{\cal B}_j|=s$ and
 $|T({\cal B}_i\setminus\pi_{s,t}{\cal B}_j)\cap T(\pi_{s,t}{\cal B}_j\setminus{\cal B}_i)|=t$. \qed

\begin{center}
\small  Table V.  \ Fine triangle intersections for $(K_{18}\setminus K_{10},G)$-designs
  \begin{tabular}{ccc|ccc|ccc}\hline
$i$ &  $j$ &  $(s,t)$ & $i$ &  $j$ &  $(s,t)$& $i$ &  $j$ &  $(s,t)$\\
\hline     1  &   1   &  $(M_{18}\cup N_{18})\setminus $\{(0,27), & 1  & 2 & (0,27)&1  &   3   &   (24,3)   \\
&&(24,3),(3,24),(6,21),& 1  &   5  &   (18,9)& 1  &   6  &   (15,12)  \\
&&(18,9),(9,18),(15,12),& 2 &    3  &   (3,24)& 2  &   4  &   (6,21)    \\
&&(12,15),(15,0),(21,0)\}& 2  &   5  &   (9,18) & 2  &   6  &   (12,15)     \\
&&& 2  &   7  &   (15,0) & 2  &   8  &   (21,0)   \\
\hline
\end{tabular}
\end{center}

\noindent {\bf Lemma $3.11$} {\em
Let $M_{19}=\{(0,0),(0,1),\ldots,(0,18),(0,21),(0,31),(31,0)\}$
and $N_{19}=\{(2,17),(3,\\28),(6,25),(9,0),(9,22),(11,20),(12,4),(12,12),(15,3),(17,14),
(18,0),(18,3),(20,11),(22,9),\\(25,0),(25,6)\}$.
Let $G$ be a kite and $(s,t)\in M_{19}\cup N_{19}$.
Then there is a pair of  maximum $(K_{19}\setminus K_{10},G)$-packings with the same vertex set
and the same subgraph $K_{10}$ removed, which intersect in $s$ blocks and
$t+s$ triangles.}

\proof Take the vertex set $X=\{0,1,\ldots,18\}$. Let

\begin{center}
\small
 \begin{tabular}{llllll}
${\cal A}_1:$ &$[10 ,18, 0 -14]$,&$[18 ,17 ,9 -13]$,&$[17 ,16 ,8- 13]$,&$[16, 7 ,15 -11]$,&$[6 ,14, 15- 0]$,\\
&$[14 ,13, 5- 10]$,&$[1 ,18, 11 -9]$,&$[0 ,17, 11- 4]$,&$[17 ,5, 15 -1]$,&$[14, 12, 2 -18]$,\\\vspace{0.1cm}
&$[13, 6, 11- 7]$,&$[16 ,12 ,0 -13]$,&$[15 ,12, 9 -14]$,&$[14 ,11, 8- 10]$,&$[18, 12, 5- 11]$;\\

${\cal A}_2:$ &$[16 ,7, 15- 1]$,&$[6 ,14, 15- 11]$,&$[14 ,13 ,5 -11]$,&$[10, 9, 16- 5]$,&$[1 ,18, 11- 7]$,\\&$[18, 6 ,16- 1]$,
&$[17 ,5, 15- 0]$,&$[13 ,6 ,11 -4]$,&$[15 ,8 ,18- 7]$,&$[7 ,14 ,17 -3]$,\\\vspace{0.1cm}&$[13 ,7, 10- 6]$,&$[18 ,12, 5 -10]$,
&$[6 ,17, 12 -7]$;\\

${\cal A}_3:$ &$[16 ,5 ,15- 11]$,&$[7 ,14 ,15- 0]$,&$[14 ,13 ,6 -10]$,&$[10, 9, 16- 6]$,&$[1 ,18, 11- 5]$,\\&$[18, 7, 16 -1]$,
&$[17 ,6, 15- 1]$,&$[13, 7, 11- 4]$,&$[15, 8, 18- 5]$,&$[5 ,14, 17- 3]$,\\\vspace{0.1cm}&$[13, 5, 10- 7]$,&$[18 ,12 ,6 -11]$,
&$[7 ,17, 12- 5 ]$;\\

${\cal A}_4:$ &$[10 ,18, 0 -13]$,&$[18 ,17 ,9 -13]$,&$[17 ,16 ,8- 13]$,&$[16, 7 ,15 -11]$,&$[6 ,14, 15- 0]$,\\&$[14 ,13, 5- 10]$,
&$[1 ,18, 11 -9]$,&$[0 ,17, 11- 4]$,&$[17 ,5, 15 -1]$,&$[14, 12, 2 -15]$,\\\vspace{0.1cm}&$[13, 6, 11- 7]$,&$[16 ,12 ,0 -14]$,
&$[15 ,12, 9 -14]$,&$[14 ,11, 8- 10]$,&$[18, 12, 5- 11]$;\\

${\cal A}_5:$ &$[10 ,18, 0 -13]$,&$[18 ,17 ,9 -13]$,&$[17 ,16 ,8- 13]$,&$[16, 7 ,15 -1]$,&$[6 ,14, 15- 11]$,\\&$[14 ,13, 5- 10]$,
&$[1 ,18, 11 -9]$,&$[0 ,17, 11- 4]$,&$[17 ,5, 15 -0]$,&$[14, 12, 2 -15]$,\\\vspace{0.1cm}&$[13, 6, 11- 7]$,&$[16 ,12 ,0 -14]$,
&$[15 ,12, 9 -14]$,&$[14 ,11, 8- 10]$,&$[18, 12, 5- 11]$;\\

${\cal A}_6:$ &$[10 ,18, 0 -13]$,&$[18 ,17 ,9 -13]$,&$[17 ,16 ,8- 13]$,&$[16, 7 ,15 -1]$,&$[6 ,14, 15- 11]$,\\&$[14 ,13, 5- 10]$,
&$[1 ,18, 11 -7]$,&$[0 ,17, 11- 9]$,&$[17 ,5, 15 -0]$,&$[14, 12, 2 -15]$,\\\vspace{0.1cm}&$[13, 6, 11- 4]$,&$[16 ,12 ,0 -14]$,
&$[15 ,12, 9 -14]$,&$[14 ,11, 8- 10]$,&$[18, 12, 5- 11]$;\\

${\cal A}_{7}:$ &$[10 ,18, 0- 13]$,&$[18 ,17 ,9 -14]$,&$[17, 16, 8- 10]$,&$[16, 7 ,15 -11]$,&$[1 ,14, 15 -0]$,\\&$[14, 13, 5- 10]$,
&$[13 ,12, 4 -18]$,&$[12, 11, 3 -16]$,&$[11, 2, 10- 15]$,&$[10 ,9, 16- 5]$,\\&$[6 ,18, 11- 9]$,&$[0 ,17 ,11- 4]$,
&$[18 ,1 ,16 -6]$,&$[17 ,5, 15- 6]$,&$[16 ,14, 4- 15]$,\\\vspace{0.1cm}&$[15, 13, 3 -10]$,&$[14, 12, 2- 18]$,&$[13, 1 ,11- 7]$;\\

${\cal A}_8:$ &$[10 ,18, 0 -13]$,&$[18 ,17 ,9 -14]$,&$[17 ,16 ,8- 10]$,&$[16, 7 ,15 -1]$,&$[6 ,14, 15- 11]$,\\&$[14 ,13, 5- 11]$,
&$[1 ,18, 11 -9]$,&$[0 ,17, 11- 4]$,&$[17 ,5, 15 -0]$,&$[14, 12, 2 -18]$,\\\vspace{0.1cm}&$[13, 6, 11- 7]$,&$[16 ,12 ,0 -14]$,
&$[15 ,12, 9 -13]$,&$[14 ,11, 8- 13]$,&$[18, 12, 5- 10]$;\\

${\cal A}_9:$ &$[10 ,18, 0 -13]$,&$[18 ,17 ,9 -14]$,&$[17 ,16 ,8- 10]$,&$[16, 7 ,15 -1]$,&$[6 ,14, 15- 11]$,\\&$[14 ,13, 5- 11]$,
&$[1 ,18, 11 -7]$,&$[0 ,17, 11- 9]$,&$[17 ,5, 15 -0]$,&$[14, 12, 2 -18]$,\\\vspace{0.1cm}&$[13, 6, 11- 4]$,&$[16 ,12 ,0 -14]$,
&$[15 ,12, 9 -13]$,&$[14 ,11, 8- 13]$,&$[18, 12, 5- 10]$;\\

${\cal A}_{10}:$ &$[18 ,17 ,9- 14]$,&$[16, 5, 15 -11]$,&$[12, 11, 3 -10]$,&$[11 ,2, 10- 15]$,&$[10, 9, 16 -6]$,\\&$[0, 17, 11 -9]$,
&$[18, 7 ,16- 1]$,&$[17, 6 ,15- 1]$,&$[16 ,14, 4- 18]$,&$[15 ,13 ,3- 16]$,\\&$[14 ,12, 2- 15]$,&$[1 ,13 ,17 -2]$,
&$[16, 12, 0- 14]$,&$[5 ,14, 17- 3]$,&$[2 ,13, 16- 11]$,\\\vspace{0.1cm}&$[3, 14 ,18 -13]$,&$[17, 4 ,10- 14]$,&$[7, 17 ,12 -5]$;\\

${\cal A}_{11}:$ &$[18 ,16 ,9- 14]$,&$[17, 5, 15 -11]$,&$[12, 11, 3 -10]$,&$[11 ,2, 10- 15]$,&$[10, 9, 17 -6]$,\\&$[0, 16, 11 -9]$,
&$[18, 7 ,17- 1]$,&$[16, 6 ,15- 1]$,&$[17 ,14, 4- 18]$,&$[15 ,13 ,3- 16]$,\\&$[14 ,12, 2- 15]$,&$[1 ,13 ,16 -2]$,
&$[17, 12, 0- 14]$,&$[5 ,16, 14- 1]$,&$[2 ,13, 17- 11]$,\\\vspace{0.1cm}&$[3, 14 ,18 -13]$,&$[16, 4 ,10- 14]$,&$[7, 16 ,12 -5]$;\\

${\cal A}_{12}:$ &$[18, 17 ,9 -14]$,&$[16, 5, 15- 11]$,&$[12 ,11, 2 -10]$,&$[11 ,3 ,10 -15]$,&$[10 ,9 ,16- 6]$,\\&$[0, 17, 11- 9]$,
&$[18, 7 ,16- 1]$,&$[17, 6, 15- 1]$,&$[16 ,14 ,4 -18]$,&$[15, 13, 2 -16]$,\\&$[14, 12, 3- 15]$,&$[1 ,13, 17 -2]$,
&$[16 ,12, 0 -14]$,&$[5 ,14, 17 -3]$,&$[3, 13 ,16 -11 ]$,\\\vspace{0.1cm}&$[2 ,14 ,18- 13]$,&$[17, 4 ,10 -14]$,&$[7 ,17, 12- 5]$;\\

${\cal A}_{13}:$ &$[13, 12, 4 -18]$,&$[12, 11 ,3 -16]$,&$[11, 2, 10- 6]$,&$[10 ,9, 16- 1]$,&$[18, 6, 16- 11]$,\\&$[16 ,14, 4 -15]$,
&$[ 15 ,13 ,3- 10]$,&$[1 ,10 ,12- 7]$,&$[1 ,13 ,17- 3]$,&$[15, 8 ,18 -13]$,\\&$[7 ,14 ,17 -2]$,&$[2 ,13, 16- 5]$,
&$[13 ,7 ,10- 14]$,&$[3 ,14 ,18- 7  ]$,&$[17 ,4 ,10 -15]$,\\\vspace{0.1cm}&$[6 ,17 ,12- 8]$;\\

${\cal A}_{14}:$ &$[10, 18 ,0 -13]$,&$[18, 17 ,9 -14]$,&$[17, 16, 8 -10]$,&$[13 ,12 ,4 -15]$,&$[12 ,11, 3- 10]$,\\&$[11 ,2 ,10 -15]$,
&$[0 ,17 ,11- 9]$,&$[16 ,14 ,4 -18]$,&$[15 ,13 ,3 -16]$,&$[14 ,12 ,2- 15]$,\\&$[1 ,10 ,12- 8]$,&$[1 ,13 ,17- 2]$,
&$[16, 12, 0 -14]$,&$[ 2 ,13 ,16- 11]$,&$[15 ,12 ,9- 13]$,\\\vspace{0.1cm}&$[14 ,11, 8 -13]$,&$[3 ,14 ,18- 13]$,&$[17 ,4, 10- 14]$;\\
\end{tabular}
\end{center}

\begin{center}
\small
 \begin{tabular}{llllll}
${\cal A}_{15}:$ &$[10, 18, 0- 13]$,&$[17 ,16, 8- 10]$,&$[7 ,14, 15- 0]$,&$[14 ,13, 6- 10]$,&$[13 ,12, 4- 15]$,\\&$[1 ,18, 11- 5]$,
&$[13, 7 ,11- 4]$,&$[1 ,10 ,12 -8]$,&$[15 ,8 ,18- 5]$,&$[15, 12 ,9 -13]$,\\\vspace{0.1cm}&$[14 ,11 ,8 -13]$,&$[13 ,5, 10 -7]$,
&$[18 ,12, 6 -11]$;\\

${\cal A}_{16}:$ &$[6 ,10 ,12- 7]$,&$[6 ,13, 17- 3]$,&$[16 ,12, 0 -14]$,&$[15, 8, 18 -13]$,&$[ 7 ,14 ,17- 2]$,\\&$[2 ,13 ,16- 11]$,
&$[15 ,12, 9 -13]$,&$[14 ,11, 8- 13]$,&$[13 ,7, 10- 1]$,&$[3 ,14, 18 -7]$,\\&$[18, 12, 5- 11]$,&$[17 ,4, 10 -14]$,
&$[1 ,17 ,12- 8]$.
\end{tabular}
\end{center}
\noindent Let ${\cal B}_{i}={\cal A}_{i}\cup {\cal A}_{13}$ when $i=1,4,5,6,8,9$. Let ${\cal B}_{i}={\cal A}_{i}\cup {\cal A}_{14}$ when $i=2,3$. Let ${\cal B}_{7}={\cal A}_{7}\cup {\cal A}_{16}$. Let ${\cal B}_{i}={\cal A}_{i}\cup {\cal A}_{15}$ when $i=10,11,12$.
Then $(X,{\cal B}_{i})$ is a  maximum $(K_{19}\setminus K_{10},G)$-packing for each $1\leq i\leq 12$,
where the removed subgraph $K_{10}$ is constructed on $Y=\{0,1,\ldots,9\}$. Consider the following permutations
on $X$.

\begin{center}
\small
\begin{tabular}{ll}
$\pi_{0,0}=(11\ 12\ 13)(15\ 18\ 16\ 17)$,&
$\pi_{0,1}=(11\ 12)(13\ 14\ 16)(15\ 17\ 18)$,\\
$\pi_{0,2}=(11\ 12)(14\ 16\ 18\ 15\ 17)$,&
$\pi_{0,3}=(11\ 12)(14\ 15\ 16\ 17\ 18)$,\\
$\pi_{0,4}=(11\ 12)(14\ 15\ 16\ 17)$,&
$\pi_{0,5}=(11\ 12)(14\ 15)(16\ 17\ 18)$,\\
$\pi_{0,6}=(11\ 12)(14\ 15\ 17\ 16)$,&
$\pi_{0,7}=(11\ 12)(14\ 16)(15\ 17)$,\\
$\pi_{0,8}=(11\ 12)(13\ 17)(14\ 15)(16\ 18)$,&
$\pi_{0,9}=(11\ 12)(13\ 18)(14\ 15)(16\ 17)$,\\
$\pi_{0,10}=(11\ 13)(12\ 15)(14\ 16)(17\ 18)$,&
$\pi_{0,11}=(10\ 11\ 13)(12\ 17)(14\ 16)(15\ 18)$,\\
$\pi_{0,12}=(11\ 17)(12\ 13)(14\ 16)(15\ 18)$,&
$\pi_{0,13}=(10\ 13)(12\ 17)(15\ 18)$,\\
$\pi_{0,14}=(10\ 13)(12\ 17)(14\ 16)(15\ 18)$,&
$\pi_{0,15}=(10\ 13)(11\ 16)(12\ 17)(15\ 18)$,\\
$\pi_{0,16}=(0\ 2\ 3\ 1\ 7\ 6)(4\ 5\ 9)$,&
$\pi_{0,17}=(0\ 1)(2\ 8\ 5)(3\ 4\ 9\ 6\ 7)$,\\
\hspace{1.2cm}$(10\ 13\ 11\ 15\ 17\ 14\ 18\ 16\ 12)$&\hspace{1.2cm}$(10\ 18\ 11\ 17\ 13\ 14\ 15\ 16\ 12)$\\
$\pi_{0,18}=(5\ 6\ 9\ 8)$,&
$\pi_{0,21}=(1)$,\\
$\pi_{0,31}=(1)$,&
$\pi_{2,17}=(0\ 9)(2\ 3)(5\ 6)(11\ 14)(15\ 16)(17\ 18)$,\\
$\pi_{3,28}=\pi_{6,25}=(1)$,&
$\pi_{9,0}=(1\ 7\ 6\ 4\ 3)(2\ 8)$,\\
$\pi_{9,22}=\pi_{11,20}=(1)$,&
$\pi_{12,4}=(0\ 5\ 8\ 1\ 4)$,\\
$\pi_{12,12}=(1)$,&
$\pi_{15,3}=(5\ 6\ 9\ 8)$,\\
$\pi_{17,14}=\pi_{18,0}=\pi_{18,3}=\pi_{20,11}=(1)$,&
$\pi_{22,9}=\pi_{25,0}=\pi_{25,6}=\pi_{31,0}=(1)$..
\end{tabular}
\end{center}

\noindent We have that for each $i,j,s,t$ in Table VI,
 $\pi_{s,t}Y=Y$, $|{\cal B}_i\cap \pi_{s,t}{\cal B}_j|=s$ and
 $|T({\cal B}_i\setminus\pi_{s,t}{\cal B}_j)\cap T(\pi_{s,t}{\cal B}_j\setminus{\cal B}_i)|=t$.\qed

\begin{center}
\small  Table VI.  \ Fine triangle intersections for maximum $(K_{19}\setminus K_{10},G)$-packings
  \begin{tabular}{ccc|ccc|ccc}\hline
$i$ &  $j$ &  $(s,t)$&$i$ &  $j$ &  $(s,t)$&$i$ &  $j$ &  $(s,t)$ \\
\hline
1  &   1   &  $(M_{19}\cup N_{19})\setminus \{(0,18),$&1  &   2   &   (0,18),(0,31)&1  &   3  &   (0,21)   \\
&&(0,21),(0,31),(3,28),&1  &   5  &   (25,6) &1  &   6  &   (22,9)    \\
&&(6,25),(9,0),(9,22),&1  &   7  &   (12,12)&1  &   8  &   (20,11)    \\
&&(11,20),(12,12),(17,14),&1  &   9  &   (17,14) &2 &    2   &   (9,0)    \\
&&(18,0),(18,3),(20,11),& 2  &   3  &   (18,3) &2  &   4  &   (3,28)    \\
&&(22,9),(25,0),(25,6)$\}$& 2  &   5  &   (6,25)& 2  &   6  &   (9,22)   \\
&&& 2  &   8  &   (11,20) &  10  &   11  &   (18,0) \\
&&&  10  &   12  &   (25,0)    \\
\hline
\end{tabular}
\end{center}

\noindent {\bf Lemma $3.12$} {\em
Let $M_{20}=\{(0,0),(0,1),\ldots,(0,17),(0,22),(0,24),(0,28),(0,36),(36,0)\}$
and $N_{20}=\{(2,16),(5,31),(7,21),(10,18),(10,26),(12,0),(12,3),(12,6),(13,15),(13,23),
(15,13),(18,10),\\(18,18),(21,1),(21,7),(23,13),(24,0),(26,10),(28,0),(31,5)\}$.
Let $G$ be a kite and $(s,t)\in M_{20}\cup N_{20}$.
Then there is a pair of  maximum $(K_{20}\setminus K_{10},G)$-packings with the same vertex set
and the same subgraph $K_{10}$ removed, which intersect in $s$ blocks and
$t+s$ triangles.}

\proof Take the vertex set $X=\{0,1,\ldots,19\}$. Let
\begin{center}
\small
 \begin{tabular}{llllll}

${\cal A}_1:$ &$[0 ,19, 10- 15]$,&$[18 ,19, 9- 12]$,&$[18 ,17 ,8 -13]$,&$[17 ,16, 7- 11]$,&$[16 ,15, 6 -12]$,\\&$[5 ,14, 15 -11]$,
&$[14 ,13, 4- 19]$,&$[11, 1, 10- 8]$,&$[10 ,17, 9 -13]$,&$[17, 19, 6- 13]$,\\&$[17 ,4 ,15- 1]$,&$[2 ,13, 15- 9]$,
&$[12 ,10, 4- 18]$,&$[11, 16 ,9 -14]$,&$[19 ,16 ,8 -14]$,\\\vspace{0.1cm}&$[18, 15 ,7- 12]$,&$[15, 12, 8- 11]$,&$[14 ,10, 7- 13]$;\\

${\cal A}_2:$ &$[ 13, 12, 3 -10]$,&$[12 ,2 ,11- 4]$,&$[1 ,12, 19- 14]$,&$[0 ,18, 11- 17]$,&$[16, 5 ,18- 13]$,\\&$[16, 14, 3 -18]$,
&$[14, 12, 0 -15]$,&$[13 ,11 ,19- 7]$,&$[13 ,17 ,0- 16]$,&$[14,2,17-1]$,\\&$[1, 13, 16 -12]$,&$[14 ,6, 11- 3]$,
&$[13 ,10, 5- 19]$,&$[14 ,1 ,18- 12]$,&$[19 ,15, 3 -17]$,\\\vspace{0.1cm}&$[2 ,10 ,16- 4]$,&$[6 ,10, 18- 2]$,&$[17 ,12, 5- 11]$;\\

${\cal A}_3:$ &$[18 ,19, 9 -12]$,&$[18 ,17, 8- 13]$,&$[17 ,16,7 -11]$,&$[16, 15, 6- 12]$,&$[14 ,13, 4- 19]$,\\&$[10 ,17, 9- 13]$,
&$[17 ,19, 6 -13]$,&$[12, 10, 4- 18]$,&$[11 ,16 ,9- 14]$,&$[19 ,16, 8- 14]$,\\\vspace{0.1cm}&$[18 ,15, 7 -12]$,&$[15 ,12 ,8 -11]$,
&$[14, 10, 7 -13]$;\\

${\cal A}_4:$ &$[0 ,19, 10- 8]$,&$[18 ,19, 9- 14]$,&$[18 ,17 ,8 -13]$,&$[17 ,16, 7- 11]$,&$[16 ,15, 6 -12]$,\\&$[5 ,14, 15 -11]$,
&$[14 ,13, 4- 19]$,&$[11, 1, 10- 15]$,&$[10 ,17, 9 -12]$,\\&$[17, 19, 6- 13]$,&$[17 ,4 ,15- 1]$,&$[2 ,13, 15- 9]$,
&$[12 ,10, 4- 18]$,&$[11, 16 ,9 -13]$,\\\vspace{0.1cm}&$[19 ,16 ,8 -14]$,&$[18, 15 ,7- 12]$,&$[15, 12, 8- 11]$,&$[14 ,10, 7- 13]$;\\

${\cal A}_5:$ &$[18 ,19, 9 -14]$,&$[18 ,17, 8- 13]$,&$[17 ,16,7 -11]$,&$[16, 15, 6- 13]$,&$[14 ,13, 4- 18]$,\\&$[10 ,17, 9- 12]$,
&$[17 ,19, 6 -12]$,&$[12, 10, 4- 19]$,&$[11 ,16 ,9- 13]$,&$[19 ,16, 8- 14]$,\\\vspace{0.1cm}&$[18 ,15, 7 -12]$,&$[15 ,12 ,8 -11]$,
&$[14, 10, 7 -13]$;\\

${\cal A}_6:$ &$[18 ,19, 9 -14]$,&$[18 ,17, 8- 11]$,&$[17 ,16,7 -11]$,&$[16, 15, 6- 13]$,&$[14 ,13, 4- 18]$,\\&$[10 ,17, 9- 12]$,
&$[17 ,19, 6 -12]$,&$[12, 10, 4- 19]$,&$[11 ,16 ,9- 13]$,&$[19 ,16, 8- 13]$,\\\vspace{0.1cm}&$[18 ,15, 7 -12]$,&$[15 ,12 ,8 -14]$,
&$[14, 10, 7 -13]$;\\

${\cal A}_7:$ &$[0 ,19, 10- 8]$,&$[18 ,19, 9- 14]$,&$[18 ,17 ,8 -11]$,&$[17 ,16, 7- 11]$,&$[16 ,15, 6 -12]$,\\&$[5 ,14, 15 -11]$,
&$[14 ,13, 4- 18]$,&$[11, 1, 10- 15]$,&$[10 ,17, 9 -12]$,&$[17, 19, 6- 13]$,\\&$[17 ,4 ,15- 1]$,&$[2 ,13, 15- 9]$,
&$[12 ,10, 4- 19]$,&$[11, 16 ,9 -13]$,&$[19 ,16 ,8 -13]$,\\\vspace{0.1cm}&$[18, 15 ,7- 12]$,&$[15, 12, 8- 14]$,&$[14 ,10, 7- 13]$;\\

${\cal A}_8:$ &$[18 ,19, 9 -14]$,&$[18 ,17, 8- 11]$,&$[17 ,16,7 -13]$,&$[16, 15, 6- 13]$,&$[14 ,13, 4- 18]$,\\&$[10 ,17, 9- 12]$,
&$[17 ,19, 6 -12]$,&$[12, 10, 4- 19]$,&$[11 ,16 ,9- 13]$,&$[19 ,16, 8- 13]$,\\\vspace{0.1cm}&$[18 ,15, 7 -11]$,&$[15 ,12 ,8 -14]$,
&$[14, 10, 7 -12]$;\\

${\cal A}_9:$ &$[0 ,19, 10- 8]$,&$[18 ,19, 9- 14]$,&$[18 ,17 ,8 -11]$,&$[17 ,16, 7- 13]$,&$[16 ,15, 6 -12]$,\\&$[5 ,14, 15 -11]$,
&$[14 ,13, 4- 18]$,&$[11, 1, 10- 15]$,&$[10 ,17, 9 -12]$,&$[17, 19, 6- 13]$,\\&$[17 ,4 ,15- 1]$,&$[2 ,13, 15- 9]$,
&$[12 ,10, 4- 19]$,&$[11, 16 ,9 -13]$,&$[19 ,16 ,8 -13]$,\\\vspace{0.1cm}&$[18, 15 ,7- 11]$,&$[15, 12, 8- 14]$,&$[14 ,10, 7- 12]$;\\

${\cal A}_{10}:$ &$[0 ,19, 10- 8]$,&$[18 ,19, 9- 14]$,&$[18 ,17 ,8 -11]$,&$[17 ,16, 7- 13]$,&$[16 ,15, 6 -13]$,\\&$[5 ,14, 15 -9]$,
&$[14 ,13, 4- 18]$,&$[11, 1, 10- 15]$,&$[10 ,17, 9 -12]$,&$[17, 19, 6- 12]$,\\&$[17 ,4 ,15- 11]$,&$[2 ,13, 15- 1]$,
&$[12 ,10, 4- 19]$,&$[11, 16 ,9 -13]$,&$[19 ,16 ,8 -13]$,\\\vspace{0.1cm}&$[18, 15 ,7- 11]$,&$[15, 12, 8- 14]$,&$[14 ,10, 7- 12]$;\\

${\cal A}_{11}:$ &$[5 ,19, 10 -15]$,&$[18 ,19, 9- 12]$,&$[18 ,17 ,8- 13]$,&$[17 ,16, 7- 11]$,&$[16 ,15, 0 -12]$,\\&$[6 ,14, 15 -11]$,
&$[14, 13, 4 -19]$,&$[13 ,12, 3 -18]$,&$[12, 2, 11 -17]$,&$[11 ,1, 10- 8]$,\\&$[10, 17, 9- 13]$,&$[1 ,12, 19- 7]$,
&$[5 ,18, 11- 3]$,&$[17, 19 ,0 -13]$,&$[16, 6 ,18- 12]$,\\\vspace{0.1cm}&$[17, 4, 15- 1]$,&$[16, 14, 3 -17]$,&$[2 ,13, 15- 9]$;\\

${\cal A}_{12}:$ &$[ 13, 12, 3 -18]$,&$[12 ,2 ,11- 17]$,&$[1 ,12, 19- 7]$,&$[0 ,18, 11- 3]$,&$[16, 5 ,18- 12]$,\\&$[16, 14, 3 -17]$,
&$[14, 12, 0 -16]$,&$[13 ,11 ,19- 14]$,&$[13 ,17 ,0- 15]$,&$[17, 14, 2 -19]$,\\&$[1, 13, 16 -4]$,&$[14 ,6, 11- 4]$,
&$[13 ,10, 5- 11]$,&$[14 ,1 ,18- 2]$,&$[19 ,15, 3 -10]$,\\\vspace{0.1cm}&$[2 ,10 ,16- 12]$,&$[6 ,10, 18- 13]$,&$[17 ,12, 5- 19]$;\\

${\cal A}_{13}:$ &$[5 ,19, 10- 15]$,&$[0 ,14 ,15 -11]$,&$[13 ,12, 3- 18]$,&$[12, 2, 11 -17]$,&$[11, 1 ,10- 8]$,\\&$[1 ,12, 19 -7]$,
&$[5 ,18, 11- 3]$,&$[16 ,0, 18- 12]$,&$[ 17, 4 ,15- 1]$,&$[16, 14, 3- 17]$,\\&$[2 ,13, 15- 9]$,&$[14, 12 ,5 -16]$,
&$[13, 11, 19- 14]$,&$[ 13 ,17 ,5- 15]$,&$[17 ,14, 2- 19]$,\\&$[1 ,13, 16 -4]$,&$[14 ,6 ,11 -4]$,&$[13 ,10, 0- 11]$,
&$[14 ,1, 18- 2]$,&$[19 ,15 ,3- 10]$,\\\vspace{0.1cm}&$[2 ,10 ,16- 12]$,&$[6 ,10 ,18- 13]$,&$[17 ,12 ,0 -19]$;\\

${\cal A}_{14}:$ &$[14, 12, 5- 16]$,&$[13, 11, 19 -14]$,&$[12 ,10 ,4 -18]$,&$[11 ,16 ,9 -14]$,&$[13, 17, 5 -15]$,\\&$[19 ,16, 8 -14]$,
&$[18 ,15 ,7- 12]$,&$[17 ,14, 2- 19]$,&$[1 ,13, 16- 4]$,&$[15 ,12 ,8- 11]$,\\&$[14 ,0, 11- 4]$,&$[13 ,10, 6- 11]$,
&$[14 ,1 ,18- 2]$,&$[19 ,15, 3- 10]$,&$[2 ,10, 16 -12]$,\\&$[14, 10, 7 -13]$,&$[0 ,10 ,18- 13]$,&$[17, 12, 6- 19 ]$.\\
\end{tabular}
\end{center}
\noindent Let ${\cal B}_{i}={\cal A}_{i}\cup {\cal A}_{12}$ when $i=1,4,7,9,10$. Let
${\cal B}_{i}={\cal A}_{i}\cup {\cal A}_{13}$ when $i=3,5,6,8$. Let ${\cal B}_{2}={\cal A}_{2}\cup {\cal A}_{10}$ and ${\cal B}_{11}={\cal A}_{11}\cup {\cal A}_{14}$.
Then $(X,{\cal B}_{i})$ is a  maximum $(K_{20}\setminus K_{10},G)$-packing for each $1\leq i\leq 11$,
where the removed subgraph $K_{10}$ is constructed on $Y=\{0,1,\ldots,9\}$. Consider the following permutations
on $X$.

\begin{center}
\small
\begin{tabular}{ll}
$\pi_{0,0}=(11\ 12\ 13)(14\ 15\ 19\ 18\ 17\ 16)$,&
$\pi_{0,1}=(11\ 12)(13\ 14\ 15)(16\ 17\ 18\ 19)$,\\
$\pi_{0,2}=(11\ 12)(14\ 15\ 16)(17\ 18\ 19)$,&
$\pi_{0,3}=(11\ 12)(14\ 15\ 16\ 17\ 18\ 19)$,\\
$\pi_{0,4}=(11\ 12)(14\ 15)(16\ 17\ 18)$,&
$\pi_{0,5}=(11\ 12)(14\ 15)(16\ 17\ 19\ 18)$,\\
$\pi_{0,6}=(11\ 12)(14\ 15\ 17\ 18)(16\ 19)$,&
$\pi_{0,7}=(11\ 12)(14\ 15)(16\ 17)(18\ 19)$,\\
$\pi_{0,8}=(11\ 12)(13\ 14)(17\ 19)$,&
$\pi_{0,9}=(11\ 12)(14\ 18)(15\ 17)(16\ 19)$,\\
$\pi_{0,10}=(11\ 12\ 13\ 14\ 15)$,&
$\pi_{0,11}=(11\ 12)(13\ 14\ 15)$,\\
$\pi_{0,12}=(11\ 13)(12\ 14\ 15)$,&
$\pi_{0,13}=(11\ 12)(13\ 17)(14\ 15)$,\\
$\pi_{0,14}=(10\ 12)(11\ 19)(13\ 14)(15\ 17)(16\ 18)$,&
$\pi_{0,15}=(10\ 12)(11\ 16)(13\ 17)(14\ 15)(18\ 19)$,\\
$\pi_{0,16}=(10\ 17)(11\ 15)(12\ 13)(14\ 16)(18\ 19)$,&
$\pi_{0,17}=(10\ 14)(11\ 18)(12\ 19)(13\ 15)(16\ 17)$,\\
$\pi_{0,22}=(0\ 3\ 4\ 6)$,&
$\pi_{0,24}=\pi_{0,28}=\pi_{0,36}=(1)$,\\
$\pi_{2,16}=(0\ 2\ 3\ 5\ 1)(4\ 7\ 6)(10\ 11\ 18\ 19\ 12\ 13)(15\ 17\ 16)$,&
$\pi_{5,31}=\pi_{7,21}=\pi_{10,18}=\pi_{10,26}=(1)$,\\
$\pi_{12,0}=(0\ 3\ 4)(1\ 6)(7\ 8)$,&
$\pi_{12,3}=(0\ 7\ 2\ 6)(3\ 4)$,\\
$\pi_{12,6}=(5\ 9)(6\ 8)(16\ 17)$,&
$\pi_{13,15}=\pi_{13,23}=\pi_{15,13}=\pi_{18,10}=(1)$,\\
$\pi_{18,18}=(1)$,&
$\pi_{21,1}=(0\ 3\ 4\ 6)$,\\
$\pi_{21,7}=\pi_{23,13}=\pi_{24,0}=\pi_{26,10}=(1)$,&
$\pi_{28,0}=\pi_{31,5}=\pi_{36,0}=(1)$.
\end{tabular}
\end{center}
\noindent We have that for each $i,j,s,t$ in Table VII,
 $\pi_{s,t}Y=Y$, $|{\cal B}_i\cap \pi_{s,t}{\cal B}_j|=s$ and
 $|T({\cal B}_i\setminus\pi_{s,t}{\cal B}_j)\cap T(\pi_{s,t}{\cal B}_j\setminus{\cal B}_i)|=t$. \qed
\begin{center}
\small  Table VII.  \ Fine triangle intersections for maximum $(K_{20}\setminus K_{10},G)$-packings
  \begin{tabular}{ccc|ccc|ccc}\hline
$i$ &  $j$ &  $(s,t)$&$i$ &  $j$ &  $(s,t)$&$i$ &  $j$ &  $(s,t)$ \\
\hline     1  &   1   &  $(M_{20}\cup N_{20})\setminus \{(0,22)$,&1  &   2   &  (0,22),(0,36)&1  &   3  &  (28,0)    \\
&&(0,24),(0,28),(0,36),& 1  &   4  &  (31,5)& 1  &   5  &   (21,7)   \\
&&(5,31),(7,21),(10,18),&1  &   6  &   (18,10)&1  &   7  &   (26,10)    \\
&&(10,26),(13,15),(13,23),&1  &   8  &   (15,13)&1  &   9  &   (23,13)    \\
&&(15,13),(18,10),(18,18),&1  &   10  &   (18,18) &1  &   11  &   (24,0)    \\
&&(21,7),(23,13),(24,0),&2  &   3  &  (0,28)   &2  &   4  &  (5,31)   \\
&&(26,10),(28,0),(31,5)$\}$&2  &   5  &   (7,21) &2  &   6  &   (10,18)  \\
&&&2  &   7  &   (10,26)&2  &   8  &   (13,15)    \\
&&&2  &   9  &   (13,23) &2  &   11  &  (0,24)     \\
\hline
\end{tabular}
\end{center}

\noindent {\bf Lemma $3.13$} {\em
Let $M_{21}=\{(0,0),(0,1),\ldots,(0,19),(0,24),(0,36),(36,0)\}$ and $N_{21}=\{(6,30),(7,\\0),\ldots,(7,7),(12,5),(12,24),(13,11),(18,18),(24,0),(24,12),(30,6)\}$.
Let $G$ be a kite and $(s,t)\in M_{21}\cup N_{21}$.
Then there is a pair of  $(K_{21}\setminus K_{12},G)$-designs with the same vertex set
and the same subgraph $K_{12}$ removed, which intersect in $s$ blocks and
$t+s$ triangles.}

\proof Take the vertex set $X=\{0,1,\ldots,20\}$. Let
\begin{center}
\small
 \begin{tabular}{llllll}
${\cal A}_1:$ &$[0 ,12, 20- 8]$,&$[20 ,9, 19- 3]$,&$[8 ,18, 19- 2]$,&$[7 ,18, 17- 0]$,&$[17, 6 ,16- 11]$,\\
&$[16, 5 ,15- 0]$,&$[4 ,14 ,15- 7]$,&$[6, 18, 20- 1]$,&$[5 ,17 ,19 -7]$,&$[18 ,4 ,16- 0]$,\\&$[3 ,15, 17 -1]$,&$[17 ,10, 20- 4]$,
&$[19, 10 ,16 -1]$,&$[15 ,2, 20- 5]$,&$[14, 1 ,19- 4]$,\\\vspace{0.1cm}&$[4 ,12, 17- 2]$,&$[9 ,13 ,17 -11]$,&$[3 ,20, 16- 9]$;\\

${\cal A}_2:$ &$[0 ,12, 20- 1]$,&$[20 ,9, 19- 2]$,&$[8 ,18, 19- 3]$,&$[7 ,18, 17- 1]$,&$[17, 6 ,16- 0]$,\\
&$[16, 5 ,15- 7]$,&$[4 ,14 ,15- 0]$,&$[6, 18, 20- 8]$,&$[5 ,17 ,19 -4]$,&$[18 ,4 ,16- 11]$,\\&$[3 ,15, 17 -0]$,&$[17 ,10, 20- 5]$,
&$[19, 10 ,16 -9]$,&$[15 ,2, 20- 4]$,&$[14, 1 ,19- 7]$,\\\vspace{0.1cm}&$[4 ,12, 17- 11]$,&$[9 ,13 ,17 -2]$,&$[3 ,20, 16- 1]$;\\
${\cal A}_3:$ &$[0 ,12, 20- 8]$,&$[8 ,18 ,19- 2]$,&$[16 ,5 ,15 -0]$,&$[4 ,14, 15 -7]$,&$[13 ,2, 12- 5]$,\\
\end{tabular}
\end{center}

\begin{center}
\small
 \begin{tabular}{llllll}

&$[19, 11, 12- 7]$,&$[0 ,19, 13 -6]$,&$[6 ,18 ,20- 1]$,&$[5 ,17, 19- 7]$,&$[1, 13 ,15- 11]$,\\\vspace{0.1cm}&$[20, 11, 13 -5]$,&$[14, 1, 19- 4]$,
&$[19, 6, 15- 8]$;\\

${\cal A}_4:$ &$[0 ,12, 20- 1]$,&$[20 ,9, 19- 2]$,&$[8 ,18, 19- 3]$,&$[7 ,18, 17- 0]$,&$[17, 6 ,16- 11]$,\\
&$[16, 5 ,15- 0]$,&$[4 ,14 ,15- 7]$,&$[6, 18, 20- 8]$,&$[5 ,17 ,19 -7]$,&$[18 ,4 ,16- 0]$,\\&$[3 ,15, 17 -1]$,&$[17 ,10, 20- 5]$,
&$[19, 10 ,16 -1]$,&$[15 ,2, 20- 4]$,&$[14, 1 ,19- 4]$,\\\vspace{0.1cm}&$[4 ,12, 17- 2]$,&$[9 ,13 ,17 -11]$,&$[3 ,20, 16- 9]$;\\

${\cal A}_5:$ &$[0 ,12, 20- 1]$,&$[20 ,9, 19- 2]$,&$[8 ,18, 19- 3]$,&$[7 ,18, 17- 1]$,&$[17, 6 ,16- 11]$,\\
&$[16, 5 ,15- 0]$,&$[4 ,14 ,15- 7]$,&$[6, 18, 20- 8]$,&$[5 ,17 ,19 -4]$,&$[18 ,4 ,16- 0]$,\\&$[3 ,15, 17 -0]$,&$[17 ,10, 20- 5]$,
&$[19, 10 ,16 -1]$,&$[15 ,2, 20- 4]$,&$[14, 1 ,19- 7]$,\\\vspace{0.1cm}&$[4 ,12, 17- 11]$,&$[9 ,13 ,17 -2]$,&$[3 ,20, 16- 9]$;\\

${\cal A}_6:$ &$[0 ,12, 20- 1]$,&$[8 ,18 ,19- 4]$,&$[16 ,5 ,15 -7]$,&$[4 ,14, 15 -0]$,&$[13 ,2, 12- 7]$,\\
&$[19, 11, 12- 5]$,&$[0 ,19, 13 -5]$,&$[6 ,18 ,20- 8]$,&$[5 ,17, 19- 2]$,&$[1, 13 ,15- 8]$,\\\vspace{0.1cm}&$[20, 11, 13 -6]$,&$[14, 1, 19- 7]$,
&$[19, 6, 15- 11]$;\\

${\cal A}_7:$ &$[0 ,12, 20- 1]$,&$[20 ,9, 19- 2]$,&$[8 ,18, 19- 3]$,&$[7 ,18, 17- 1]$,&$[17, 6 ,16- 0]$,\\
&$[16, 5 ,15- 7]$,&$[4 ,14 ,15- 0]$,&$[6, 18, 20- 8]$,&$[5 ,17 ,19 -4]$,&$[18 ,4 ,16- 11]$,\\&$[3 ,15, 17 -0]$,&$[17 ,10, 20- 5]$,
&$[19, 10 ,16 -9]$,&$[15 ,2, 20- 4]$,&$[14, 1 ,19- 7]$,\\\vspace{0.1cm}&$[4 ,12, 17- 11]$,&$[9 ,13 ,17 -2]$,&$[3 ,20, 16- 1]$;\\

${\cal A}_8:$ &$ [14 ,3, 13 -4]$,&$[ 13 ,2 ,12- 5]$,&$[19 ,11, 12 -7]$,&$[0, 19 ,13- 6]$,&$[ 2, 16, 14- 6]$,\\
&$[1 ,13 ,15 -11]$,&$[9 ,12, 14- 10]$,&$[ 20, 11, 13- 5]$,&$[9 ,15 ,18- 1]$,&$[17 ,8 ,14 -0]$,\\&$[16 ,7, 13- 8]$,&$[10 ,15, 12- 6]$,
&$[14, 11 ,18- 2]$,&$[13, 10, 18 -5]$,&$[20 ,7 ,14- 5]$,\\\vspace{0.1cm}&$[3 ,12 ,18 -0]$,&$[16, 8, 12 -1]$,&$[19, 6, 15- 8]$;\\

${\cal A}_9:$ &$[ 20, 10, 19- 9]$,&$[7 ,18, 17- 0]$,&$[17, 6, 16 -11]$,&$[14 ,9, 13- 4]$,&$[ 18, 4, 16 -0]$,\\&$[9 ,15, 17- 1]$,
&$[2 ,16, 14- 6]$,&$[ 10, 12, 14 -3]$,&$[17, 3, 20- 4]$,&$[19 ,3, 16 -1]$,\\&$[10 ,15, 18- 1]$,&$[17, 8 ,14 -0]$,
&$[16 ,7 ,13- 8]$,&$[3 ,15 ,12- 6]$,&$[14 ,11 ,18 -2]$,\\&$[13 ,3 ,18- 5]$,&$[15 ,2, 20- 5]$,&$[20, 7 ,14 -5]$,
&$[9 ,12 ,18 -0]$,&$[16, 8 ,12- 1]$,\\\vspace{0.1cm}&$[4 ,12, 17- 2]$,&$[10, 13, 17- 11]$,&$[9 ,20, 16- 10]$;\\

${\cal A}_{10}:$ &$ [14 ,3, 13 -6]$,&$[ 13 ,2 ,12- 7]$,&$[19 ,11, 12 -5]$,&$[0, 19 ,13- 4]$,&$[ 2, 16, 14- 10]$,\\&$[1 ,13 ,15 -8]$,
&$[9 ,12, 14- 6]$,&$[ 20, 11, 13- 8]$,&$[9 ,15 ,18- 2]$,&$[17 ,8 ,14 -5]$,\\&$[16 ,7, 13- 5]$,&$[10 ,15, 12- 1]$,
&$[14, 11 ,18- 1]$,&$[13, 10, 18 -0]$,&$[20 ,7 ,14- 0]$,\\&$[3 ,12 ,18 -5]$,&$[16, 8, 12 -6]$,&$[19, 6, 15- 11]$.
\end{tabular}
\end{center}
\noindent Let ${\cal B}_{i}={\cal A}_{i}\cup {\cal A}_{8}$ when $i=1,4,5,7$. Let ${\cal B}_{i}={\cal A}_{i}\cup {\cal A}_{9}$ when $i=3,6$. Let ${\cal B}_{2}={\cal A}_{2}\cup {\cal A}_{10}$. Then $(X,{\cal B}_{i})$ is a  $(K_{21}\setminus K_{12},G)$-design for each $1\leq i\leq 7$, where the removed subgraph $K_{12}$ is constructed on $Y=\{0,1,\ldots,11\}$. Consider the following permutations on $X$.
\begin{center}\small\tabcolsep 0.04in
\begin{tabular}{ll}
$\pi_{0,0}=(0\ 11\ 4\ 9\ 6\ 5\ 2)(1\ 3\ 10\ 8\ 7)$,&
$\pi_{0,1}=(13\ 14)(15\ 16\ 18\ 19\ 20\ 17)$,\\
$\pi_{0,2}=(0\ 10\ 5\ 11\ 6\ 7\ 3\ 2\ 9\ 8\ 4)$,&
$\pi_{0,3}=(0\ 6\ 8\ 4\ 2\ 3\ 9\ 10)(1\ 7\ 5)$,\\
$\pi_{0,4}=(0\ 6\ 8\ 2\ 1)(3\ 4\ 9)(5\ 11\ 7)$,&
$\pi_{0,5}=(1\ 3\ 7\ 5)(2\ 10\ 4)(6\ 8\ 9)$,\\
$\pi_{0,6}=(0\ 3\ 8\ 2\ 10\ 9\ 7)(1\ 6\ 5)$,&
$\pi_{0,7}=(0\ 10\ 5\ 2\ 6\ 8\ 1\ 7)(4\ 9)$,\\
$\pi_{0,8}=(0\ 11\ 1\ 6\ 5)(2\ 4\ 9\ 8\ 7)$,&
$\pi_{0,9}=(1\ 5\ 7\ 9\ 2)(3\ 6\ 4\ 8)$,\\
$\pi_{0,10}=(0\ 6\ 9\ 5\ 4\ 1\ 7\ 2\ 10)$,&
$\pi_{0,11}=(0\ 2\ 9\ 4\ 8)(1\ 11\ 6\ 5)$,\\
$\pi_{0,12}=(0\ 7\ 2\ 1\ 10\ 6\ 11\ 3)$,&
$\pi_{0,13}=(12\ 13)(14\ 18)(15\ 16)(17\ 19\ 20)$,\\
$\pi_{0,14}=(0\ 6\ 9\ 2\ 5)(1\ 4\ 7)$,&
$\pi_{0,15}=(0\ 10\ 1)(12\ 17)(13\ 20)(14\ 16)(15\ 18)$,\\
$\pi_{0,16}=(12\ 16)(13\ 14)(15\ 19)(18\ 20)$,&
$\pi_{7,0}=(0\ 3\ 6)(1\ 4\ 8\ 9)(15\ 16\ 17)$,\\
$\pi_{0,17}=(0\ 6\ 3\ 7\ 1\ 5\ 4\ 2\ 9\ 8)(12\ 20\ 13\ 19\ 15\ 17\ 18\ 14\ 16)$,&
$\pi_{7,1}=(3\ 10\ 6\ 11\ 9)(16\ 19\ 18)$,\\
$\pi_{0,18}=(0\ 7\ 1\ 6\ 11\ 2\ 4\ 8\ 9)(12\ 16\ 14\ 15)(13\ 17\ 18\ 20)$,&
$\pi_{7,2}=(1\ 3\ 9\ 6\ 11\ 7\ 8)(14\ 19)$,\\
$\pi_{0,19}=(0\ 11)(2\ 9)(3\ 8)(5\ 10\ 6\ 7)(14\ 15)(16\ 18)(19\ 20)$,&
$\pi_{7,3}=(0\ 9)(2\ 6)(8\ 10\ 11)(12\ 20)(16\ 18)$,\\
$\pi_{0,24}=\pi_{0,36}=\pi_{6,30}=(1)$,&
$\pi_{7,4}=(0\ 8)(2\ 10\ 6)(12\ 17\ 19)$,\\
$\pi_{7,5}=(0\ 6\ 8\ 10\ 1)(14\ 18\ 15\ 16)$,&
$\pi_{7,6}=(0\ 1\ 4)(2\ 8\ 6)(15\ 17)$,\\
$\pi_{7,7}=(0\ 5\ 10\ 4\ 3)(15\ 17)$,&
$\pi_{12,5}=(0\ 1\ 3)(15\ 16)$,\\
$\pi_{12,24}=\pi_{13,11}=\pi_{18,18}=\pi_{24,0}=(1)$,&
$\pi_{24,12}=\pi_{30,6}=\pi_{36,0}=(1)$.
\end{tabular}
\end{center}
\noindent We have that for each $i,j,s,t$ in Table VIII,
 $\pi_{s,t}Y=Y$, $|{\cal B}_i\cap \pi_{s,t}{\cal B}_j|=s$ and
 $|T({\cal B}_i\setminus\pi_{s,t}{\cal B}_j)\cap T(\pi_{s,t}{\cal B}_j\setminus{\cal B}_i)|=t$. \qed

\begin{center}
\small  Table VIII.  \ Fine triangle intersections for $(K_{21}\setminus K_{12},G)$-designs
  \begin{tabular}{ccc|ccc}\hline
$i$ &  $j$ &  $(s,t)$ & $i$ &  $j$ &  $(s,t)$ \\
\hline     1  &   1   &  $(M_{21}\cup N_{21})\setminus \{(0,36),$ & 1 & 2  &  (0,36)   \\
&&(24,0),(0,24),(30,6),&  1 &    3   &   (24,0)    \\
&&(6,30),(24,12),(12,24),& 1  &   4  &   (30,6)    \\
&&(13,11),(18,18)\} & 1  &   5  &   (24,12)    \\
&&&           1  &   6  &   (13,11)    \\
&&&           1  &   7  &   (18,18)    \\
&&&           2  &   3  &   (0,24)   \\
&&&           2  &   4  &   (6,30)    \\
&&&           2  &   5  &   (12,24)    \\
\hline
\end{tabular}
\end{center}

\noindent {\bf Lemma $3.14$} {\em
Let $M_{22}=\{(0,0),(0,1),\ldots,(0,20),(10,0),(10,1)\ldots,(10,6),(13,7),(40,0)\}$ and
$N_{22}=\{(0,28),(0,40),(7,33),(11,17),(14,26),(17,11),(21,19),(26,14),(28,0),(33,7)\}$.
Let $(s,t)$ $\in M_{22}\cup N_{22}$.
Then there is a pair of  kite-GDDs of type $8^{2}6^{1}$ with the same group set,
which intersect in $s$ blocks and $t+s$ triangles.}

\proof Take the vertex set $X=\{0,1,\ldots,21\}$. Let
\begin{center}
\small
 \begin{tabular}{llllll}
${\cal A}_1:$ &$[0,16,8-6]$,&$[1,17,8-7]$,&$[3,18,8-2]$,&$[1,16,9-3]$,&$[4,20,9-6]$,\\
&$[5,19,9-7]$,&$[6,19,10-0]$,&$[21,4,10-5]$,&$[2,16,11-6]$,&$[3,19,11-7]$,\\
&$[4,18,11-5]$,&$[18,6,15-2]$,&$[7,20,15-4]$,&$[19,0,15-3]$,&$[17,5,15-1]$,\\\vspace{0.1cm}
&$[0,11,20-2]$,&$[13,5,20-6]$,&$[1,12,20-8]$,&$[3,10,20-14]$;\\

${\cal A}_2:$ &$[0,16,8-2]$,&$[1,17,8-6]$,&$[3,18,8-7]$,&$[1,16,9-3]$,&$[4,20,9-6]$,\\
&$[5,19,9-7]$,&$[6,19,10-0]$,&$[21,4,10-5]$,&$[2,16,11-6]$,&$[3,19,11-7]$,\\
&$[4,18,11-5]$,&$[18,6,15-4]$,&$[7,20,15-2]$,&$[19,0,15-1]$,&$[17,5,15-3]$,\\\vspace{0.1cm}
&$[0,11,20-2]$,&$[13,5,20-6]$,&$[1,12,20-8]$,&$[3,10,20-14]$;\\

${\cal A}_3:$ &$[0,16,8-2]$,&$[1,17,8-6]$,&$[3,18,8-7]$,&$[1,16,9-7]$,&$[4,20,9-3]$,\\
&$[5,19,9-6]$,&$[6,19,10-0]$,&$[21,4,10-5]$,&$[2,16,11-6]$,&$[3,19,11-7]$,\\
&$[4,18,11-5]$,&$[18,6,15-4]$,&$[7,20,15-2]$,&$[19,0,15-1]$,&$[17,5,15-3]$,\\\vspace{0.1cm}
&$[0,11,20-6]$,&$[13,5,20-2]$,&$[1,12,20-14]$,&$[3,10,20-8]$;\\

${\cal A}_4:$ &$[0,16,8-2]$,&$[1,17,8-6]$,&$[3,18,8-7]$,&$[1,16,9-7]$,&$[4,20,9-3]$,\\
&$[5,19,9-6]$,&$[6,19,10-5]$,&$[21,4,10-0]$,&$[2,16,11-5]$,&$[3,19,11-6]$,\\
&$[4,18,11-7]$,&$[18,6,15-4]$,&$[7,20,15-2]$,&$[19,0,15-1]$,&$[17,5,15-3]$,\\\vspace{0.1cm}
&$[0,11,20-6]$,&$[13,5,20-2]$,&$[1,12,20-14]$,&$[3,10,20-8]$;\\

${\cal A}_5:$ &$[0,16,8-6]$,&$[1,17,8-7]$,&$[3,18,8-2]$,&$[1,16,9-3]$,&$[4,21,9-6]$,\\
&$[5,19,9-7]$,&$[6,19,10-0]$,&$[20,4,10-5]$,&$[2,16,11-6]$,&$[3,19,11-7]$,\\\vspace{0.1cm}
&$[4,18,11-5]$,&$[2,19,12-6]$,&$[4,17,12-0]$,&$[20,7,12-5]$;\\

${\cal A}_6:$ &$[0,16,8-2]$,&$[1,17,8-6]$,&$[3,18,8-7]$,&$[1,16,9-7]$,&$[4,21,9-3]$,\\
&$[5,19,9-6]$,&$[6,19,10-5]$,&$[20,4,10-0]$,&$[2,16,11-5]$,&$[3,19,11-6]$,\\\vspace{0.1cm}
&$[4,18,11-7]$,&$[2,19,12-5]$,&$[4,17,12-6]$,&$[20,7,12-0]$;\\

${\cal A}_7:$ &$[0,16,8-2]$,&$[1,17,8-6]$,&$[3,18,8-7]$,&$[1,16,9-7]$,&$[4,20,9-3]$,\\
&$[5,19,9-6]$,&$[6,19,10-5]$,&$[21,4,10-0]$,&$[2,16,11-5]$,&$[3,19,11-6]$,\\
&$[4,18,11-7]$,&$[2,19,12-5]$,&$[4,17,12-6]$,&$[21,7,12-0]$,&$[3,17,13-2]$,\\
&$[4,16,13-1]$,&$[7,18,13-6]$,&$[5,16,14-7]$,&$[6,17,14-4]$,&$[0,18,14-3]$,\\
&$[18,6,15-4]$,&$[7,20,15-2]$,&$[19,0,15-1]$,&$[17,5,15-3]$,&$[7,10,16-6]$,\\
&$[3,12,16-15]$,&$[0,9,17-11]$,&$[2,10,17-7]$,&$[1,10,18-12]$,&$[2,9,18-5]$,\\
\end{tabular}
\end{center}

\begin{center}
\small
 \begin{tabular}{llllll}

&$[4,8,19-7]$,&$[14,1,19-13]$,&$[0,11,20-6]$,&$[13,5,20-2]$,&$[1,12,20-14]$,\\\vspace{0.1cm}
&$[3,10,20-8]$,&$[0,13,21-9]$,&$[1,11,21-3]$,&$[5,8,21-6]$,&$[14,2,21-15]$;\\

${\cal A}_8:$ &$[2,19,12-6]$,&$[4,17,12-0]$,&$[21,7,12-5]$,&$[3,17,13-1]$,&$[4,16,13-6]$,\\
&$[7,18,13-2]$,&$[5,16,14-4]$,&$[6,17,14-3]$,&$[0,18,14-7]$,&$[7,10,16-15]$,\\
&$[3,12,16-6]$,&$[0,9,17-7]$,&$[2,10,17-11]$,&$[1,10,18-5]$,&$[2,9,18-12]$,\\
&$[4,8,19-13]$,&$[14,1,19-7]$,&$[0,13,21-3]$,&$[1,11,21-9]$,&$[5,8,21-15]$,\\\vspace{0.1cm}
&$[14,2,21-6]$;\\

${\cal A}_9:$ &$[3,17,13-1]$,&$[4,16,13-6]$,&$[7,18,13-2]$,&$[5,16,14-4]$,&$[6,17,14-3]$,\\
&$[0,18,14-7]$,&$[18,6,15-2]$,&$[7,21,15-4]$,&$[19,0,15-3]$,&$[17,5,15-1]$,\\
&$[7,10,16-15]$,&$[3,12,16-6]$,&$[0,9,17-7]$,&$[2,10,17-11]$,&$[1,10,18-5]$,\\
&$[2,9,18-12]$,&$[4,8,19-13]$,&$[14,1,19-7]$,&$[0,11,21-2]$,&$[13,5,21-6]$,\\
&$[1,12,21-8]$,&$[3,10,21-14]$,&$[0,13,20-3]$,&$[1,11,20-9]$,&$[5,8,20-15]$,\\\vspace{0.1cm}
&$[14,2,20-6]$.\\
\end{tabular}
\end{center}

\noindent Let ${\cal B}_{i}={\cal A}_{i}\cup {\cal A}_{8}$ when $1\leq i\leq 4$. Let ${\cal B}_{i}={\cal A}_{i}\cup {\cal A}_{9}$ when $i=5,6$. Let ${\cal B}_{7}={\cal A}_{7}$.
Then $(X,{\cal G},{\cal B}_{i})$ is a kite-GDDs of type $8^{2}6^{1}$ for each $1\leq i\leq 7$, where ${\cal G}=\{G_1,G_2,G_3\}=\{\{0,1,\ldots,7\},\{8,9,\ldots,15\},\{16,17,\ldots,21\}\}$.
Now take permutation $\pi_{s,t}=(1)$ if $(s,t)\in N_{22}$. For $(s,t)\in M_{22}$, we take

\begin{center}
\small
\begin{tabular}{ll}
$\pi_{0,0}=(0\ 1)(2\ 3)(4\ 5)(6\ 7)$,&
$\pi_{0,1}=(2\ 3)(4\ 5\ 6\ 7)(16\ 17)(18\ 19\ 20\ 21)$,\\
$\pi_{0,2}=(2\ 3)(4\ 7)(5\ 6)(16\ 17)(18\ 21)(19\ 20)$,&
$\pi_{0,3}=(3\ 4\ 6\ 5\ 7)(17\ 18\ 20\ 19\ 21)$,\\
$\pi_{0,4}=(3\ 4\ 7)(5\ 6)(17\ 18\ 21)(19\ 20)$,&
$\pi_{0,5}=(3\ 4)(5\ 6\ 7)(17\ 18)(19\ 20\ 21)$,\\
$\pi_{0,6}=(3\ 4\ 5\ 6\ 7)(17\ 18\ 19\ 20\ 21)$,&
$\pi_{0,7}=(3\ 4\ 6\ 7\ 5)(17\ 18\ 20\ 21\ 19)$,\\
$\pi_{0,8}=(3\ 4\ 7\ 6)(17\ 18\ 21\ 20)$,&
$\pi_{0,9}=(2\ 4\ 6)(3\ 7\ 5)(16\ 18\ 20)(17\ 21\ 19)$,\\
$\pi_{0,10}=(3\ 4)(6\ 7)(17\ 18)(20\ 21)$,&
$\pi_{0,11}=(0\ 1\ 4\ 2\ 6\ 7\ 3\ 5)(16\ 20\ 21\ 18\ 17\ 19)$,\\
$\pi_{0,12}=(0\ 1\ 3\ 4\ 5\ 6)(16\ 17\ 18\ 19\ 20)$,&
$\pi_{0,13}=(0\ 2\ 1\ 3\ 4\ 5\ 6)(16\ 17\ 18\ 19\ 20)$,\\
$\pi_{0,14}=(9\ 11\ 14\ 12)(10\ 15\ 13)(16\ 21\ 19\ 17\ 20\ 18)$,&
$\pi_{0,15}=(2\ 4\ 6)(3\ 7)(10\ 12\ 14)(11\ 15)$,\\
$\pi_{0,16}=(0\ 1)(2\ 6)(5\ 7)(8\ 9)(10\ 14)(13\ 15)$,&
$\pi_{0,17}=(0\ 1)(2\ 4\ 3\ 5\ 7\ 6)(8\ 9)(10\ 12\ 11\ 13\ 15\ 14)$,\\
$\pi_{0,18}=(0\ 7\ 4\ 1\ 2\ 5\ 3)(8\ 13)(9\ 10\ 14\ 12)$&
$\pi_{0,19}=(0\ 5\ 3\ 4\ 1)(2\ 7)(8\ 11\ 15\ 10\ 14)$\\
\hspace{1.2cm}$(16\ 17)(18\ 21\ 20)$,&\hspace{1.2cm}$(9\ 13)(16\ 21\ 19)(17\ 20)$,\\
$\pi_{0,20}=(0\ 1\ 5\ 7\ 6\ 2\ 4)(8\ 14\ 12\ 15\ 10\ 13\ 11)$&
$\pi_{10,0}=(2\ 4)(8\ 15\ 13)(16\ 20)$,\\
\hspace{1.2cm}$(16\ 19\ 18\ 20\ 17)$,&$\pi_{10,1}=(0\ 7\ 1)(10\ 13)(17\ 19)$,\\
$\pi_{10,2}=(8\ 12\ 11\ 10\ 9)(20\ 21)$,&
$\pi_{10,3}=(2\ 5)(8\ 12)(16\ 18\ 21)$,\\
$\pi_{10,4}=(1\ 7)(9\ 12\ 10)(18\ 19)$,&
$\pi_{10,5}=(1\ 5\ 3)(11\ 15\ 14)$,\\
$\pi_{10,6}=(3\ 7)(8\ 10\ 15\ 14)$,&
$\pi_{13,7}=(1\ 3\ 7)(10\ 15)$.
\end{tabular}
\end{center}
\noindent We have that for each $i,j,s,t$ in Table IX,
 $\pi_{s,t}G_{l}=G_{l}$ for $1\leq l\leq 3$, $|{\cal B}_i\cap \pi_{s,t}{\cal B}_j|=s$ and
 $|T({\cal B}_i\setminus\pi_{s,t}{\cal B}_j)\cap T(\pi_{s,t}{\cal B}_j\setminus{\cal B}_i)|=t$. \qed

\begin{center}
\small  Table IX.  \ Fine triangle intersections for kite-GDDs of type $8^{2}6^{1}$

\begin{tabular}{ccc|ccc|ccc}\hline
$i$ &  $j$ &  $(s,t)$ & $i$ &  $j$ &  $(s,t)$&$i$ &  $j$ &  $(s,t)$ \\
\hline     1  &   1   &  $M_{22}$ & 1 &   2&  (33,7)& 1 &  3  & (26,14)  \\
 1 &  4   & (21,19)& 1&   5&  (28,0)  &1&   6&  (17,11)\\
 1&   7&  (0,40)&  7&   2& (7,33)   &7&   3& (14,26)\\
 7 &  5 &  (0,28)  &  7 &  6  &  (11,17)   \\
\hline
\end{tabular}
\end{center}

\noindent {\bf Lemma $3.15$} {\em
Let $M_{23}=\{(0,0),(0,1),\ldots,(0,18),(7,0),(7,6),(7,7),(7,10),(10,6),(52,0)\}$ and
$N_{23}=\{(0,22),(0,26),(0,34),(0,42),(0,52),(3,49),(5,21),(5,29),(5,37),(6,46),
(7,19),(9,43),\\(10,24),(10,32),(12,14),(12,40),(14,0),(14,12),
(15,19),(15,27),(15,37),(18,8),(18,34),(19,\\15),(20,22),(21,5),(21,31),(22,0),
(22,20),(24,10),(24,28),(26,0),(27,15),(27,25),(28,24),(29,\\5),(31,21),(32,10),(34,0),(34,18),
(37,5),(37,15),(40,12),(42,0),(43,9),(46,6),(49,3)\}$. Let $G$ be a kite and $(s,t)\in M_{23}\cup N_{23}$. Then there is a pair of  $(K_{23}\setminus K_{10},G)$-designs with the same vertex set and the same subgraph $K_{10}$ removed, which intersect in $s$ blocks and $t+s$ triangles.}

\proof Take the vertex set $X=\{0,1,\ldots,22\}$. Let

\begin{center}
\small
\begin{tabular}{llllll}
${\cal A}_{1}:$ &$[0,22,12-17]$,&$[22,21,11-8]$,&$[21,20,10-8]$,&$[15,14,4-10]$,&$[3,13,14-21]$,\\
&$[12, 2, 13 -22]$,&$[1 ,11, 12 -21]$,&$[0 ,10 ,11- 7]$,&$[22 ,9 ,10 -7]$,&$[0 ,13, 21 -1]$,\\
&$[15, 1 ,13 -18]$,&$[12, 6, 10 -5]$,&$[1 ,22 ,14 -11]$,&$[5 ,18, 21 -2]$,&$[20 ,17 ,4 -21]$,\\
&$[19 ,3 ,16- 0]$,&$[20 ,1, 16- 10]$,&$[10, 18, 14- 9]$,&$[17 ,13, 10- 3]$,&$[8 ,12 ,16 -11]$,\\\vspace{0.1cm}&$[1 ,19, 10 -2]$,&$[19, 12, 4 -22]$,&$[6 ,15, 21- 3]$,&$[20 ,5, 13 -9]$;\\

${\cal A}_{2}:$ &$[0 ,22, 12- 17]$,&$[22 ,21 ,11- 8]$,&$[21, 20, 10- 8]$,&$[15, 14 ,4- 10]$,&$[3 ,13, 14- 21]$,\\&$[12, 2, 13 -18]$,
&$[1 ,11, 12 -21]$,&$[0 ,10 ,11- 7]$,&$[22 ,9 ,10 -7]$,&$[0 ,13, 21 -1]$,\\&$[15, 1 ,13 -9]$,&$[12, 6, 10 -5]$,
&$[1 ,22 ,14 -11]$,&$[5 ,18, 21 -2]$,&$[20 ,17 ,4 -21]$,\\&$[19 ,3 ,16- 0]$,&$[20 ,1, 16- 10]$,&$[10, 18, 14- 9]$,
&$[17 ,13, 10- 3]$,&$[8 ,12 ,16 -11]$,\\\vspace{0.1cm}&$[1 ,19, 10 -2]$,&$[19, 12, 4 -22]$,&$[6 ,15, 21- 3]$,&$[20 ,5, 13 -22]$;\\

${\cal A}_{3}:$ &$[0 ,22, 12- 17]$,&$[22 ,21 ,11- 8]$,&$[21, 20, 10- 8]$,&$[15, 14 ,4- 10]$,&$[3 ,13, 14- 9]$,\\&$[12, 2, 13 -18]$,
&$[1 ,11, 12 -21]$,&$[0 ,10 ,11- 7]$,&$[22 ,9 ,10 -7]$,&$[0 ,13, 21 -1]$,\\&$[15, 1 ,13 -9]$,&$[12, 6, 10 -5]$,
&$[1 ,22 ,14 -21]$,&$[5 ,18, 21 -2]$,&$[20 ,17 ,4 -21]$,\\&$[19 ,3 ,16- 0]$,&$[20 ,1, 16- 10]$,&$[10, 18, 14- 11]$,
&$[17 ,13, 10- 3]$,&$[8 ,12 ,16 -11]$,\\\vspace{0.1cm}&$[1 ,19, 10 -2]$,&$[19, 12, 4 -22]$,&$[6 ,15, 21- 3]$,&$[20 ,5, 13 -22]$;\\

${\cal A}_{4}:$ &$[0 ,22, 12- 17]$,&$[22 ,21 ,11- 8]$,&$[21, 20, 10- 8]$,&$[15, 14 ,4- 10]$,&$[3 ,13, 14- 9]$,\\&$[12, 2, 13 -18]$,
&$[1 ,11, 12 -21]$,&$[0 ,10 ,11- 7]$,&$[22 ,9 ,10 -7]$,&$[0 ,13, 21 -3]$,\\&$[15, 1 ,13 -9]$,&$[12, 6, 10 -5]$,
&$[1 ,22 ,14 -21]$,&$[5 ,18, 21 -1]$,&$[20 ,17 ,4 -21]$,\\&$[19 ,3 ,16- 0]$,&$[20 ,1, 16- 10]$,&$[10, 18, 14- 11]$,
&$[17 ,13, 10- 3]$,&$[8 ,12 ,16 -11]$,\\\vspace{0.1cm}&$[1 ,19, 10 -2]$,&$[19, 12, 4 -22]$,&$[6 ,15, 21- 2]$,&$[20 ,5, 13 -22]$;\\

${\cal A}_{5}:$ &$[0 ,22, 12- 17]$,&$[22 ,21 ,11- 8]$,&$[21, 20, 10- 8]$,&$[15, 14 ,4- 22]$,&$[3 ,13, 14- 9]$,\\&$[12, 2, 13 -18]$,
&$[1 ,11, 12 -21]$,&$[0 ,10 ,11- 7]$,&$[22 ,9 ,10 -7]$,&$[0 ,13, 21 -3]$,\\&$[15, 1 ,13 -9]$,&$[12, 6, 10 -5]$,
&$[1 ,22 ,14 -21]$,&$[5 ,18, 21 -1]$,&$[20 ,17 ,4 -10]$,\\&$[19 ,3 ,16- 0]$,&$[20 ,1, 16- 10]$,&$[10, 18, 14- 11]$,
&$[17 ,13, 10- 3]$,&$[8 ,12 ,16 -11]$,\\\vspace{0.1cm}&$[1 ,19, 10 -2]$,&$[19, 12, 4 -21]$,&$[6 ,15, 21- 2]$,&$[20 ,5, 13 -22]$;\\

${\cal A}_{6}:$ &$[0 ,22, 12- 17]$,&$[22 ,21 ,11- 8]$,&$[21, 20, 10- 8]$,&$[15, 14 ,4- 22]$,&$[3 ,13, 14- 9]$,\\&$[12, 2, 13 -18]$,
&$[1 ,11, 12 -21]$,&$[0 ,10 ,11- 7]$,&$[22 ,9 ,10 -7]$,&$[0 ,13, 21 -3]$,\\&$[15, 1 ,13 -9]$,&$[12, 6, 10 -5]$,
&$[1 ,22 ,14 -21]$,&$[5 ,18, 21 -1]$,&$[20 ,17 ,4 -10]$,\\&$[19 ,3 ,16- 11]$,&$[20 ,1, 16- 0]$,&$[10, 18, 14- 11]$,
&$[17 ,13, 10- 3]$,&$[8 ,12 ,16 -10]$,\\\vspace{0.1cm}&$[1 ,19, 10 -2]$,&$[19, 12, 4 -21]$,&$[6 ,15, 21- 2]$,&$[20 ,5, 13 -22]$;\\

${\cal A}_{7}:$ &$[0 ,22, 12- 17]$,&$[22 ,21 ,11- 8]$,&$[21, 20, 10- 8]$,&$[15, 14 ,4- 22]$,&$[3 ,13, 14- 9]$,\\&$[12, 2, 13 -18]$,
&$[1 ,11, 12 -21]$,&$[0 ,10 ,11- 7]$,&$[22 ,9 ,10 -7]$,&$[0 ,13, 21 -3]$,\\&$[15, 1 ,13 -9]$,&$[12, 6, 10 -2]$,
&$[1 ,22 ,14 -21]$,&$[5 ,18, 21 -1]$,&$[20 ,17 ,4 -10]$,\\&$[19 ,3 ,16- 11]$,&$[20 ,1, 16- 0]$,&$[10, 18, 14- 11]$,
&$[17 ,13, 10- 5]$,&$[8 ,12 ,16 -10]$,\\\vspace{0.1cm}&$[1 ,19, 10 -3]$,&$[19, 12, 4 -21]$,&$[6 ,15, 21- 2]$,&$[20 ,5, 13 -22]$;\\

${\cal A}_{8}:$ &$[0 ,22, 12- 21]$,&$[22 ,21 ,11- 7]$,&$[21, 20, 10- 7]$,&$[15, 14 ,4- 22]$,&$[3 ,13, 14- 9]$,\\&$[12, 2, 13 -18]$,
&$[1 ,11, 12 -17]$,&$[0 ,10 ,11- 8]$,&$[22 ,9 ,10 -8]$,&$[0 ,13, 21 -3]$,\\&$[15, 1 ,13 -9]$,&$[12, 6, 10 -5]$,
&$[1 ,22 ,14 -21]$,&$[5 ,18, 21 -1]$,&$[20 ,17 ,4 -10]$,\\&$[19 ,3 ,16- 11]$,&$[20 ,1, 16- 0]$,&$[10, 18, 14- 11]$,
&$[17 ,13, 10- 3]$,&$[8 ,12 ,16 -10]$,\\\vspace{0.1cm}&$[1 ,19, 10 -2]$,&$[19, 12, 4 -21]$,&$[6 ,15, 21- 2]$,&$[20 ,5, 13 -22]$;\\

${\cal A}_{9}:$ &$[0 ,22, 12- 21]$,&$[22 ,21 ,11- 7]$,&$[21, 20, 10- 7]$,&$[15, 14 ,4- 22]$,&$[3 ,13, 14- 9]$,\\&$[12, 2, 13 -18]$,
&$[1 ,11, 12 -17]$,&$[0 ,10 ,11- 8]$,&$[22 ,9 ,10 -8]$,&$[0 ,13, 21 -3]$,\\&$[15, 1 ,13 -9]$,&$[12, 6, 10 -2]$,
&$[1 ,22 ,14 -21]$,&$[5 ,18, 21 -1]$,&$[20 ,17 ,4 -10]$,\\
&$[19 ,3 ,16- 11]$,&$[20 ,1, 16- 0]$,&$[10, 18, 14- 11]$,
&$[17 ,13, 10- 5]$,&$[8 ,12 ,16 -10]$,\\
\end{tabular}
\end{center}

\begin{center}
\small
\begin{tabular}{llllll}
\vspace{0.1cm}&$[1 ,19, 10 -3]$,&$[19, 12, 4 -21]$,&$[6 ,15, 21- 2]$,&$[20 ,5, 13 -22]$;\\

${\cal A}_{10}:$ &$[16 ,17 ,6- 13]$,&$[16, 5 ,15 -10]$,&$[22 ,8, 20- 7]$,&$[21 ,19 ,7- 14]$,&$[20 ,18 ,6 -14]$,\\&$[5, 17 ,19- 13]$,
&$[18 ,4 ,16- 21]$,&$[3 ,15, 17- 1]$,&$[2 ,14, 16 -22]$,&$[14 ,5 ,12- 18]$,\\&$[13 ,4, 11- 5]$,&$[12 ,6, 10 -2]$,
&$[18 ,9 ,11 -6]$,&$[0 ,14, 20 -2]$,&$[6 ,19, 22- 5]$,\\&$[16, 13, 7- 22]$,&$[7, 12, 15 -8]$,&$[15 ,0, 19 -14]$,
&$[3 ,18, 22- 15]$,&$[17, 13 ,10- 5]$,\\\vspace{0.1cm}&$[15 ,20 ,11 -17]$,&$[1 ,19, 10 -3]$,&$[2 ,22, 17- 0]$,&$[2 ,19, 11- 3]$,&$[20, 3, 12- 9]$;\\

${\cal A}_{11}:$ &$[16 ,17 ,6- 14]$,&$[16, 5 ,15 -8]$,&$[22 ,8, 20- 2]$,&$[21 ,19 ,7- 22]$,&$[20 ,18 ,6 -13]$,\\&$[5, 17 ,19- 14]$,
&$[18 ,4 ,16- 22]$,&$[3 ,15, 17- 0]$,&$[2 ,14, 16 -21]$,&$[14 ,5 ,12- 9]$,\\&$[13 ,4, 11- 6]$,&$[12 ,6, 10 -5]$,
&$[18 ,9 ,11 -5]$,&$[0 ,14, 20 -7]$,&$[6 ,19, 22- 15]$,\\&$[16, 13, 7- 14]$,&$[7, 12, 15 -10]$,&$[15 ,0, 19 -13]$,
&$[3 ,18, 22- 5]$,&$[17, 13 ,10- 3]$,\\\vspace{0.1cm}&$[15 ,20 ,11 -3]$,&$[1 ,19, 10 -2]$,&$[2 ,22, 17- 1]$,&$[2 ,19, 11- 17]$,&$[20, 3, 12- 18]$;\\

 ${\cal A}_{12}:$ &$[0 ,22, 12- 17]$,&$[22 ,21, 11- 8]$,&$[21, 20, 10- 8]$,&$[20, 19 ,9- 15]$,&$[3 ,13, 14- 21]$,\\&$[12, 2, 13- 22]$,
&$[1, 11, 12- 21]$,&$[0 ,10, 11- 7]$,&$[22, 9 ,10- 7]$,&$[22 ,8, 20 -2]$,\\&$[18, 6 ,16- 22]$,&$[2 ,14 ,16 -21]$,
&$[15 ,1, 13 -18]$,&$[1 ,22, 14- 11]$,&$[0 ,14 ,20- 7]$,\\&$[19 ,3, 16 -0]$,&$[21 ,17, 9- 16]$,&$[20 ,1 ,16- 10]$,
&$[10, 18 ,14- 9]$,&$[8 ,12, 16 -11]$,\\\vspace{0.1cm}&$[20, 5 ,13 -9]$;\\

${\cal A}_{13}:$ &$[0 ,22, 12- 21]$,&$[22 ,21, 11- 7]$,&$[21, 20, 10- 7]$,&$[20, 19 ,9- 16]$,&$[3 ,13, 14- 9]$,\\&$[12, 2, 13- 18]$,
&$[1, 11, 12- 17]$,&$[0 ,10, 11- 8]$,&$[22, 9 ,10- 8]$,&$[22 ,8, 20 -7]$,\\&$[18, 6 ,16- 21]$,&$[2 ,14 ,16 -22]$,
&$[15 ,1, 13 -9]$,&$[1 ,22, 14- 21]$,&$[0 ,14 ,20- 2]$,\\&$[19 ,3, 16 -11]$,&$[21 ,17, 9- 15]$,&$[20 ,1 ,16- 0]$,
&$[10, 18 ,14- 11]$,&$[8 ,12, 16 -10]$,\\\vspace{0.1cm}&$[20, 5 ,13 -22]$;\\

${\cal A}_{14}:$ &$[0 ,22, 12- 21]$,&$[22 ,21, 11- 8]$,&$[21, 20, 10- 8]$,&$[20, 19 ,9- 15]$,&$[3 ,13, 14- 9]$,\\&$[12, 2, 13- 22]$,
&$[1, 11, 12- 17]$,&$[0 ,10, 11- 7]$,&$[22, 9 ,10- 7]$,&$[22 ,8, 20 -2]$,\\&$[18, 6 ,16- 22]$,&$[2 ,14 ,16 -21]$,
&$[15 ,1, 13 -18]$,&$[1 ,22, 14- 21]$,&$[0 ,14 ,20- 7]$,\\&$[19 ,3, 16 -0]$,&$[21 ,17, 9- 16]$,&$[20 ,1 ,16- 10]$,
&$[10, 18 ,14- 11]$,&$[8 ,12, 16 -11]$,\\\vspace{0.1cm}&$[20, 5 ,13 -9]$;\\

${\cal A}_{15}:$ &$[0 ,22, 12- 21]$,&$[22 ,21, 11- 7]$,&$[21, 20, 10- 8]$,&$[20, 19 ,9- 15]$,&$[3 ,13, 14- 9]$,\\&$[12, 2, 13- 18]$,
&$[1, 11, 12- 17]$,&$[0 ,10, 11- 8]$,&$[22, 9 ,10- 7]$,&$[22 ,8, 20 -7]$,\\&$[18, 6 ,16- 22]$,&$[2 ,14 ,16 -21]$,
&$[15 ,1, 13 -9]$,&$[1 ,22, 14- 21]$,&$[0 ,14 ,20- 2]$,\\&$[19 ,3, 16 -11]$,&$[21 ,17, 9- 16]$,&$[20 ,1 ,16- 0]$,
&$[10, 18 ,14- 11]$,&$[8 ,12, 16 -10]$,\\\vspace{0.1cm}&$[20, 5 ,13 -22]$;\\

${\cal A}_{16}:$ &$[0 ,22, 12- 21]$,&$[22 ,21, 11- 7]$,&$[21, 20, 10- 8]$,&$[20, 19 ,9- 15]$,&$[3 ,13, 14- 9]$,\\&$[12, 2, 13- 18]$,
&$[1, 11, 12- 17]$,&$[0 ,10, 11- 8]$,&$[22, 9 ,10- 7]$,&$[22 ,8, 20 -2]$,\\&$[18, 6 ,16- 22]$,&$[2 ,14 ,16 -21]$,
&$[15 ,1, 13 -9]$,&$[1 ,22, 14- 21]$,&$[0 ,14 ,20- 7]$,\\&$[19 ,3, 16 -0]$,&$[21 ,17, 9- 16]$,&$[20 ,1 ,16- 10]$,
&$[10, 18 ,14- 11]$,&$[8 ,12, 16 -11]$,\\\vspace{0.1cm}&$[20, 5 ,13 -22]$;\\

${\cal A}_{17}:$ &$[22 ,21 ,11 -7]$,&$[21 ,20, 10 -7]$,&$[19, 12, 8- 21]$,&$[16, 17, 6- 13]$,&$[0 ,10, 11- 8]$,\\&$[22, 9 ,10 -8]$,
&$[22 ,8, 20- 7]$,&$[21 ,19, 7- 14]$,&$[20, 12, 6 -14]$,&$[5 ,17 ,19 -13]$,\\&$[13 ,4 ,11- 5]$,&$[12 ,9 ,11- 6]$,
&$[0 ,14 ,20- 2]$,&$[17 ,14, 8 -13]$,&$[16, 13, 7- 22]$,\\\vspace{0.1cm}&$[15, 0 ,19 -14]$,&$[15 ,20, 11 -17]$,&$[2 ,19, 11- 3]$;\\

${\cal A}_{18}:$ &$[22 ,21 ,11 -8]$,&$[21 ,20, 10 -8]$,&$[19, 12, 8- 13]$,&$[16, 17, 6- 14]$,&$[0 ,10, 11- 7]$,\\&$[22, 9 ,10 -7]$,
&$[22 ,8, 20- 2]$,&$[21 ,19, 7- 22]$,&$[20, 12, 6 -13]$,&$[5 ,17 ,19 -14]$,\\&$[13 ,4 ,11- 6]$,&$[12 ,9 ,11- 5]$,
&$[0 ,14 ,20- 7]$,&$[17 ,14, 8 -21]$,&$[16, 13, 7- 14]$,\\\vspace{0.1cm}&$[15, 0 ,19 -13]$,&$[15 ,20, 11 -3]$,&$[2 ,19, 11- 17]$;\\

${\cal A}_{19}:$ &$[22 ,21 ,11 -8]$,&$[21 ,20, 10 -8]$,&$[19, 12, 8- 13]$,&$[16, 17, 6- 13]$,&$[0 ,10, 11- 7]$,\\&$[22, 9 ,10 -7]$,
&$[22 ,8, 20- 2]$,&$[21 ,19, 7- 14]$,&$[20, 12, 6 -14]$,&$[5 ,17 ,19 -13]$,\\&$[13 ,4 ,11- 6]$,&$[12 ,9 ,11- 5]$,
&$[0 ,14 ,20- 7]$,&$[17 ,14, 8 -21]$,&$[16, 13, 7- 22]$,\\\vspace{0.1cm}&$[15, 0 ,19 -14]$,&$[15 ,20, 11 -3]$,&$[2 ,19, 11- 17]$;\\

${\cal A}_{20}:$ &$[22 ,21 ,11 -8]$,&$[21 ,20, 10 -7]$,&$[19, 12, 8- 21]$,&$[16, 17, 6- 13]$,&$[0 ,10, 11- 7]$,\\&$[22, 9 ,10 -8]$,
&$[22 ,8, 20- 7]$,&$[21 ,19, 7- 14]$,&$[20, 12, 6 -14]$,&$[5 ,17 ,19 -13]$,\\&$[13 ,4 ,11- 6]$,&$[12 ,9 ,11- 5]$,
&$[0 ,14 ,20- 2]$,&$[17 ,14, 8 -13]$,&$[16, 13, 7- 22]$,\\\vspace{0.1cm}&$[15, 0 ,19 -14]$,&$[15 ,20, 11 -3]$,&$[2 ,19, 11- 17]$;\\

${\cal A}_{21}:$ &$[0 ,22 ,12 -17]$,&$[22 ,21, 11- 8]$,&$[21 ,15 ,10- 8]$,&$[20 ,18 ,8 -13]$,&$[1 ,11, 12- 21]$,\\&$[0 ,10, 11- 7]$,
&$[22, 9, 10- 7]$,&$[18 ,4 ,16 -22]$,&$[2 ,14 ,16- 21]$,&$[14, 5, 12- 9]$,\\&$[13, 4 ,11- 6]$,&$[12 ,6, 10- 5]$,
&$[18, 9 ,11- 5]$,&$[17 ,14 ,8 -21]$,&$[17, 13 ,10- 3]$,\\\vspace{0.1cm}&$[1 ,20, 10- 2]$,&$[15, 3 ,12- 18]$;\\
\end{tabular}
\end{center}

\begin{center}
\small
 \begin{tabular}{llllll}
${\cal A}_{22}:$ &$[0 ,22 ,12 -21]$,&$[22 ,21, 11- 7]$,&$[21 ,15 ,10- 8]$,&$[20 ,18 ,8 -21]$,&$[1 ,11, 12- 17]$,\\&$[0 ,10, 11- 8]$,
&$[22, 9, 10- 7]$,&$[18 ,4 ,16 -21]$,&$[2 ,14 ,16- 22]$,&$[14, 5, 12- 18]$,\\&$[13, 4 ,11- 5]$,&$[12 ,6, 10- 2]$,
&$[18, 9 ,11- 6]$,&$[17 ,14 ,8 -13]$,&$[17, 13 ,10- 5]$,\\\vspace{0.1cm}&$[1 ,20, 10- 3]$,&$[15, 3 ,12- 9]$;\\

${\cal A}_{23}:$ &$[0 ,22 ,12 -21]$,&$[22 ,21, 11- 7]$,&$[21 ,15 ,10- 8]$,&$[20 ,18 ,8 -13]$,&$[1 ,11, 12- 17]$,\\
&$[0 ,10, 11- 8]$,&$[22, 9, 10- 7]$,&$[18 ,4 ,16 -21]$,&$[2 ,14 ,16- 22]$,&$[14, 5, 12- 18]$,\\&$[13, 4 ,11- 6]$,&$[12 ,6, 10- 5]$,
&$[18, 9 ,11- 5]$,&$[17 ,14 ,8 -21]$,&$[17, 13 ,10- 3]$,\\\vspace{0.1cm}&$[1 ,20, 10- 2]$,&$[15, 3 ,12- 9]$;\\

${\cal A}_{24}:$ &$[0 ,22 ,12 -21]$,&$[22 ,21, 11- 7]$,&$[21 ,15 ,10- 8]$,&$[20 ,18 ,8 -13]$,&$[1 ,11, 12- 17]$,\\
&$[0 ,10, 11- 8]$,&$[22, 9, 10- 7]$,&$[18 ,4 ,16 -21]$,&$[2 ,14 ,16- 22]$,&$[14, 5, 12- 9]$,\\&$[13, 4 ,11- 5]$,&$[12 ,6, 10- 5]$,
&$[18, 9 ,11- 6]$,&$[17 ,14 ,8 -21]$,&$[17, 13 ,10- 3]$,\\\vspace{0.1cm}&$[1 ,20, 10- 2]$,&$[15, 3 ,12- 18]$;\\

${\cal A}_{25}:$ &$[0 ,22 ,12 -21]$,&$[22 ,21, 11- 8]$,&$[21 ,15 ,10- 7]$,&$[20 ,18 ,8 -13]$,&$[1 ,11, 12- 17]$,\\
&$[0 ,10, 11- 7]$,&$[22, 9, 10- 8]$,&$[18 ,4 ,16 -21]$,&$[2 ,14 ,16- 22]$,&$[14, 5, 12- 9]$,\\&$[13, 4 ,11- 6]$,&$[12 ,6, 10- 5]$,
&$[18, 9 ,11- 5]$,&$[17 ,14 ,8 -21]$,&$[17, 13 ,10- 3]$,\\\vspace{0.1cm}&$[1 ,20, 10- 2]$,&$[15, 3 ,12- 18]$;\\

${\cal A}_{26}:$ &$[0 ,22, 12- 17]$,&$[20 ,18, 8- 13]$,&$[7 ,17, 18 -1]$,&$[16 ,17, 4 -14]$,&$[19, 14 ,6 -10]$,\\
&$[12, 2, 13- 22]$,&$[1 ,11 ,12 -21]$,&$[15 ,18, 4 -13]$,&$[18 ,6 ,16- 22]$,&$[19, 1 ,13- 18]$,\\&$[14, 5 ,12- 9]$,&$[13, 6, 11- 4]$,
&$[12 ,4 ,10- 5]$,&$[18, 9 ,11- 5]$,&$[4 ,20, 22- 19]$,\\&$[5 ,18, 21- 2]$,&$[15 ,17 ,6- 21]$,&$[2 ,19 ,18 -0]$,
&$[7 ,12 ,19- 10]$,&$[3 ,18, 22 -5]$,\\\vspace{0.1cm}&$[10 ,18, 14 -9]$,&$[8 ,12 ,16 -11 ]$,&$[15, 3 ,12 -18]$,&$[20 ,12 ,6- 22]$,
&$[4 ,19, 21- 3]$;\\

${\cal A}_{27}:$ &$[0 ,22, 18- 17]$,&$[20, 12, 8- 13]$,&$[7 ,17 ,12- 1]$,&$[16 ,17, 6 -14]$,&$[19 ,14, 4- 10]$,\\
&$[18 ,2 ,13 -22]$,&$[1 ,11, 18- 21]$,&$[15 ,12, 6 -13]$,&$[12 ,4 ,16 -22]$,&$[19 ,1 ,13 -12]$,\\&$[14, 5 ,18 -9]$,&$[13, 4 ,11- 6]$,
&$[18 ,6, 10- 5]$,&$[12, 9 ,11- 5]$,&$[6 ,20 ,22- 19]$,\\&$[5 ,12, 21- 2]$,&$[15 ,17 ,4 -21]$,&$[2 ,19, 12 -0]$,
&$[7 ,18, 19 -10]$,&$[3 ,12, 22- 5]$,\\\vspace{0.1cm}&$[10, 12 ,14- 9]$,&$[8, 18, 16 -11]$,&$[15 ,3, 18 -12]$,&$[20 ,18, 4- 22]$,
&$[6 ,19 ,21 -3]$;\\

${\cal A}_{28}:$ &$[20 ,19, 9 -15]$,&$[ 19 ,18, 8- 13]$,&$[7 ,17 ,18 -1]$,&$[16, 17 ,6 -14]$,&$[16 ,5 ,15 -8]$,\\
&$[22 ,8 ,20- 2]$,&$[21 ,19 ,7- 22]$,&$[20 ,18 ,6- 13]$,&$[5, 17, 19 -14]$,&$[18, 4 ,16- 22]$,\\&$[3, 15 ,17- 0]$,&$[2 ,14 ,16- 21]$,
&$[ 14 ,5, 12- 9]$,&$[13, 4, 11- 6]$,&$[18 ,9, 11- 5]$,\\&$[0 ,14 ,20 -7]$,&$[6, 19, 22- 15]$,&$[2 ,15, 18 -0]$,
&$[17, 14, 8 -21]$,&$[16, 13, 7 -14]$,\\&$[7 ,12 ,15 -10]$,&$[15 ,0, 19 -13]$,&$[3 ,18 ,22 -5]$,&$[21 ,17, 9- 16]$,
&$[15 ,20, 11 -3]$,\\\vspace{0.1cm}&$[2 ,22, 17- 1]$,&$[2 ,19 ,11- 17]$,&$[20, 3, 12- 18]$;\\

${\cal A}_{29}:$ &$[0 ,22, 12- 21]$,&$[ 22 ,21 ,11- 7]$,&$[21 ,20, 10- 7]$,&$[20 ,19, 9 -16]$,&$[19, 18, 8- 21]$,\\
&$[ 7 ,17, 18- 0]$,&$[ 15 ,14 ,4 -22]$,&$[3 ,13 ,14- 9]$,&$[12, 2 ,13- 18]$,&$[1 ,11 ,12- 17]$,\\&$[0 ,10, 11- 8]$,&$[22, 9, 10- 8]$,
&$[ 0 ,13 ,21- 3]$,&$[ 15, 1 ,13- 9]$,&$[1 ,22 ,14- 21]$,\\&$[5 ,18 ,21- 1]$,&$[20 ,17, 4- 10]$,&$[19 ,3, 16 -11]$,
&$[2 ,15, 18 -1]$,&$[17 ,14, 8 -13]$,\\&$[21, 17 ,9 -15]$,&$[20, 1 ,16- 0]$,&$[ 10, 18 ,14- 11]$,&$[8 ,12 ,16 -10]$,
&$[19, 12 ,4 -21]$,\\\vspace{0.1cm}&$[6 ,15, 21 -2]$,&$[20 ,5, 13- 22]$;\\

${\cal A}_{30}:$ &$[19 ,18, 8- 13]$,&$[7 ,17, 18 -1]$,&$[ 16 ,17, 4- 14]$,&$[16 ,5, 15- 8]$,&$[15, 14 ,6- 10]$,\\
&$[0 ,13, 21- 1]$,&$[21, 19, 7- 22]$,&$[20, 18, 4- 13]$,&$[5 ,17, 19- 14]$,&$[ 3 ,15 ,17 -0]$,\\&$[14 ,5, 12- 9]$,&$[13 ,6, 11- 4]$,
&$[ 12, 4, 10- 5]$,&$[18 ,9, 11- 5]$,&$[4 ,19 ,22 -15]$,\\&$[5 ,18, 21 -2]$,&$[20, 17 ,6 -21]$,&$[2, 15, 18- 0]$,
&$[17, 14, 8- 21]$,&$[16 ,13 ,7- 14]$,\\&$[7 ,12 ,15- 10]$,&$[15 ,0 ,19- 13]$,&$[3 ,18, 22 -5]$,&$[17, 13, 10- 3]$,
&$[15 ,20 ,11- 3]$,\\&$[1 ,19 ,10- 2]$,&$[2 ,22, 17- 1]$,&$[2, 19, 11 -17]$,&$[20 ,3 ,12- 18]$,&$[19 ,12, 6 -22]$,\\\vspace{0.1cm}
&$[4 ,15, 21- 3]$;\\

${\cal A}_{31}:$ &$[0 ,22 ,18 -21]$,&$[20 ,19, 9 -16]$,&$[7, 17, 12- 0]$,&$[16, 5 ,15 -10]$,&$[15 ,14, 4 -22]$,\\&$[ 3 ,13 ,14- 9]$,
&$[18 ,2, 13- 12]$,&$[1 ,11 ,18- 17]$,&$[0 ,13 ,21- 3]$,&$[12, 4, 16 -21]$,\\&$[3 ,15, 17- 1]$,&$[2 ,14, 16- 22]$,
&$[15, 1, 13 -9]$,&$[14 ,5 ,18- 12]$,&$[18, 6, 10 -2]$,\\&$[1 ,22 ,14 -21]$,&$[6, 19, 22- 5]$,&$[ 5 ,12, 21 -1]$,
&$[ 20, 17, 4 -10]$,&$[19 ,3 ,16- 11]$,\\&$[2 ,15, 12- 1]$,&$[7 ,18 ,15- 8]$,&$[3 ,12, 22- 15]$,&$[21, 17 ,9- 15]$,
&$[20, 1, 16 -0]$,\\&$[10 ,12 ,14- 11]$,&$[17 ,13, 10- 5]$,&$[8 ,18 ,16- 10]$,&$[1 ,19, 10- 3]$,&$[2 ,22 ,17 -0]$,\\\vspace{0.1cm}
&$[20, 3, 18- 9]$,&$[19 ,18 ,4 -21]$,&$[6 ,15 ,21- 2]$,&$[20 ,5 ,13 -22]$;\\

 ${\cal A}_{32}:$ &$[15 ,20, 9- 19]$,&$[7 ,17, 18- 1]$,&$[16 ,17 ,6 -14]$,&$[16 ,5 ,19- 8]$,&$[19, 14, 4- 10]$,\\
&$[3 ,13 ,14 -21]$,&$[12, 2, 13- 22]$,&$[0 ,13, 21 -1]$,&$[22 ,8, 15- 2]$,&$[21, 20, 7- 22]$,\\
 \end{tabular}
\end{center}

\begin{center}
\small
 \begin{tabular}{llllll}
&$[15 ,18, 6 -13]$,&$[5 ,17, 20- 14]$,
&$[3 ,19, 17- 0]$,&$[ 19 ,1, 13- 18]$,&$[1, 22 ,14 -11]$,\\&$[0 ,14 ,15- 7]$,&$[6 ,20, 22- 19]$,&$[5 ,18, 21- 2]$,
&$[15, 17, 4 -21]$,&$[20, 3, 16- 0]$,\\&$[2 ,19, 18- 0]$,&$[16, 13, 7- 14]$,&$[7 ,12, 19- 10]$,&$[19 ,0, 20- 13]$,
&$[3 ,18, 22- 5]$,\\&$[21 ,17, 9 -16]$,&$[15 ,1, 16- 10]$,&$[10, 18, 14- 9]$,&$[8 ,12, 16 -11]$,&$[19 ,15, 11- 3]$,\\\vspace{0.1cm}
&$[2 ,22, 17- 1]$,&$[2, 20, 11- 17]$,&$[20, 12 ,4- 22]$,&$[6 ,19, 21 -3]$,&$[15, 5 ,13- 9]$;\\

${\cal A}_{33}:$ &$[22, 21, 11- 8]$,&$[ 21 ,15 ,10- 8]$,&$[15 ,20, 9 -19]$,&$[16, 5 ,19- 8]$,&$[3 ,13 ,14 -21]$,\\
&$[0 ,10 ,11- 7]$,&$[22 ,9, 10- 7]$,&$[0 ,13 ,21- 1]$,&$[22, 8 ,15- 2]$,&$[21, 20, 7- 22]$,\\&$[5 ,17, 20- 14]$,&$[3 ,19, 17 -0]$,
&$[ 2 ,14 ,16- 21]$,&$[1 ,22, 14- 11]$,&$[0 ,14, 15- 7]$,\\&$[20, 3, 16- 0]$,&$[17 ,14 ,8 -21]$,&$[16 ,13, 7- 14]$,
&$[19 ,0, 20- 13]$,&$[21, 17, 9- 16]$,\\&$[15, 1, 16 -10]$,&$[17 ,13 ,10- 3]$,&$[19 ,15, 11- 3]$,&$[1 ,20, 10- 2]$,
&$[ 2, 22, 17 -1]$,\\&$[2 ,20, 11 -17]$,&$[15 ,5, 13- 9]$.\\
  \end{tabular}
\end{center}

\noindent Let $1\leq i\leq9$ and ${\cal B}_{i}={\cal A}_{i}\cup {\cal A}_{28}$. Let $i=10,11$
and ${\cal B}_{i}={\cal A}_{i}\cup {\cal A}_{29}$. Let $12\leq i\leq 16$ and ${\cal B}_{i}={\cal A}_{i}\cup {\cal A}_{30}$.
Let $17\leq i\leq 20$ and ${\cal B}_{i}={\cal A}_{i}\cup {\cal A}_{31}$. Let $21\leq i\leq 25$ and ${\cal B}_{i}={\cal A}_{i}\cup {\cal A}_{32}$.
Let $i=26,27$ and ${\cal B}_{i}={\cal A}_{i}\cup {\cal A}_{33}$.
Then $(X,{\cal B}_{j})$ is a $(K_{23}\setminus K_{10},G)$-design for each $1\leq j\leq 27$, where the removed subgraph $K_{10}$ is constructed on $Y=\{0,1,\ldots,9\}$. Now take permutation $\pi_{s,t}=(1)$ if $(s,t)\in N_{23}$. For $(s,t)\in M_{23}$, we take

\begin{center}
\small\tabcolsep 0.015in
\begin{tabular}{ll}
$\pi_{0,0}=(10\ 15\ 19\ 22)(11\ 13\ 18\ 17\ 21\ 20\ 12\ 16\ 14)$,&
$\pi_{0,1}=(10\ 12\ 13\ 17\ 18\ 22\ 20\ 21\ 19)(11\ 14\ 16)$,\\
$\pi_{0,2}=(10\ 14)(11\ 19\ 22\ 21\ 16\ 13\ 15\ 20\ 17\ 18\ 12)$,&
$\pi_{0,3}=(10\ 11)(12\ 22\ 16\ 13\ 19\ 15\ 14)(17\ 21\ 18\ 20)$,\\
$\pi_{0,4}=(10\ 15\ 16)(11\ 19\ 20\ 21\ 17\ 22\ 13\ 12\ 14\ 18)$,&
$\pi_{0,5}=(0\ 1)(2\ 3)(4\ 5)(6\ 7)(8\ 9)$,\\
$\pi_{0,6}=(10\ 17\ 18\ 16\ 12)(11\ 21\ 15\ 14\ 13\ 20)$,&
$\pi_{0,7}=(10\ 19\ 16\ 15\ 11\ 21\ 17)(12\ 22\ 20\ 13)(14\ 18)$,\\
$\pi_{0,8}=(10\ 14\ 12\ 20\ 22\ 15)(11\ 16\ 19\ 21)(17\ 18)$,&
$\pi_{0,9}=(10\ 21\ 14\ 18\ 22\ 20\ 19\ 11\ 13\ 17\ 12\ 15\ 16)$,\\
$\pi_{0,10}=(10\ 21\ 16\ 19\ 14\ 11\ 13\ 20)(12\ 18\ 15\ 22\ 17)$,&
$\pi_{0,11}=(10\ 22\ 14\ 17\ 16\ 18\ 13\ 15\ 21\ 20)(12\ 19)$,\\
$\pi_{0,12}=(11\ 21\ 20\ 12\ 16\ 13\ 14\ 19)(15\ 22)(17\ 18)$,&
$\pi_{0,13}=(10\ 20\ 13)(11\ 14)(12\ 18)(15\ 17\ 22\ 21)(16\ 19)$,\\
$\pi_{0,14}=(10\ 18\ 17\ 14)(12\ 22)(13\ 15)(16\ 20)(19\ 21)$,&
$\pi_{0,15}=(10\ 20\ 22\ 17\ 11\ 14\ 21)(12\ 18\ 13\ 15)(16\ 19)$,\\
$\pi_{0,16}=(10\ 21)(11\ 13\ 18\ 12\ 15)(14\ 20)$&
$\pi_{0,17}=(0\ 4\ 1\ 2\ 7)(3\ 8)(5\ 6\ 9)$\\
\hspace{1.2cm}$(16\ 19)(17\ 22)$,&\hspace{1.2cm}$(10\ 11\ 13\ 15\ 18\ 17\ 19\ 22\ 14\ 16\ 20\ 21)$,\\
$\pi_{0,18}=(0\ 3\ 9\ 1)(2\ 7\ 8)(4\ 6\ 5)(10\ 14\ 21\ 22)$&
$\pi_{7,0}=(0\ 7)(3\ 4)(6\ 9)(10\ 13\ 12)(15\ 18\ 16)$,\\
\hspace{1.2cm}$(11\ 13\ 18\ 15\ 12)(16\ 19\ 20)$,&
$\pi_{7,6}=(1\ 6\ 3\ 8\ 4)(10\ 18\ 22\ 15\ 11)$,\\
$\pi_{7,7}=(0\ 9\ 4\ 7)(10\ 13\ 20\ 12)(11\ 22)$,&
$\pi_{7,10}=(1\ 8\ 5\ 4\ 2\ 3)(12\ 14\ 18\ 17)$,\\
$\pi_{10,6}=(2\ 5)(4\ 8\ 6)(13\ 17)(16\ 21)$.
\end{tabular}
\end{center}

\noindent We have that for each $i,j,s,t$ in Table X,
 $\pi_{s,t}Y=Y$, $|{\cal B}_i\cap \pi_{s,t}{\cal B}_j|=s$ and
 $|T({\cal B}_i\setminus\pi_{s,t}{\cal B}_j)\cap T(\pi_{s,t}{\cal B}_j\setminus{\cal B}_i)|=t$. \qed

\begin{center}
\small  Table X.  \ Fine triangle intersections for $(K_{23}\setminus K_{10},G)$-designs

  \begin{tabular}{ccc|ccc|ccc|ccc}\hline
$i$ &  $j$ &  $(s,t)$ & $i$ &  $j$ &  $(s,t)$ & $i$ &  $j$ &  $(s,t)$ & $i$ &  $j$ &  $(s,t)$   \\\hline
1 & 1 &$M_{23}$ & 1&2 &   (49,3) & 1 & 3  & (46,6)&1 &  4   &   (43,9) \\
1 &    5   &   (40,12) &  1 &    6   &  (37,15) & 1 &7& (34,18) & 1 & 8 & (31,21)\\
1 &    9   &   (28,24) & 1  &   10   & (0,52) &1 & 12   & (42,0)&1 & 13   & (22,20)\\
1 &  14 & (37,5) & 1 &    15   &   (27,15) &1 &    16 & (32,10)&1 & 17   & (0,34) \\
1 & 18&   (15,19)&1 &    19   &   (10,24)&1 &    20   &   (5,29)  &1 & 21   & (26,0)\\
1 &  22 &   (14,12)&1 &  24   &   (18,8)&1 &    25   &   (21,5)&1 & 26   &   (22,0) \\
1 &    27   &   (14,0)& 10 &    2   & (3,49)&10  &   3  &  (6,46) &10  &   4  & (9,43) \\
10  &  5  & (12,40)&10  &   6  &   (15,37)&10  &   7  &   (18,34) &10  &   8  & (21,31)\\
10  &   9  &   (24,28) &  10 &11& (27,25) & 10  & 12  &   (0,42) &10  & 13  & (20,22)\\
10  &   14  &   (5,37) &10  &   15  &   (15,27)&10  &   16  &   (10,32)&10  &   17  &   (34,0)   \\
10  &   18  &   (19,15)&10  &   19  &   (24,10) &10  &   20  &   (29,5) & 10  &   21  &   (0,26)   \\
10  &   22  &   (12,14)& 10  &   23  &   (7,19) &10  &   25  &   (5,21)& 10  &   26  &   (0,22)\\
\hline
\end{tabular}
\end{center}

\end{document}